\documentclass[10pt]{article} 
\usepackage[utf8]{inputenc} 
\usepackage{graphicx} 
\usepackage{booktabs} 
\usepackage{array} 
\usepackage{enumerate}
\usepackage{amsmath}
\usepackage{graphicx}
\usepackage{float}
\usepackage{subfig}
\usepackage{tensor}
\usepackage{tikz}
\usepackage{lmodern}
\usetikzlibrary{arrows,calc,decorations.markings}

\newcommand\myline[1][]{%
  \,\tikz[baseline]\draw[thin,#1](0,-\dp\strutbox)--(0,\ht\strutbox);\,%
}

\makeatletter
\newcommand\resetsubfigs{\setcounter{sub\@captype}{0}}
\makeatother

\begin{document}

\pgfarrowsdeclarecombine{dimarrow}{dimarrow}{latex}{latex}{}{}
\def\Dimline[#1][#2][#3][#4]{
    \begin{scope}[>=latex] 
        \draw[|-|,
        decoration={markings, 
  		       mark=at position 0 with {\arrowreversed[scale=1]{dimarrow}};,
                mark=at position .5 with {\node[#4] at (0,-0.35) {\tiny{#3}};},
			  mark=at position 1 with {\arrow[scale=1]{dimarrow}};,
        },
        postaction=decorate] #1 -- #2 ;
    \end{scope}
}
\allowdisplaybreaks
\allowdisplaybreaks

\title{A New Method for Signal and Image Analysis: The Square Wave Method}
\author{
Osvaldo Skliar\thanks{oskliar@costarricense.cr; Escuela de Inform\'atica, Universidad Nacional, Costa Rica.}
\and Ricardo E. Monge\thanks{ricardo@mogap.net; Escuela de Ciencias de la Computaci\'on e Inform\'atica, Universidad de Costa Rica, Costa Rica.} \and Sherry Gapper\thanks{sherry.gapper.morrow@una.cr; Universidad Nacional, Costa Rica.}}
\maketitle
\begin{abstract}
A brief review is provided of the use of the Square Wave Method (SWM) in the field of signal and image analysis and it is specified how results thus obtained are expressed using the Square Wave Transform (SWT), in the frequency domain. To illustrate the new approach introduced in this field, the results of two cases are analyzed: a) a sequence of samples (that is, measured values) of an electromyographic recording; and b) the classic image of Lenna.
\end{abstract}

Mathematics Subject Classification: 94A12, 65F99

\textbf{Keywords:} signal and image analysis, Square Wave Method (SWM), Square Wave Transform (SWT).

\section{Introduction}

It was previously shown how a new method, the Square Wave Method (SWM), for the analysis of signals depending on one variable \cite{b1} can be presented in the frequency domain by using a mathematical tool called Square Wave Transform (SWT) \cite{b2} \cite{b3}. The SWM was then generalized quite naturally and directly for image analysis \cite{b4}.

The objectives of this paper are the following:
\begin{enumerate}
\item To provide a brief review of the use of the SWM for the analysis of signals and specify the relations existing between a) the sampling frequency $f_s$, with which the successive values of recordings of biomedical signals (such as those of an electrocardiogram, electromyogram or electroencephalogram) are measured, and b) the frequencies $f_1,f_2, \ldots,f_n$, corresponding respectively to the different trains of square waves $S_1, S_2, \ldots, S_n$ obtained using the SWM;
\item To indicate how it also is possible to present in the frequency domain using the SWT, the results of the analysis of images obtained with the SWM.
\end{enumerate}

The application of the SWM in the field of signal and image analysis is exemplified with the results of an analysis of a) a sequence of samples (that is, measured values from an electromyographic recording); and b) the classic image of Lenna, using the SWT \cite{b5}.

\section{Analysis of a Function of One Variable}

Consider a function of time ($t$), in the interval $\Delta t$, satisfying the conditions of Dirichlet \cite{b6}:
\begin{equation}
f(t)=(6-t)(2\cos(2 \pi 4 t)+5\cos(2\pi 6 t)) \quad 0\leq t\leq 4 \: \mathrm{s}
\label{e1}
\end{equation}

Suppose that the time interval in which the function characterized in equation~\eqref{e1} will be analyzed ($\Delta t=4\:\mathrm{s}$) has been divided into 18 equal sub-intervals. In this case, it will be seen that function (1) can be approximated in $\Delta t$, using the sum of the parts corresponding to $\Delta t$ of 18 trains of square waves. These trains of square waves will be called $S_1,S_2,S_3,\ldots,S_{18}$; the ``$S$'' being based on the word ``square'' in the expression ``train of square waves''.

If $\Delta t$ has been divided into 100 equal sub-intervals, the approximation to the function~\eqref{e1} in interval $\Delta t$  will be carried out by adding the parts corresponding to $\Delta t$ of 100 trains of square waves: $S_1,S_2,S_3,\ldots,S_{100}$. In general, if $\Delta t$ is divided into any natural number $n$ of equal sub-intervals, the approximation in $\Delta t$ to function~\eqref{e1} will be obtained by adding the parts corresponding to $\Delta t$ of $n$ trains of square waves: $S_1,S_2,S_3,\ldots,S_{n}$. The Square Wave Method (SWM) described in this section makes it possible to determine those  trains of square waves unambiguously.  Therefore, each $S_i$ (where $i=1,2,\ldots,n$) of those trains of square waves will be characterized by a specific frequency $f_i$ (i.e., consideration is given to the number of waves in the train of square waves, which is contained in the unit of time 1 s) and a particular coefficient $C_i$, whose absolute value is the amplitude of the corresponding train.

\begin{figure}[H]
\centering
\includegraphics[width=5in]{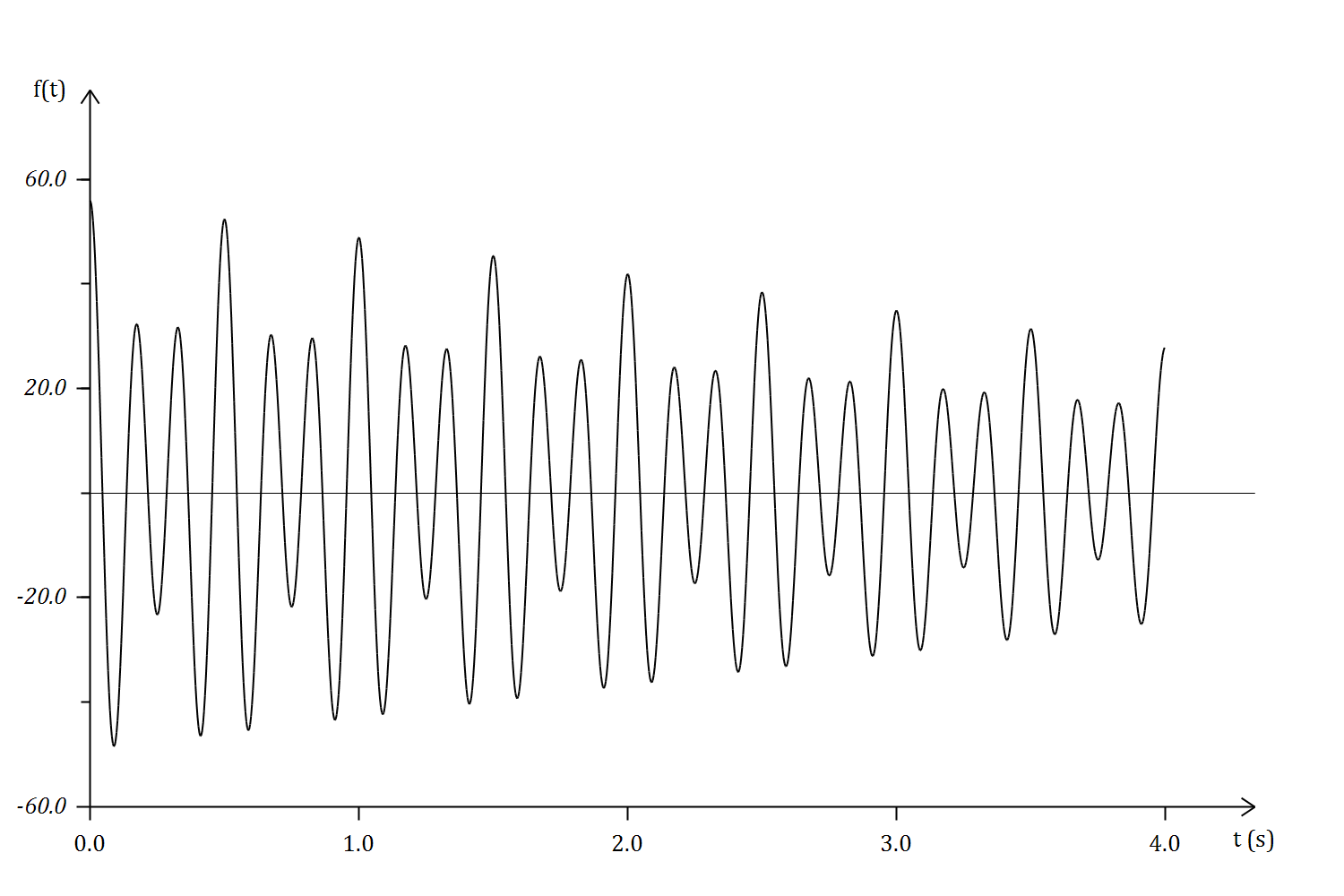}
\caption{$f(t)=(6-t)(2\cos (2 \pi 4 t)+5\cos (2\pi 6 t)) \quad 0\leq t\leq 4 \: \mathrm{s}$.}
\label{f1}
\end{figure}

The function $f(t)$ specified in~\eqref{e1} is shown in figure~\ref{f1}.

\begin{figure}[H]
\centering
\includegraphics[width=4.8in]{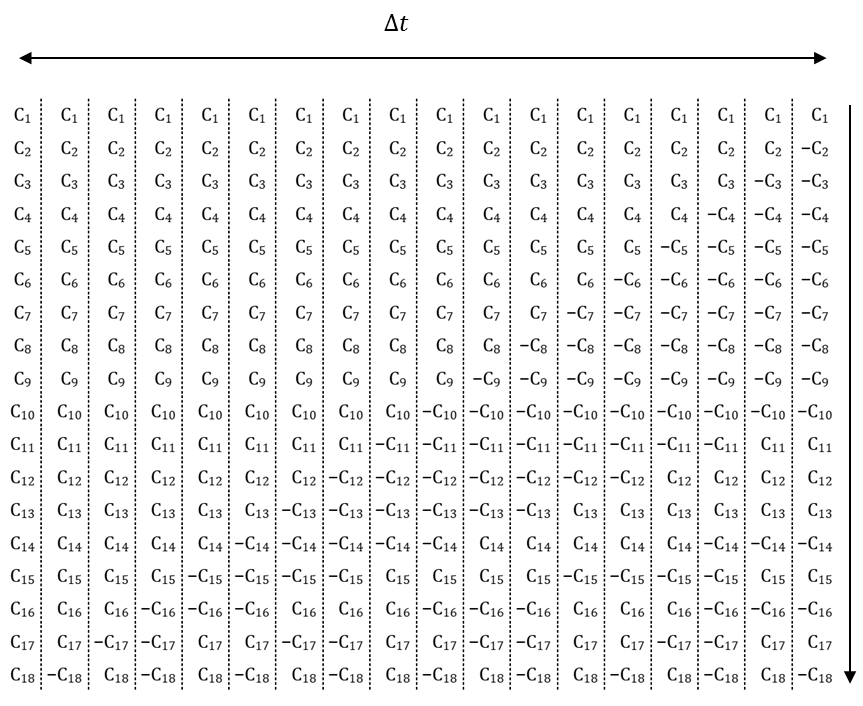}
\caption{How to apply the SWM to the analysis of the function represented in figure~\ref{f1}. (See indications in text.)}
\label{f2}
\end{figure}

For the case considered here, $n=18$, a description will be provided below of how the frequencies $f_i$ (where $i=1,2,\ldots,18$) and the values of the coefficients $C_i$ (where $i=1,2,\ldots,18$) corresponding to the different trains of square waves $S_i$ (where $i=1,2,\ldots,18$) are determined; see figure~\ref{f2}.

The first row of figure~\ref{f2} (with coefficients $C_1$) represents half a square wave, the first semi-wave of the train of square waves $S_1$. The frequency of $S_1$ (i.e., $f_1$) is clearly equal to the number of square waves per unit of time (1 s). To obtain $f_1$, the part of $S_1$ which occupies $\Delta t$ (the half-wave) is divided by $\Delta t$.

\begin{equation*}
f_1=\frac{\tfrac{1}{2}}{\Delta t}=\frac{\tfrac{1}{2}}{4\:\mathrm{s}}=\frac{1}{8}\:\mathrm{s}^{-1}
\end{equation*}

To compute $f_2$, note that $\Delta t$ is occupied by the sum of that half-wave of the train of square waves $S_2$ and the fraction $\tfrac{1}{17}$ of the second semi-wave of the first square wave of $S_2$. That fraction is represented by the symbol $C_2$ in the second row of figure~\ref{f2}. Thus the following value is obtained for $f_2$:
\begin{equation*}
f_2=\frac{\tfrac{1}{2} + \left( \tfrac{1}{17}\cdot\tfrac{1}{2}\right)}{\Delta t} 
= \frac{\tfrac{1}{2}\left( 1+\tfrac{1}{17}\right)}{\Delta t}
=\frac{1}{2}\cdot\frac{\tfrac{18}{17}}{\Delta t}
=\frac{1}{2\Delta t}\left(\frac{18}{18-1}\right)
=\frac{1}{8}\left(\frac{18}{18-1}\right)\:\mathrm{s}^{-1}
\end{equation*}

To compute $f_3$, note that $\Delta t$ is occupied by the sum of that half-wave of the train of square waves $S_3$ and the fraction $\tfrac{2}{16}$ of the second semi-wave of the first square wave of $S_3$. This fraction (in the third row of figure~\ref{f2}) is represented by the sequence of symbols $-C_3$\myline[densely dashed]$-C_3$. Therefore, the following value is obtained for $f_3$:

\begin{equation*}
f_3=\frac{\tfrac{1}{2} + \left( \tfrac{2}{16}\cdot\tfrac{1}{2}\right)}{\Delta t} 
= \frac{\tfrac{1}{2}\left( 1+\tfrac{2}{16}\right)}{\Delta t}
=\frac{1}{2\Delta t}\left(\frac{18}{18-2}\right)
=\frac{1}{8}\left(\frac{18}{18-2}\right)\:\mathrm{s}^{-1}
\end{equation*}

With a precision of 7 decimal places, the values are given below not only for $f_1$, $f_2$,$f_3$, but also for those corresponding to $f_4,f_5,\ldots,f_{17},f_{18}$.

\begin{equation*}
f_{1}=\frac{1}{2\Delta t}\left(\frac{18}{18-0}\right)=\frac{1}{8}\left(\frac{18}{18-0}\right)\:\mathrm{s}^{-1}=0.1250000\:\mathrm{s}^{-1}
\end{equation*}
\begin{equation*}
f_{2}=\frac{1}{2\Delta t}\left(\frac{18}{18-1}\right)=\frac{1}{8}\left(\frac{18}{18-1}\right)\:\mathrm{s}^{-1}=0.1323529\:\mathrm{s}^{-1}
\end{equation*}
\begin{equation*}
f_{3}=\frac{1}{2\Delta t}\left(\frac{18}{18-2}\right)=\frac{1}{8}\left(\frac{18}{18-2}\right)\:\mathrm{s}^{-1}=0.1406250\:\mathrm{s}^{-1}
\end{equation*}
\begin{equation*}
f_{4}=\frac{1}{2\Delta t}\left(\frac{18}{18-3}\right)=\frac{1}{8}\left(\frac{18}{18-3}\right)\:\mathrm{s}^{-1}=0.1500000\:\mathrm{s}^{-1}
\end{equation*}
\begin{equation*}
f_{5}=\frac{1}{2\Delta t}\left(\frac{18}{18-4}\right)=\frac{1}{8}\left(\frac{18}{18-4}\right)\:\mathrm{s}^{-1}=0.1607143\:\mathrm{s}^{-1}
\end{equation*}
\begin{equation*}
f_{6}=\frac{1}{2\Delta t}\left(\frac{18}{18-5}\right)=\frac{1}{8}\left(\frac{18}{18-5}\right)\:\mathrm{s}^{-1}=0.1730769\:\mathrm{s}^{-1}
\end{equation*}
\begin{equation*}
f_{7}=\frac{1}{2\Delta t}\left(\frac{18}{18-6}\right)=\frac{1}{8}\left(\frac{18}{18-6}\right)\:\mathrm{s}^{-1}=0.1875000\:\mathrm{s}^{-1}
\end{equation*}
\begin{equation*}
f_{8}=\frac{1}{2\Delta t}\left(\frac{18}{18-7}\right)=\frac{1}{8}\left(\frac{18}{18-7}\right)\:\mathrm{s}^{-1}=0.2045455\:\mathrm{s}^{-1}
\end{equation*}
\begin{equation*}
f_{9}=\frac{1}{2\Delta t}\left(\frac{18}{18-8}\right)=\frac{1}{8}\left(\frac{18}{18-8}\right)\:\mathrm{s}^{-1}=0.2250000\:\mathrm{s}^{-1}
\end{equation*}
\begin{equation*}
f_{10}=\frac{1}{2\Delta t}\left(\frac{18}{18-9}\right)=\frac{1}{8}\left(\frac{18}{18-9}\right)\:\mathrm{s}^{-1}=0.2500000\:\mathrm{s}^{-1}
\end{equation*}
\begin{equation*}
f_{11}=\frac{1}{2\Delta t}\left(\frac{18}{18-10}\right)=\frac{1}{8}\left(\frac{18}{18-10}\right)\:\mathrm{s}^{-1}=0.2812500\:\mathrm{s}^{-1}
\end{equation*}
\begin{equation*}
f_{12}=\frac{1}{2\Delta t}\left(\frac{18}{18-11}\right)=\frac{1}{8}\left(\frac{18}{18-11}\right)\:\mathrm{s}^{-1}=0.3214286\:\mathrm{s}^{-1}
\end{equation*}
\begin{equation*}
f_{13}=\frac{1}{2\Delta t}\left(\frac{18}{18-12}\right)=\frac{1}{8}\left(\frac{18}{18-12}\right)\:\mathrm{s}^{-1}=0.3750000\:\mathrm{s}^{-1}
\end{equation*}
\begin{equation*}
f_{14}=\frac{1}{2\Delta t}\left(\frac{18}{18-13}\right)=\frac{1}{8}\left(\frac{18}{18-13}\right)\:\mathrm{s}^{-1}=0.4500000\:\mathrm{s}^{-1}
\end{equation*}
\begin{equation*}
f_{15}=\frac{1}{2\Delta t}\left(\frac{18}{18-14}\right)=\frac{1}{8}\left(\frac{18}{18-14}\right)\:\mathrm{s}^{-1}=0.5625000\:\mathrm{s}^{-1}
\end{equation*}
\begin{equation*}
f_{16}=\frac{1}{2\Delta t}\left(\frac{18}{18-15}\right)=\frac{1}{8}\left(\frac{18}{18-15}\right)\:\mathrm{s}^{-1}=0.7500000\:\mathrm{s}^{-1}
\end{equation*}
\begin{equation*}
f_{17}=\frac{1}{2\Delta t}\left(\frac{18}{18-16}\right)=\frac{1}{8}\left(\frac{18}{18-16}\right)\:\mathrm{s}^{-1}=1.1250000\:\mathrm{s}^{-1}
\end{equation*}
\begin{equation*}
f_{18}=\frac{1}{2\Delta t}\left(\frac{18}{18-17}\right)=\frac{1}{8}\left(\frac{18}{18-17}\right)\:\mathrm{s}^{-1}=2.2500000\:\mathrm{s}^{-1}
\end{equation*}

Observe that any of the 18 values of $f_i$, where $i=1,2,\ldots,18$, can be computed with the following equation:
\begin{equation*}
f_{i}=\frac{1}{2\Delta t}\left(\frac{18}{18-(i-1)}\right)=\frac{1}{8}\left(\frac{18}{18-(i-1)}\right)\:\mathrm{s}^{-1};\quad i=1,2,\ldots,18 
\end{equation*}

In general, if the interval $\Delta t$, whose value, of course, may be different from 4 s, is divided into $n$ equal sub-intervals, the frequencies corresponding to each of the $n$ trains of square waves are as follows:
\begin{equation}
f_{i}=\frac{1}{2\Delta t}\left(\frac{n}{n-(i-1)}\right)\:\mathrm{s}^{-1};\quad i=1,2,\ldots,n 
\label{e2}
\end{equation}

It has been explained how to compute each $f_i$ corresponding to each $S_i$, where $i=1,2,\ldots,18$, for the case of the approximation to $f(t)$ specified in~\eqref{e1} when dividing $\Delta t$ into 18 equal sub-intervals ($n=18$). Indications will now be given on how to compute the $C_i$ for each $S_i$.

The vertical arrow pointing down at the right of figure 2 indicates how to add the terms corresponding to each of the 18 sub-intervals of $\Delta t$. Thus, to obtain the values of the coefficients $C_1,C_2,\ldots$, $C_{17}$ and $C_{18}$, corresponding to $S_1,S_2,\ldots$, $S_{17}$ and $S_{18}$, the following system of linear equations must be solved.

{
 \setlength{\abovedisplayskip}{6pt}
  \setlength{\belowdisplayskip}{\abovedisplayskip}
  \setlength{\abovedisplayshortskip}{0pt}
  \setlength{\belowdisplayshortskip}{3pt}
\begin{equation}
  \left.\begin{aligned}
C_{1} & +C_{2}+C_{3}+C_{4}+C_{5}+C_{6}+C_{7}+C_{8}+C_{9}+C_{10}& \\&+C_{11}+C_{12}+C_{13}+C_{14}+C_{15}+C_{16}+C_{17}+C_{18} &= V_{1} \\
C_{1}&+C_{2}+C_{3}+C_{4}+C_{5}+C_{6}+C_{7}+C_{8}+C_{9}+C_{10}& \\&+C_{11}+C_{12} +C_{13}+C_{14}+C_{15}+C_{16}+C_{17}-C_{18} &= V_{2} \\
C_{1}&+C_{2}+C_{3}+C_{4}+C_{5}+C_{6}+C_{7}+C_{8}+C_{9}+C_{10}& \\&+C_{11}+C_{12}+C_{13}+C_{14}+C_{15}+C_{16}-C_{17}-C_{18} &= V_{3} \\
C_{1}&+C_{2}+C_{3}+C_{4}+C_{5}+C_{6}+C_{7}+C_{8}+C_{9}+C_{10}& \\&+C_{11}+C_{12} +C_{13}+C_{14}+C_{15}-C_{16}-C_{17}-C_{18} &= V_{4} \\
C_{1}&+C_{2}+C_{3}+C_{4}+C_{5}+C_{6}+C_{7}+C_{8}+C_{9}+C_{10}& \\&+C_{11}+C_{12} +C_{13}+C_{14}-C_{15}-C_{16}-C_{17}-C_{18} &= V_{5} \\
C_{1}&+C_{2}+C_{3}+C_{4}+C_{5}+C_{6}+C_{7}+C_{8}+C_{9}+C_{10}& \\&+C_{11}+C_{12} +C_{13}-C_{14}-C_{15}-C_{16}-C_{17}-C_{18} &= V_{6} \\
C_{1}&+C_{2}+C_{3}+C_{4}+C_{5}+C_{6}+C_{7}+C_{8}+C_{9}+C_{10}& \\&+C_{11}+C_{12} -C_{13}-C_{14}-C_{15}-C_{16}-C_{17}-C_{18} &= V_{7} \\
C_{1}&+C_{2}+C_{3}+C_{4}+C_{5}+C_{6}+C_{7}+C_{8}+C_{9}+C_{10}& \\&+C_{11}-C_{12} -C_{13}-C_{14}-C_{15}-C_{16}-C_{17}-C_{18} &= V_{8} \\
C_{1}&+C_{2}+C_{3}+C_{4}+C_{5}+C_{6}+C_{7}+C_{8}+C_{9}+C_{10}& \\&-C_{11}-C_{12} -C_{13}-C_{14}-C_{15}-C_{16}-C_{17}-C_{18} &= V_{9} \\
C_{1}&+C_{2}+C_{3}+C_{4}+C_{5}+C_{6}+C_{7}+C_{8}+C_{9}-C_{10}& \\&-C_{11}-C_{12} -C_{13}-C_{14}-C_{15}-C_{16}-C_{17}-C_{18} &= V_{10} \\
C_{1}&+C_{2}+C_{3}+C_{4}+C_{5}+C_{6}+C_{7}+C_{8}-C_{9}-C_{10}& \\&-C_{11}-C_{12}-C_{13}-C_{14}-C_{15}-C_{16}+C_{17}+C_{18} &= V_{11} \\
C_{1}&+C_{2}+C_{3}+C_{4}+C_{5}+C_{6}+C_{7}-C_{8}-C_{9}-C_{10}& \\&-C_{11}-C_{12}-C_{13}-C_{14}+C_{15}+C_{16}+C_{17}+C_{18} &= V_{12} \\
C_{1}&+C_{2}+C_{3}+C_{4}+C_{5}+C_{6}-C_{7}-C_{8}-C_{9}-C_{10}& \\&-C_{11}-C_{12} +C_{13}+C_{14}+C_{15}+C_{16}+C_{17}+C_{18} &= V_{13} \\
C_{1}&+C_{2}+C_{3}+C_{4}+C_{5}-C_{6}-C_{7}-C_{8}-C_{9}-C_{10}& \\&+C_{11}+C_{12} +C_{13}+C_{14}+C_{15}-C_{16}-C_{17}-C_{18} &= V_{14} \\
C_{1}&+C_{2}+C_{3}+C_{4}-C_{5}-C_{6}-C_{7}-C_{8}+C_{9}+C_{10}& \\&+C_{11}+C_{12} -C_{13}-C_{14}-C_{15}-C_{16}+C_{17}+C_{18} &= V_{15} \\
C_{1}&+C_{2}+C_{3}-C_{4}-C_{5}-C_{6}+C_{7}+C_{8}+C_{9}-C_{10}& \\&-C_{11}-C_{12} +C_{13}+C_{14}+C_{15}-C_{16}-C_{17}-C_{18} &= V_{16} \\
C_{1}&+C_{2}-C_{3}-C_{4}+C_{5}+C_{6}-C_{7}-C_{8}+C_{9}+C_{10}& \\&-C_{11}-C_{12}+C_{13}+C_{14}-C_{15}-C_{16}+C_{17}+C_{18} &= V_{17} \\
C_{1}&-C_{2}+C_{3}-C_{4}+C_{5}-C_{6}+C_{7}-C_{8}+C_{9}-C_{10}& \\&+C_{11}-C_{12} +C_{13}-C_{14}+C_{15}-C_{16}+C_{17}-C_{18} &= V_{18} 
   \end{aligned}
  \right\}  
\label{e3}
\end{equation}
}

In the preceding system of linear algebraic equations \eqref{e3}, $V_1,V_2,\ldots, V_{17}$ and $V_{18}$  are the values for $f(t)$ as specified in~\eqref{e1} at the midpoints of the first, second, third, $\ldots$, seventeenth and eighteenth sub-intervals, respectively, of interval $\Delta t$, in which $f(t)$ is analyzed. It follows that the values $V_i$ (where $i=1,2,3,\ldots,17$ and $18$) can be computed given that $f(t)$ has been specified in \eqref{e1}. These values are as follows:
\begin{align*}
& V_{1}=-34.5484836 & & V_{10}=-25.7897131\\
& V_{2}=30.6666667 & & V_{11}=22.6666667\\
& V_{3}=-16.0256827 & & V_{12}=-11.7202754\\
& V_{4}=-6.9904692 & & V_{13}=-5.0546469\\
& V_{5}=49.0000000 & & V_{14}=35.0000000\\
& V_{6}=-6.5602864 & & V_{15}=-4.6244642\\
& V_{7}=-14.1121683 & & V_{16}=-9.8067610\\
& V_{8}=25.3333333 & & V_{17}=17.3333333\\
& V_{9}=-26.7629098 & & V_{18}=-18.0041393
\end{align*}

Each of the 18 values of $V_i$, where $i=1,2,\ldots,18$, has been computed with a precision of seven decimal digits.

The 18 unknowns of the systems of equations specified in \eqref{e3} are $C_1,C_2,\ldots$, $C_{17}$, and $C_{18}$. Thus $|C_i|$ refers to the amplitude of the train of square waves $S_i$, where $i=1,2,\ldots,18$. The (constant) value of each positive square semi-wave of the train of square waves $S_i$ is $|C_i|$ and the (constant) value of each negative square semi-wave of that $S_i$ is $-|C_i|$.

The system of equations \eqref{e3} has been solved by using LAPACK \cite{b7}, and the following results were obtained for the unknowns:
\begin{align*}
& C_{1}=117.12980 & & C_{10}=4.03101\\
& C_{2}=50.27631 & & C_{11}=-85.68506\\
& C_{3}=-210.98830 & & C_{12}=12.88482\\
& C_{4}=-53.27896 & & C_{13}=8.51973\\
& C_{5}=9.35088 & & C_{14}=60.38772\\
& C_{6}=12.58025 & & C_{15}=-69.86421\\
& C_{7}=61.27212 & & C_{16}=28.08997\\
& C_{8}=49.80105 & & C_{17}=-9.26140\\
& C_{9}=12.81335 & & C_{18}=-32.60758
\end{align*}

The trains of square waves $S_1,S_2,S_3,\ldots,S_{17}$ and $S_{18}$ have been shown for interval $\Delta t$ in figures \ref{f3a}, \ref{f3b}, \ref{f3c}, \ldots, \ref{f3r}, respectively.
\begin{figure}
\centering
\subfloat[: $S_{1}(t)$]{\includegraphics[width=4.8in]{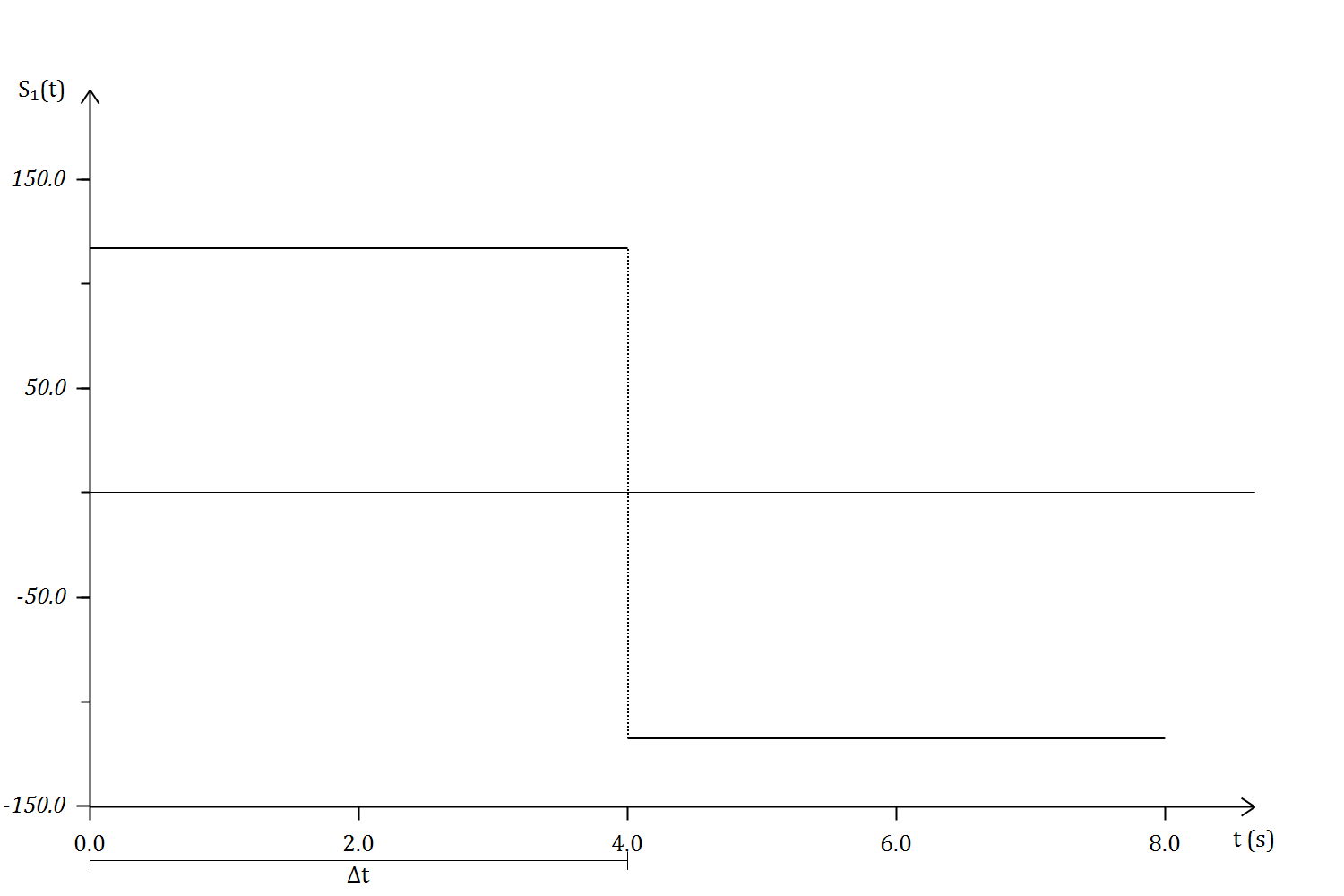}\label{f3a}}\qquad
\subfloat[: $S_{2}(t)$]{\includegraphics[width=4.8in]{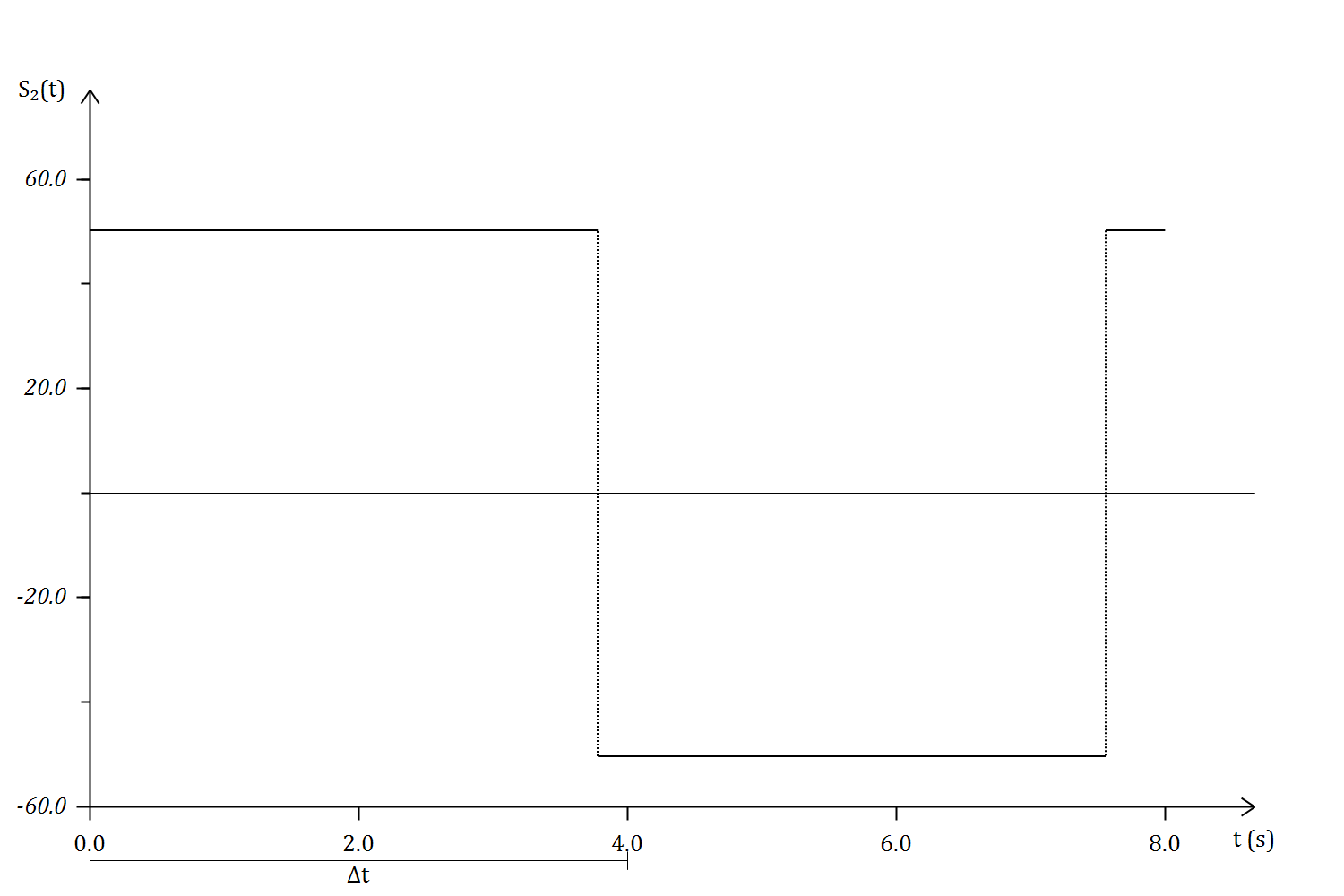}\label{f3b}}\qquad
\caption[]{}\end{figure}\begin{figure}\ContinuedFloat
\centering
\subfloat[: $S_{3}(t)$]{\includegraphics[width=4.8in]{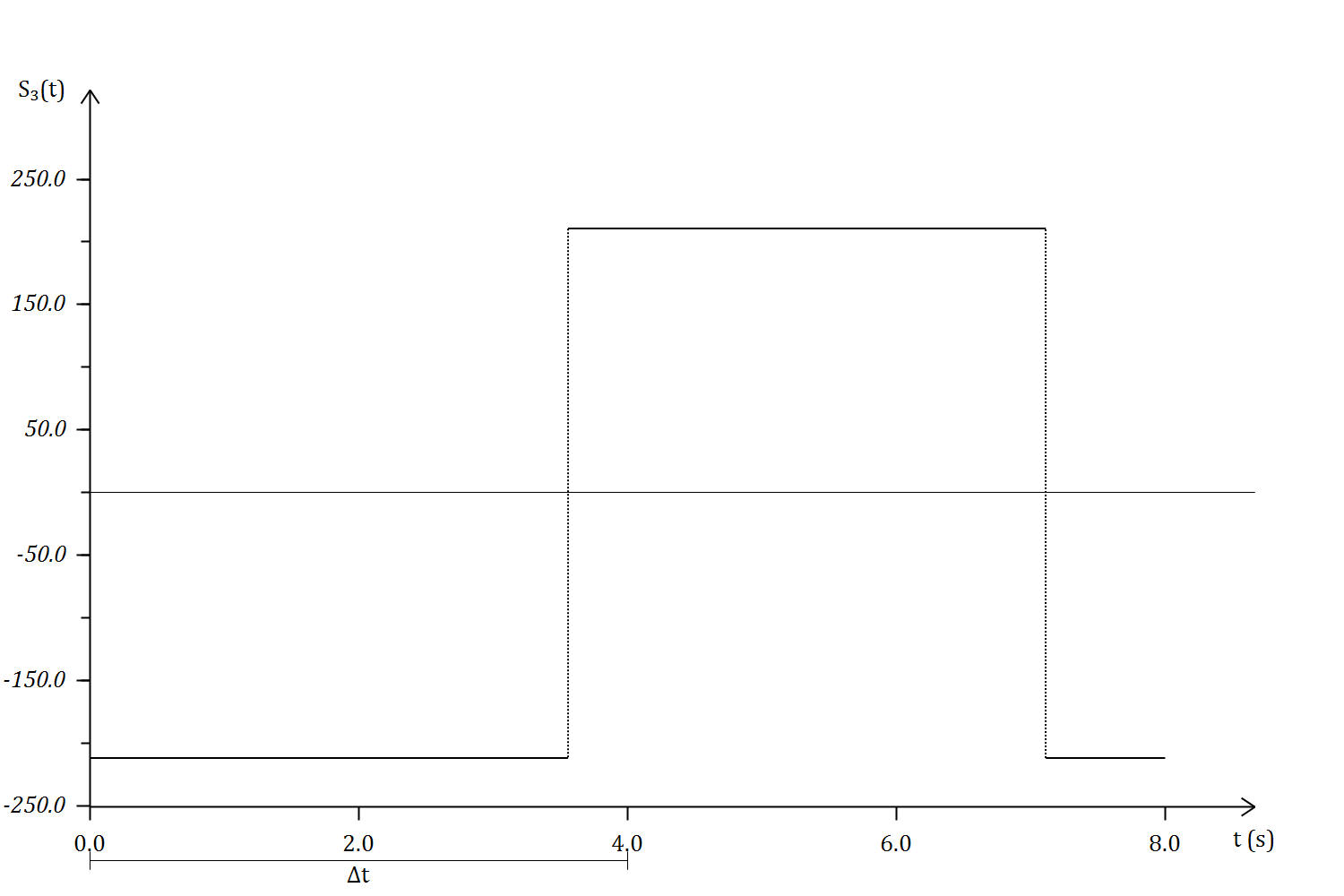}\label{f3c}}\qquad
\subfloat[: $S_{4}(t)$]{\includegraphics[width=4.8in]{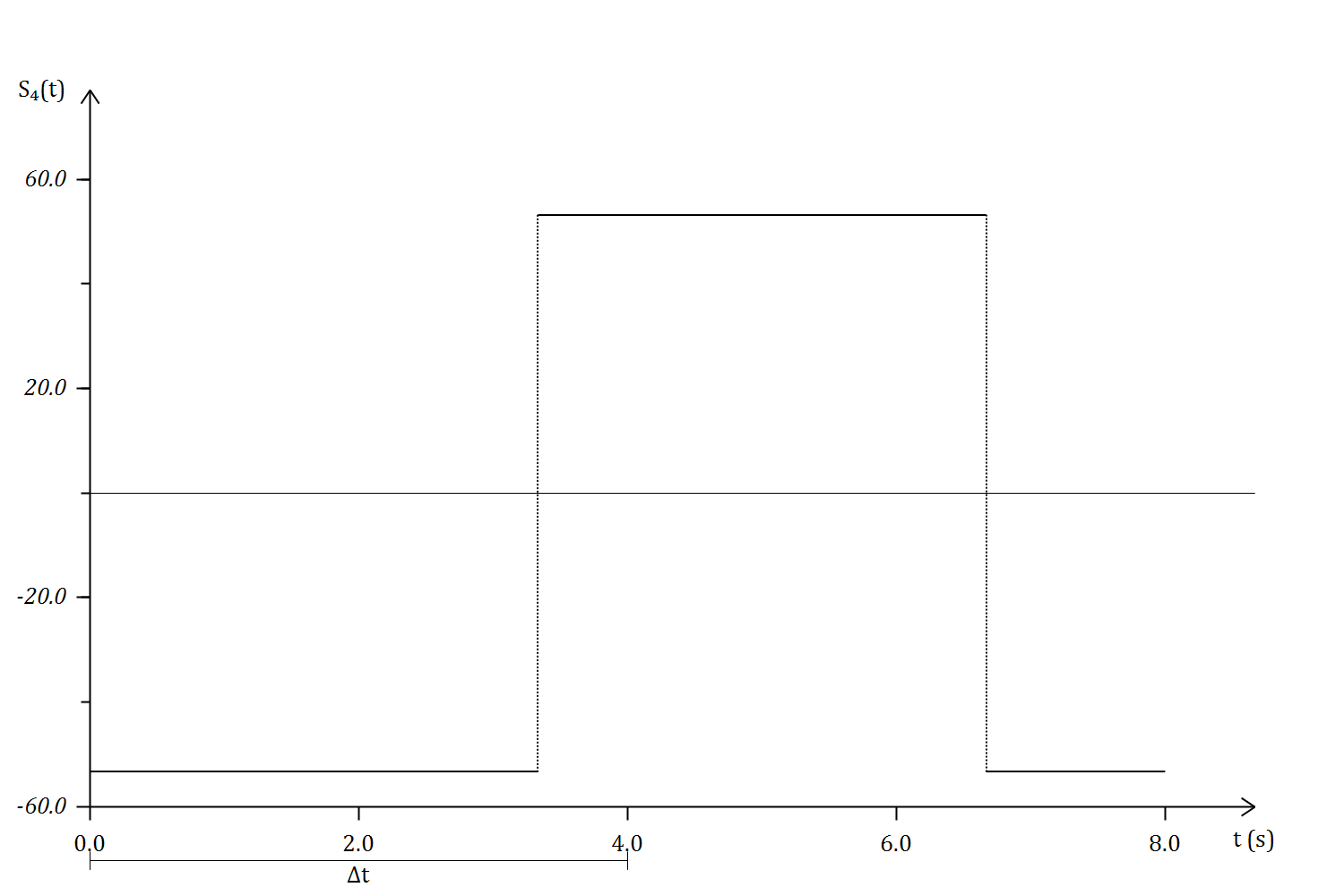}\label{f3d}}\qquad
\caption[]{}\end{figure}\begin{figure}\ContinuedFloat
\centering
\subfloat[: $S_{5}(t)$]{\includegraphics[width=4.8in]{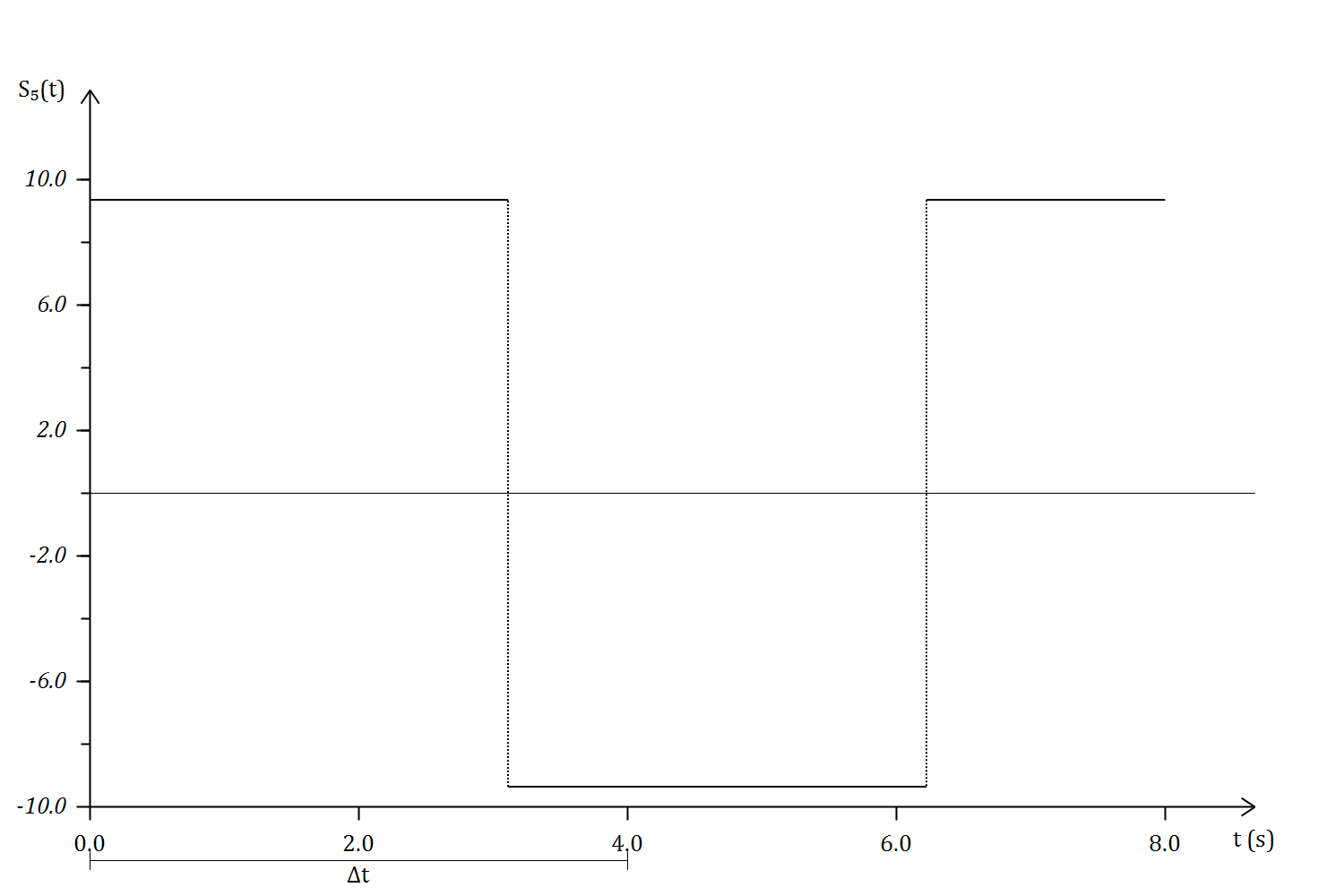}\label{f3e}}\qquad
\subfloat[: $S_{6}(t)$]{\includegraphics[width=4.8in]{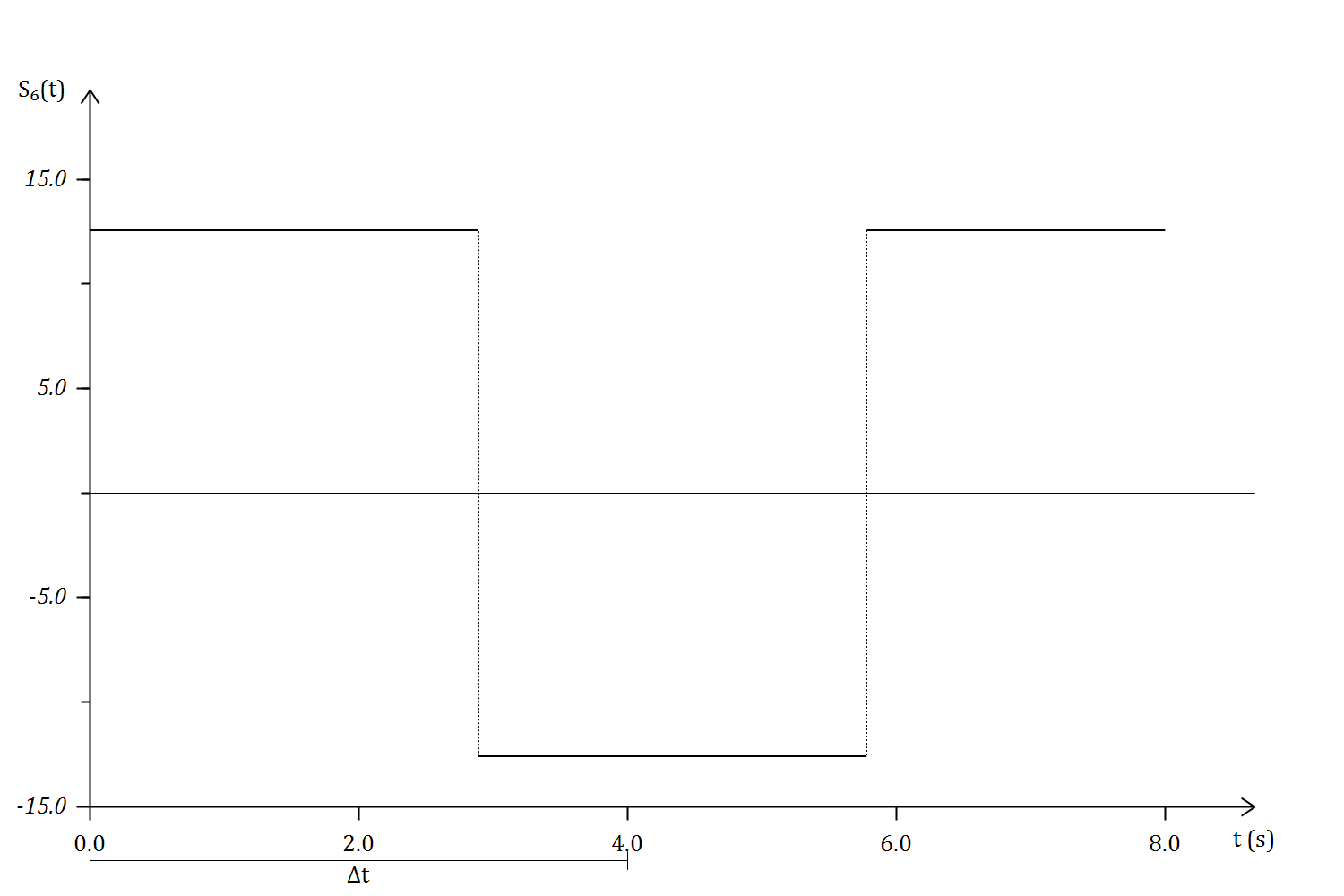}\label{f3f}}\qquad
\caption[]{}\end{figure}\begin{figure}\ContinuedFloat
\centering
\subfloat[: $S_{7}(t)$]{\includegraphics[width=4.8in]{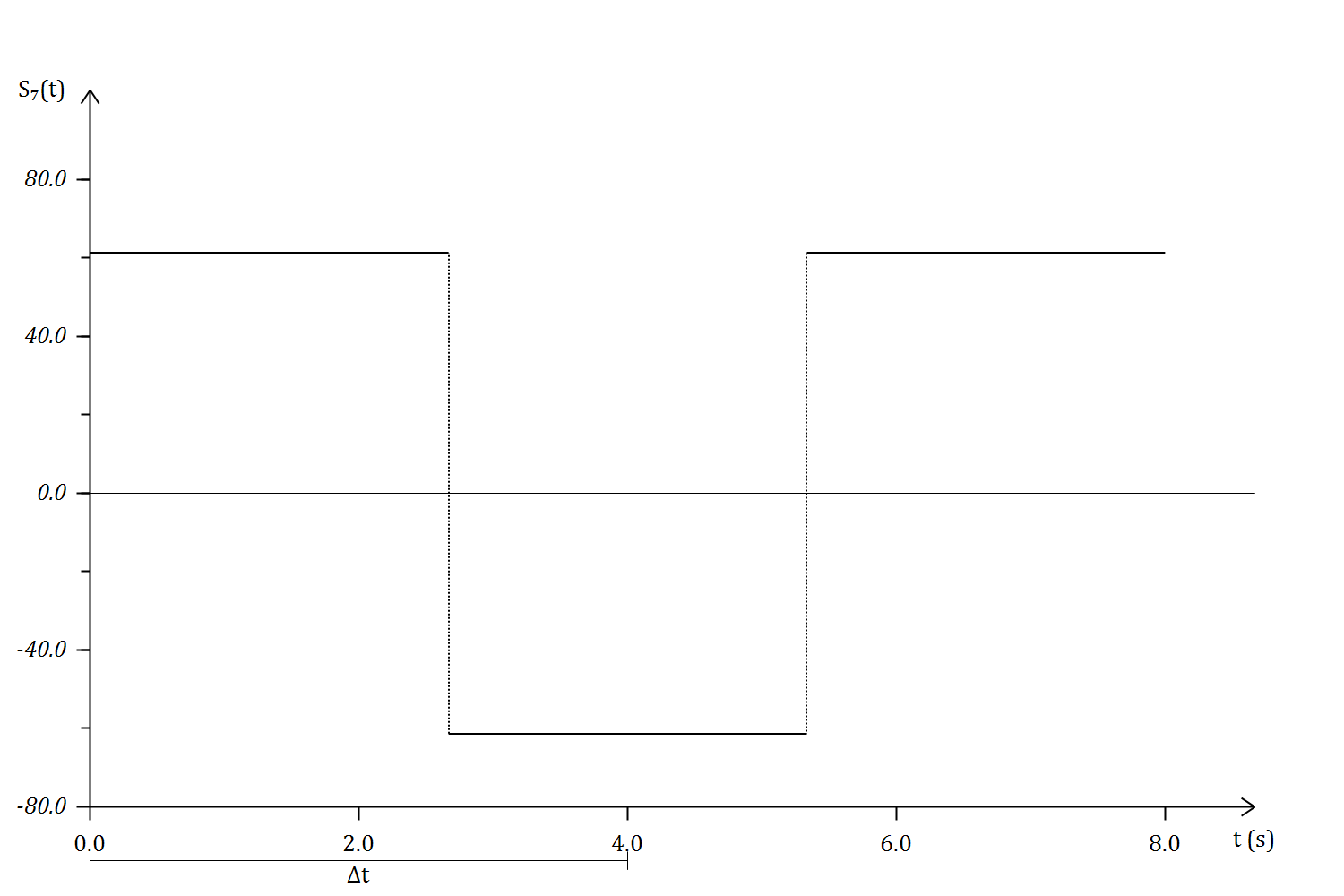}\label{f3g}}\qquad
\subfloat[: $S_{8}(t)$]{\includegraphics[width=4.8in]{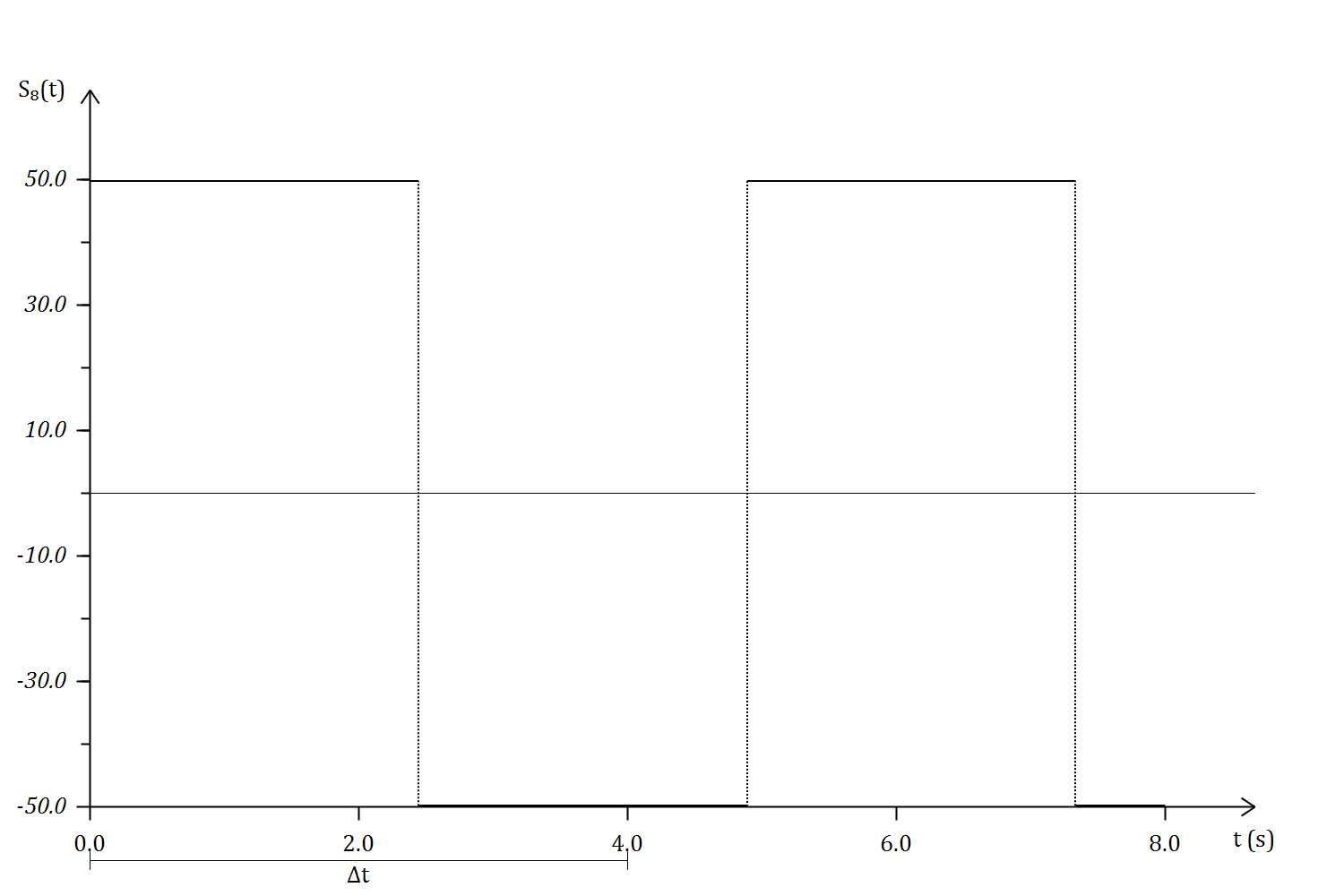}\label{f3h}}\qquad
\caption[]{}\end{figure}\begin{figure}\ContinuedFloat
\centering
\subfloat[: $S_{9}(t)$]{\includegraphics[width=4.8in]{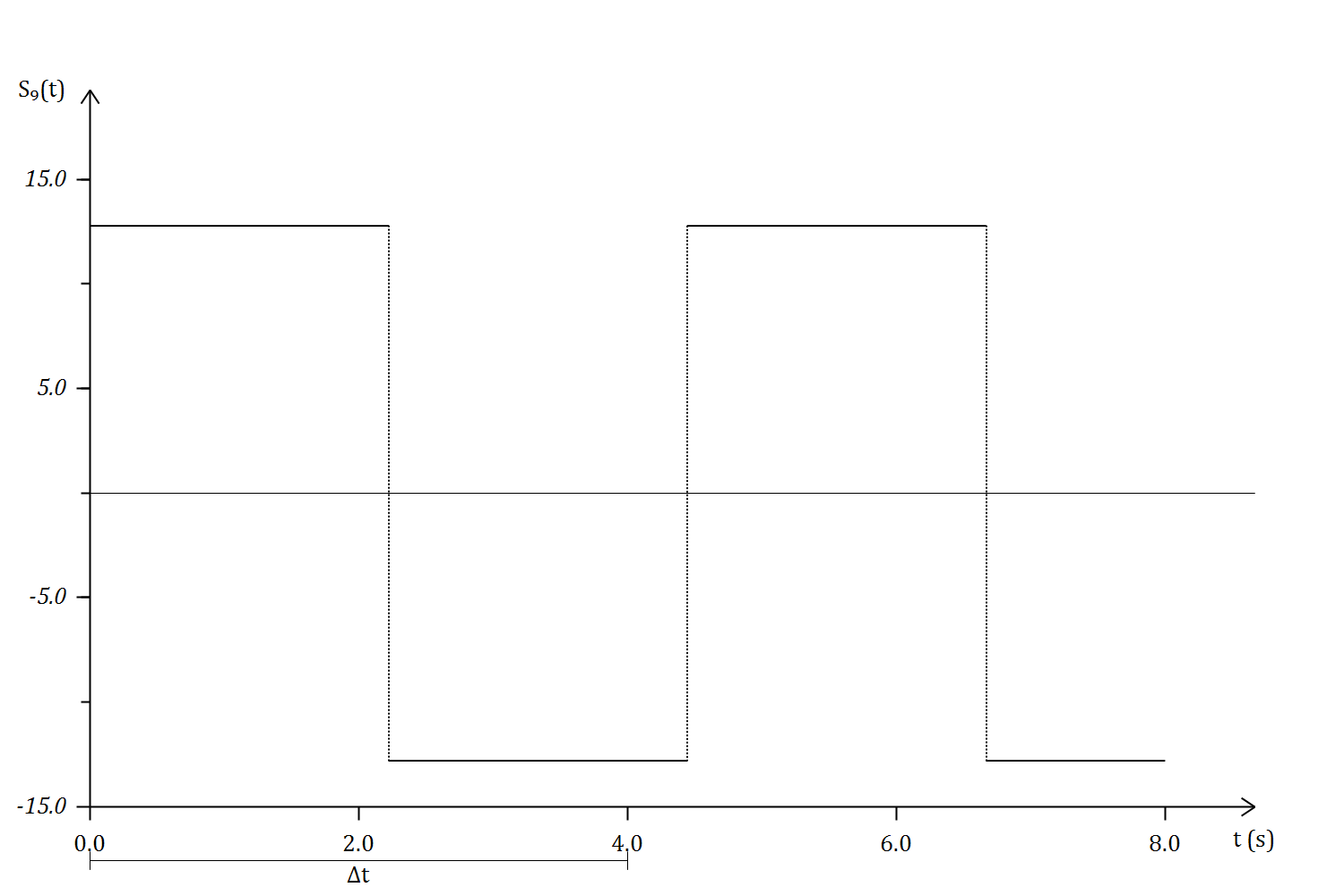}\label{f3i}}\qquad
\subfloat[: $S_{10}(t)$]{\includegraphics[width=4.8in]{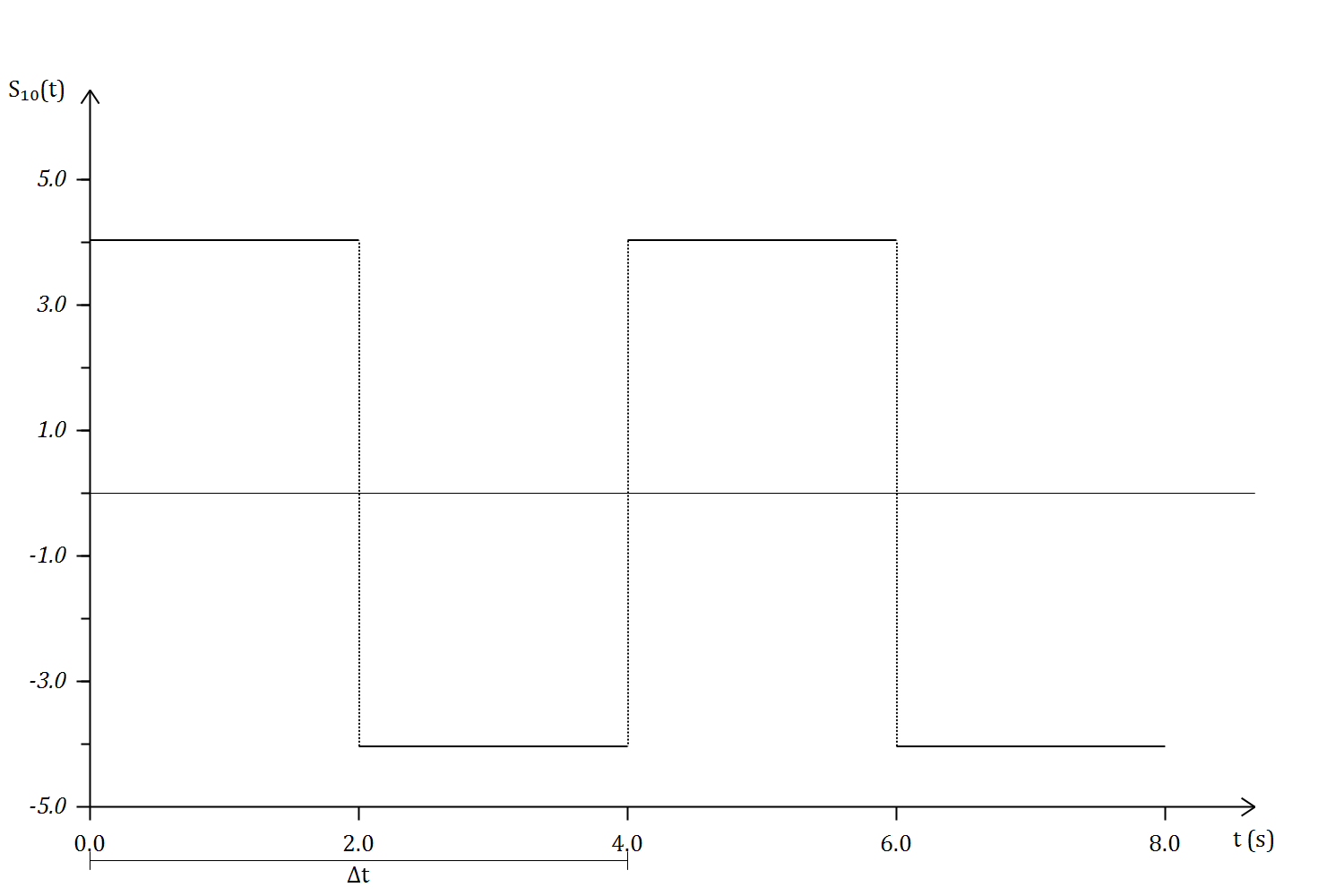}\label{f3j}}\qquad
\caption[]{}\end{figure}\begin{figure}\ContinuedFloat
\centering
\subfloat[: $S_{11}(t)$]{\includegraphics[width=4.8in]{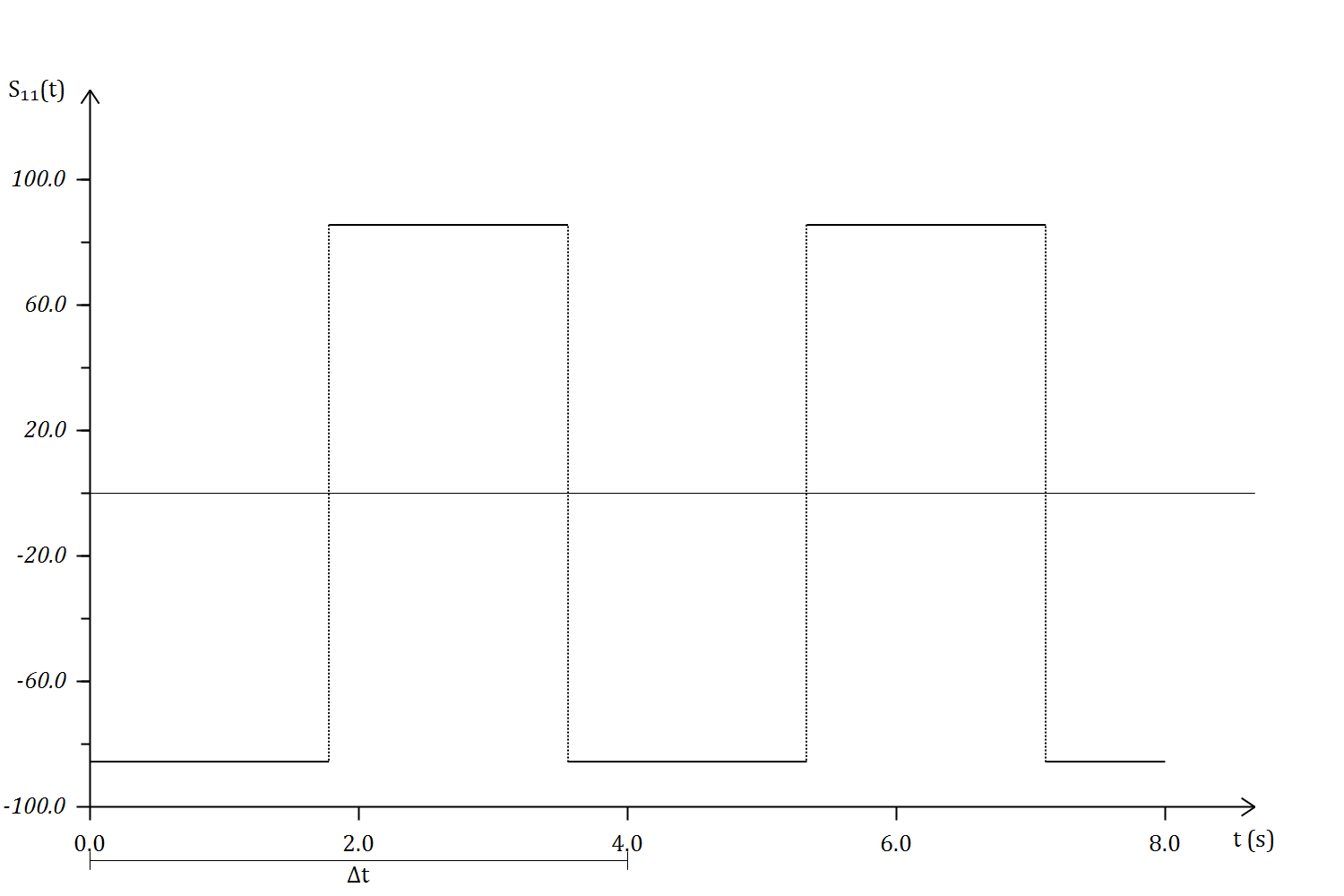}\label{f3k}}\qquad
\subfloat[: $S_{12}(t)$]{\includegraphics[width=4.8in]{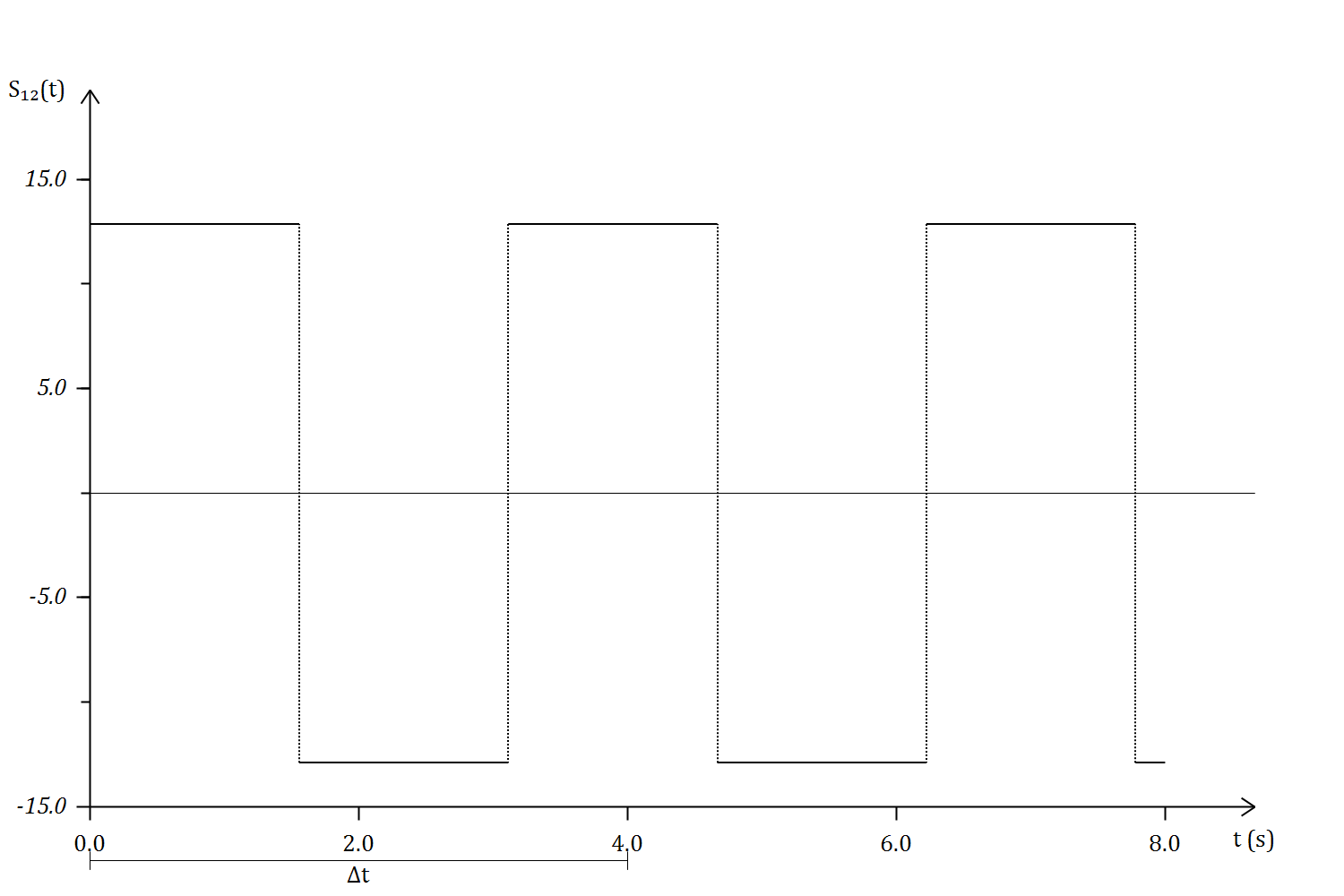}\label{f3l}}\qquad
\caption[]{}\end{figure}\begin{figure}\ContinuedFloat
\centering
\subfloat[: $S_{13}(t)$]{\includegraphics[width=4.8in]{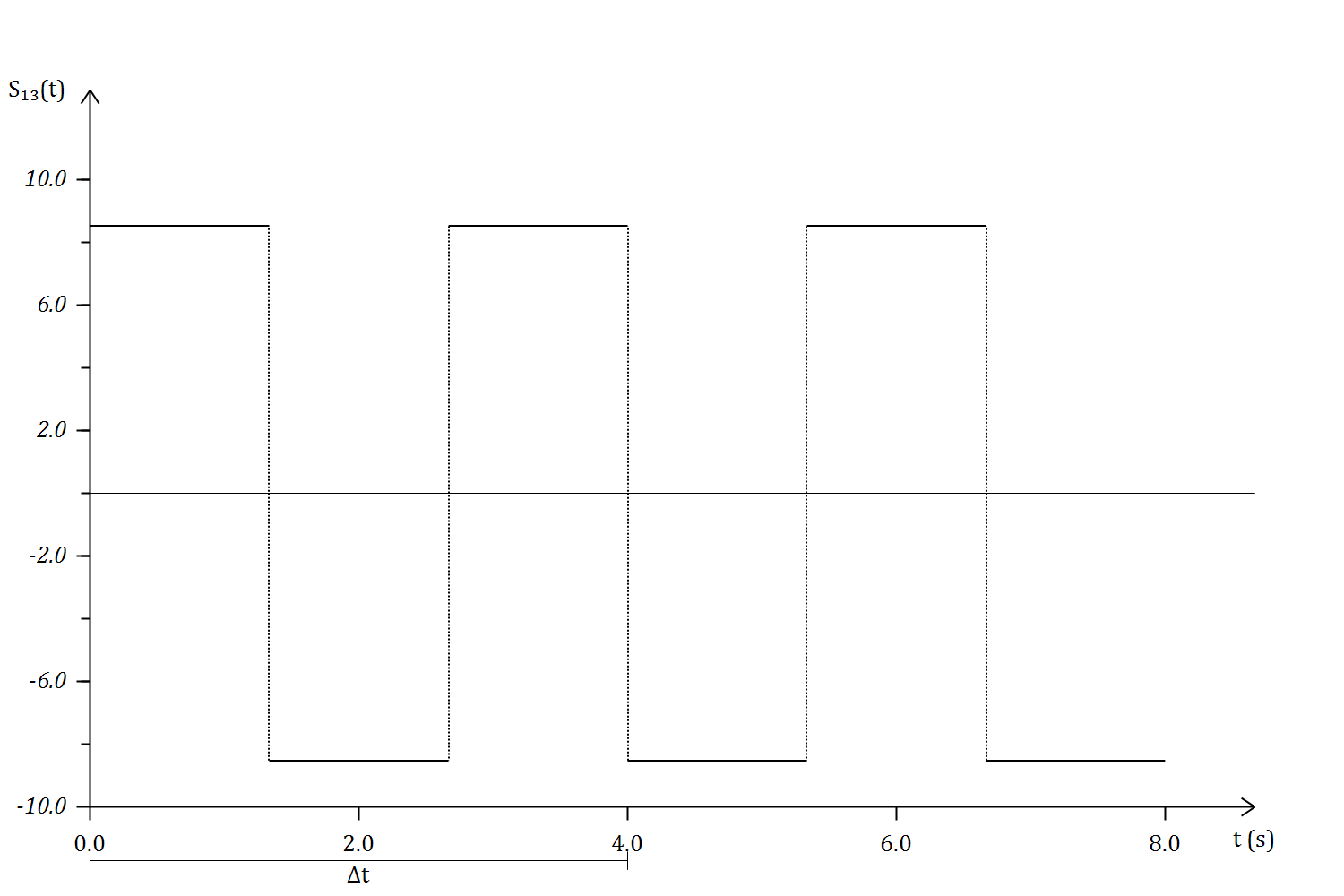}\label{f3m}}\qquad
\subfloat[: $S_{14}(t)$]{\includegraphics[width=4.8in]{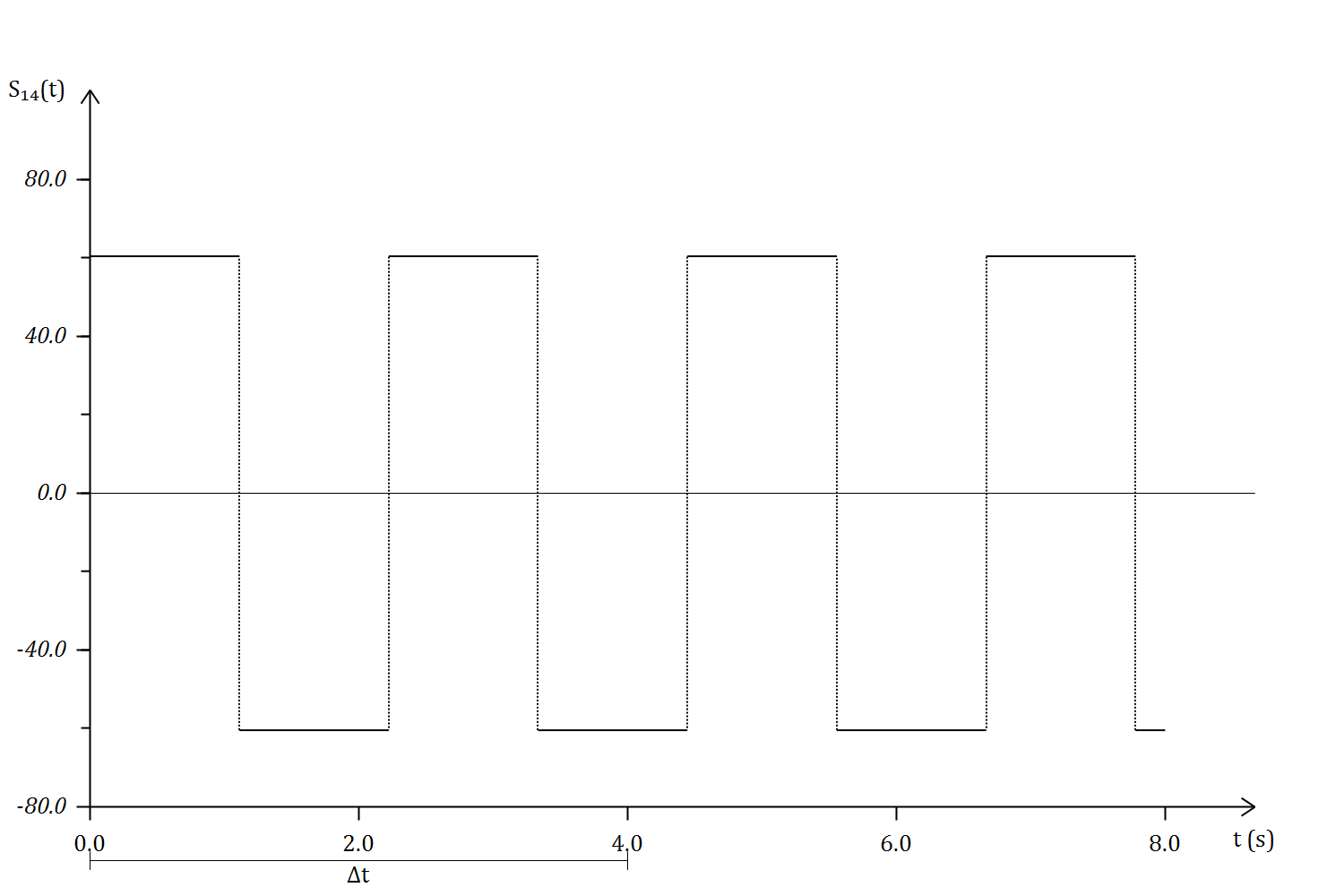}\label{f3n}}\qquad
\caption[]{}\end{figure}\begin{figure}\ContinuedFloat
\centering
\subfloat[: $S_{15}(x)$]{\includegraphics[width=4.8in]{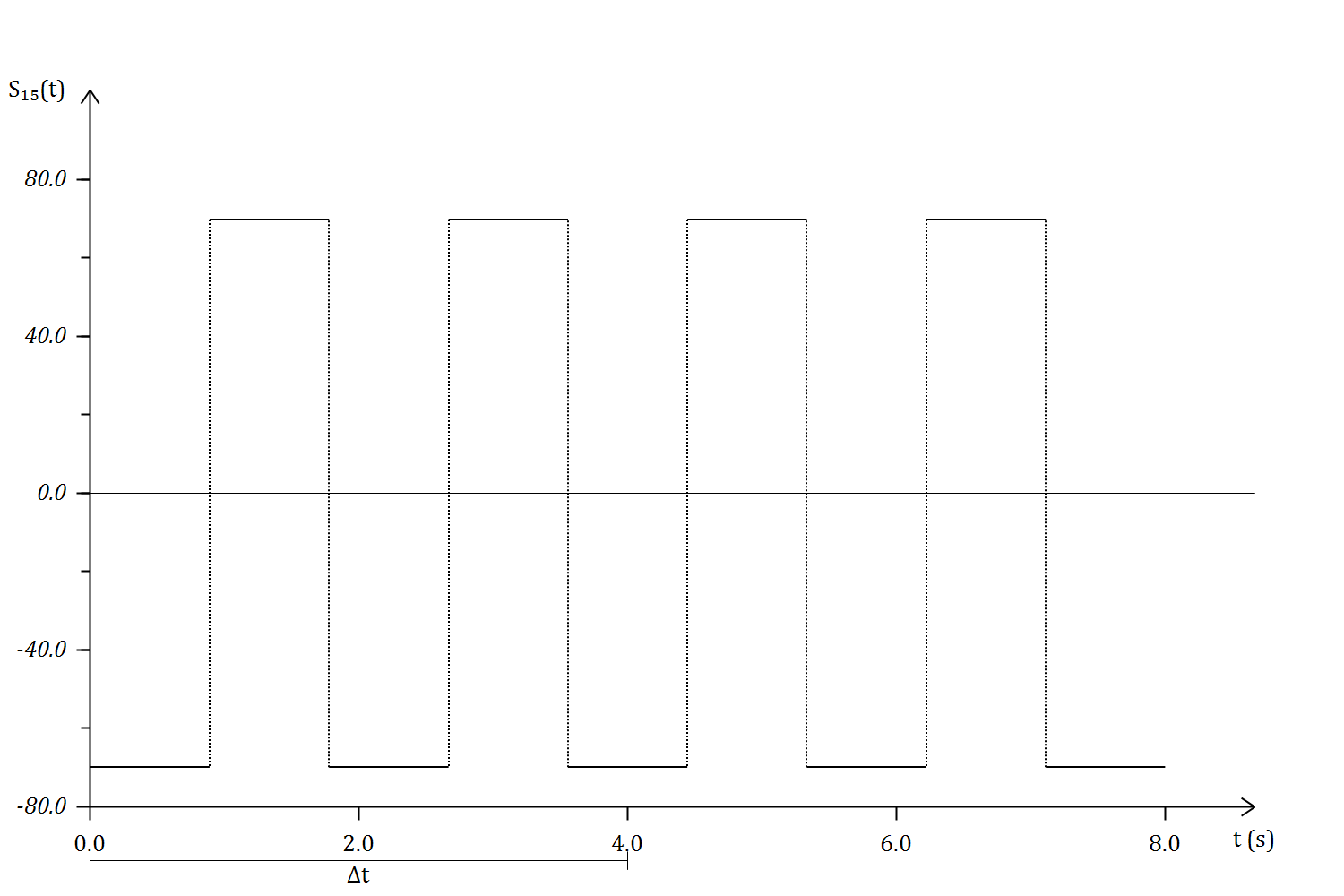}\label{f3o}}\qquad
\subfloat[: $S_{16}(t)$]{\includegraphics[width=4.8in]{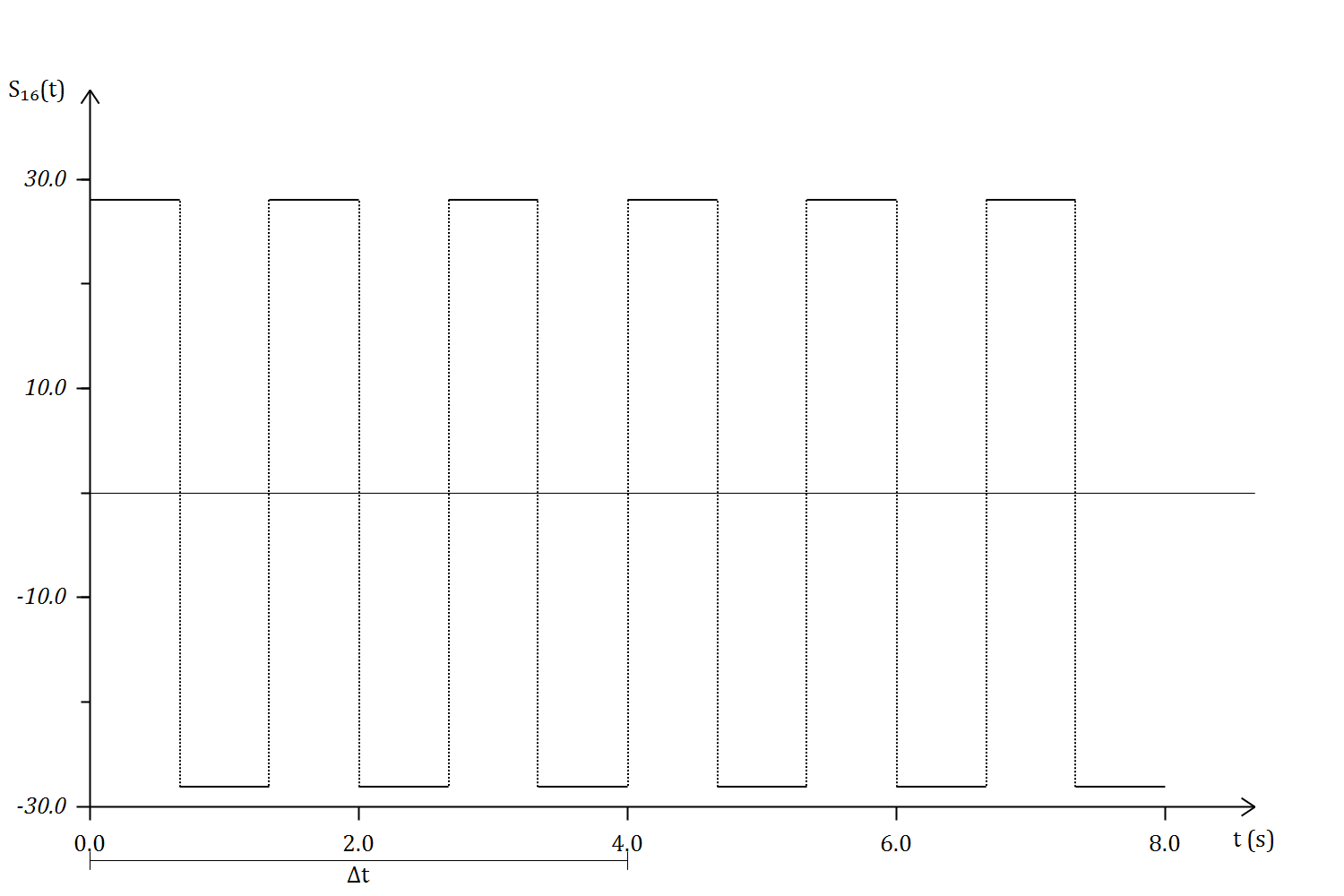}\label{f3p}}\qquad
\caption[]{}\end{figure}\begin{figure}\ContinuedFloat
\centering
\subfloat[: $S_{17}(t)$]{\includegraphics[width=4.8in]{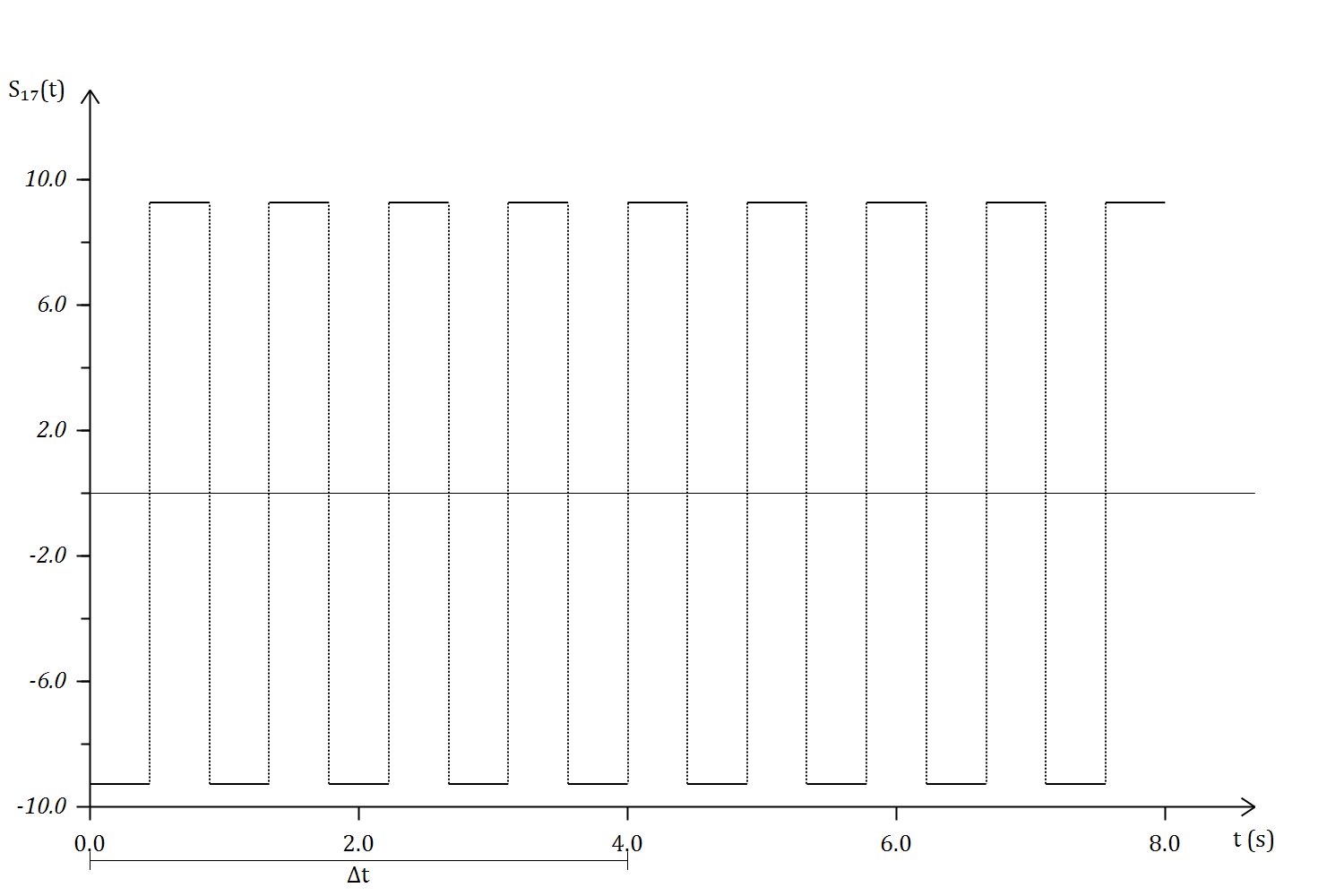}\label{f3q}}\qquad
\subfloat[: $S_{18}(t)$]{\includegraphics[width=4.8in]{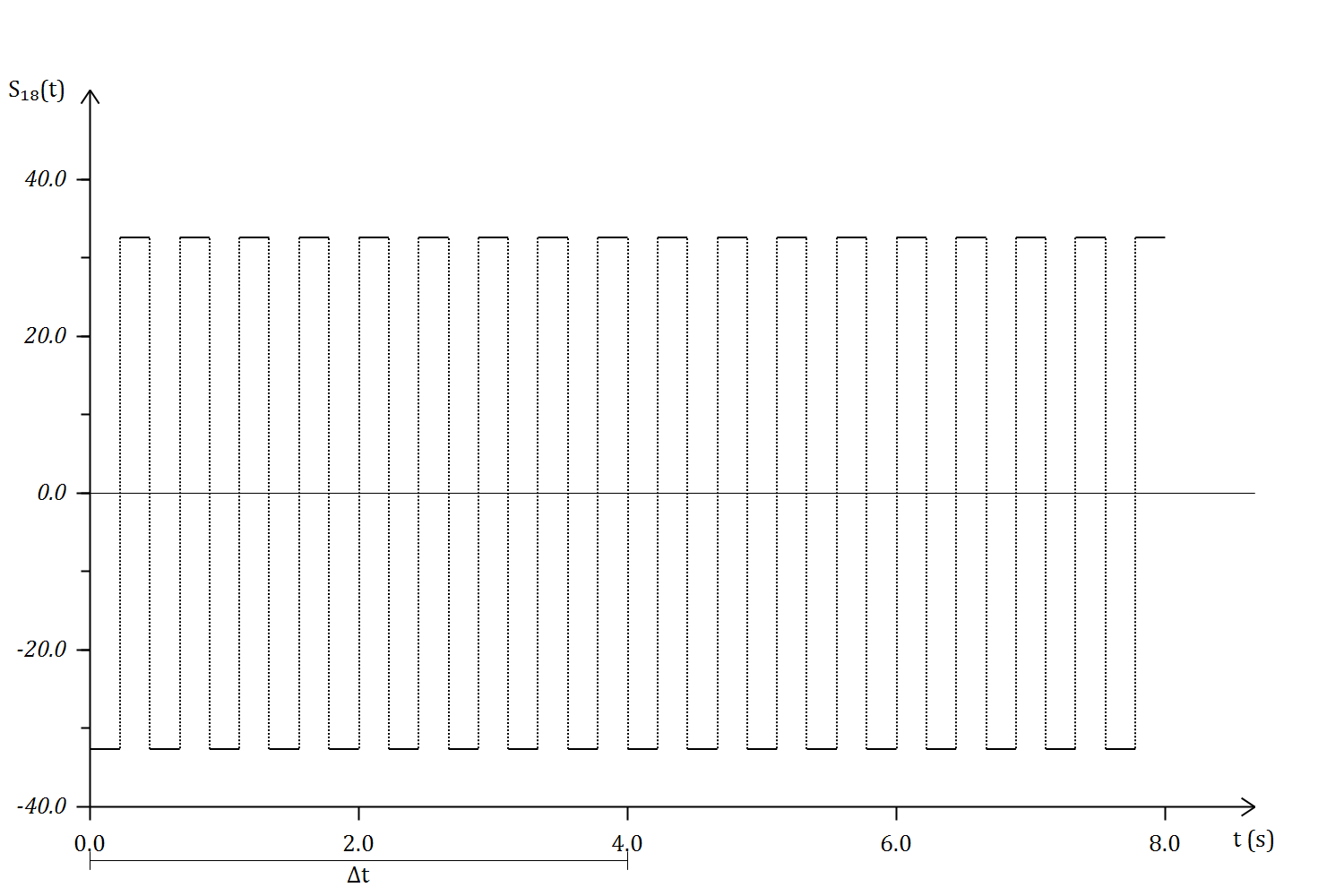}\label{f3r}}\qquad
\caption[]{Trains of square waves $S_1,S_2,\ldots$, $S_{17}$ and $S_{18}$}\end{figure}

\begin{figure}[H]
\centering
\includegraphics[width=5in]{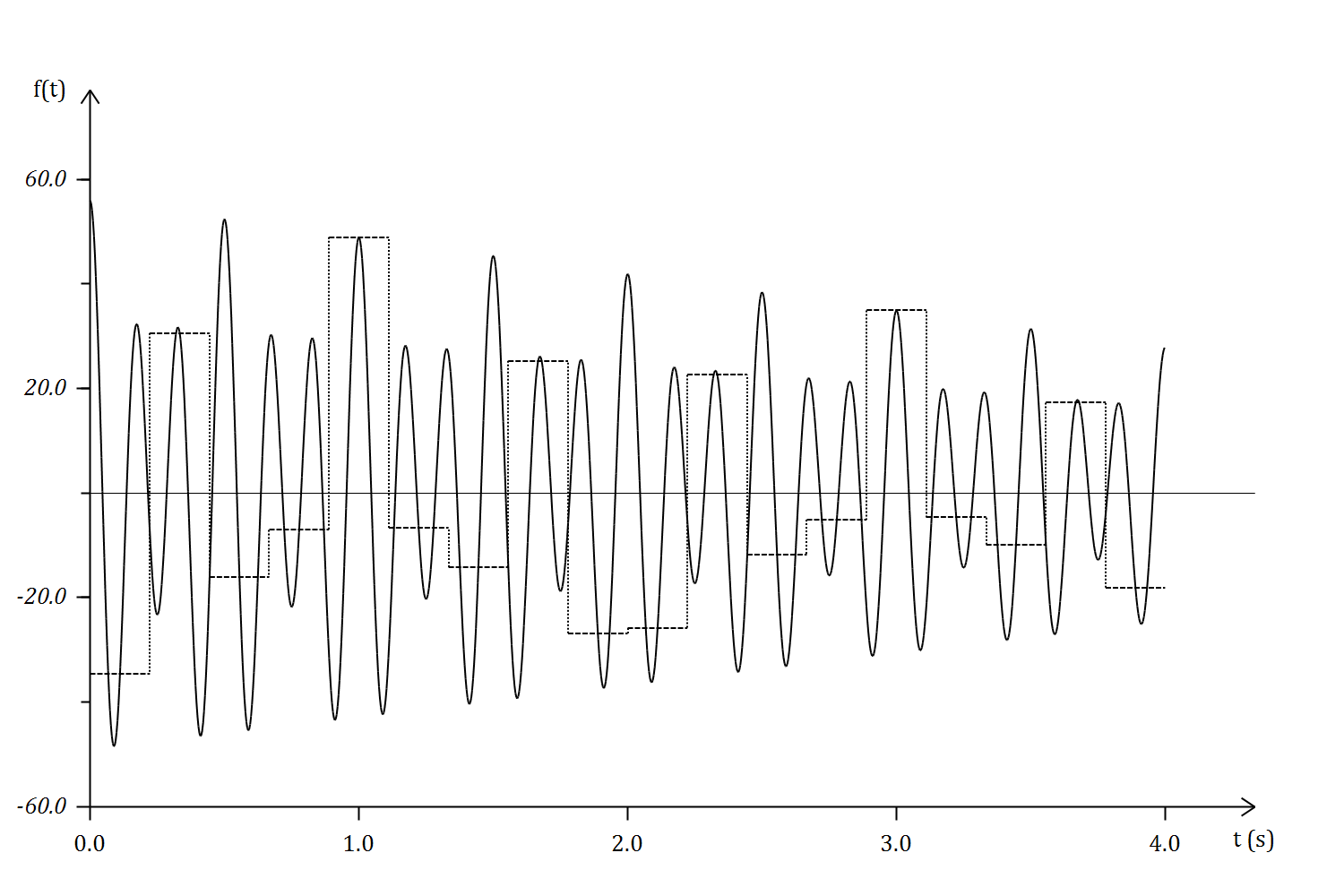}
\caption{The dashed line indicates the approximation to $f(t)$, specified in~\eqref{e1}, by $\displaystyle\sum\limits_{i=1}^{18} S_i(t)$.}
\label{f4}
\end{figure}

\newpage 
The approximation obtained for $f(t)$ (as specified in~\eqref{e1}, in interval $\Delta t$, by adding the 18 trains of square waves) is displayed in figure~\ref{f4}.

If one requires a better approximation to $f(t)$, by adding the trains of square waves, then $\Delta t$ should be divided into a larger number of equal sub-intervals. The higher the number of sub-intervals, the better the approximation.

Suppose that interval $\Delta t$ is divided into $n$ sub-intervals of equal duration. In equation~\eqref{e2} it was specified how to compute each $f_i$ (where $i=1,2,3,\ldots,n$) for each train of square waves $S_1, S_2, S_3, \ldots, S_n$, which must be added in $\Delta t$ to obtain, in that interval, the corresponding approximation to $f(t)$ specified in~\eqref{e1}.

To obtain the coefficients $C_1,C_2, C_3,\ldots, C_n$ corresponding respectively to those square waves, a system of linear algebraic equations must be solved. This system can be obtained by using the same type of approach as that used to obtain the system of equations specified in~\eqref{e3}.

\begin{figure}
\subfloat[: The dashed line indicates the approximation to $f(t)$, specified in~\eqref{e1}, by $\displaystyle\sum\limits_{i=1}^{100} S_i(t)$.]{\includegraphics[width=5in]{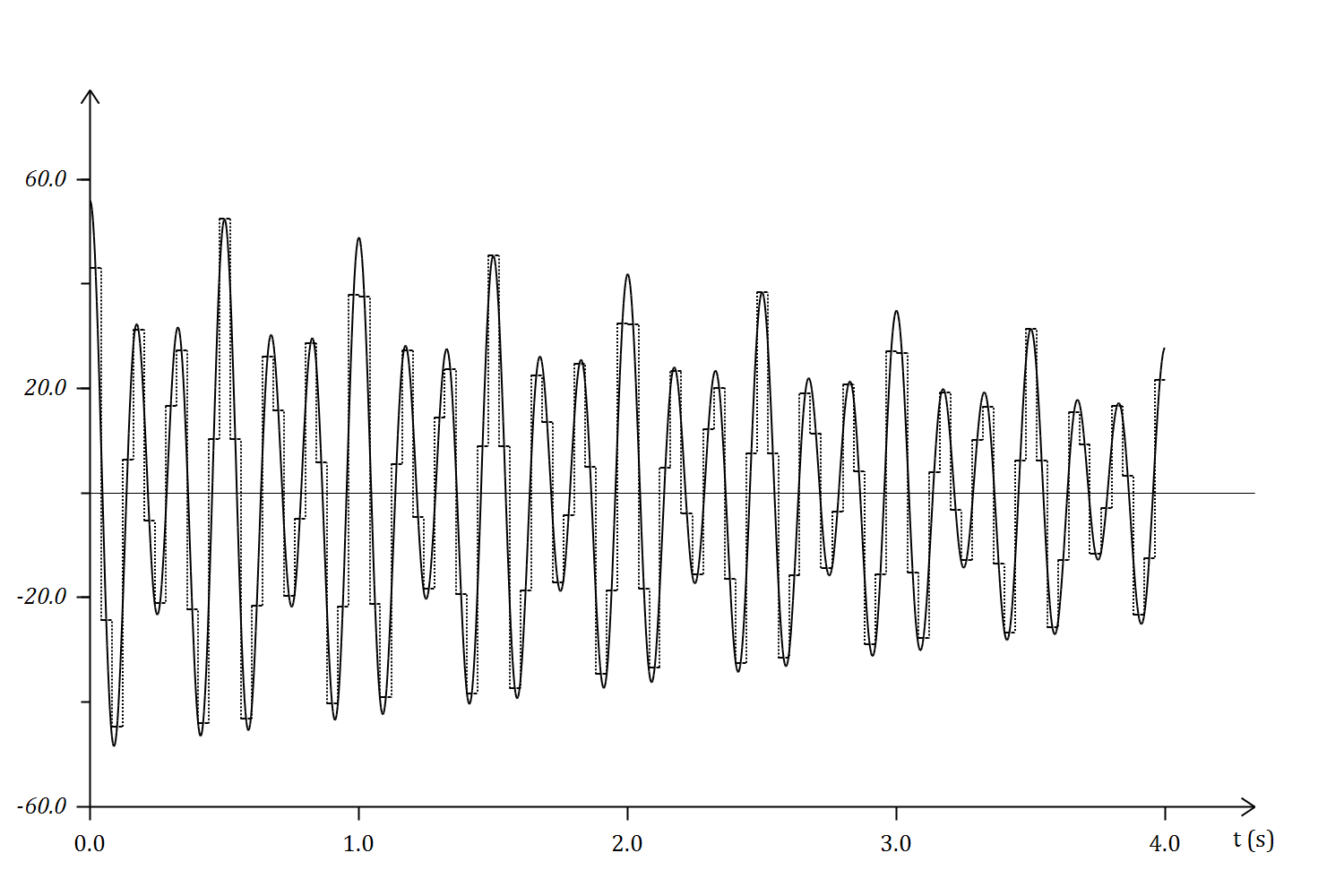}\label{f5a}}\qquad
\subfloat[: The dashed line indicates the approximation to $f(t)$, specified in~\eqref{e1} by $\displaystyle\sum\limits_{i=1}^{1000} S_i(t)$.]{\includegraphics[width=5in]{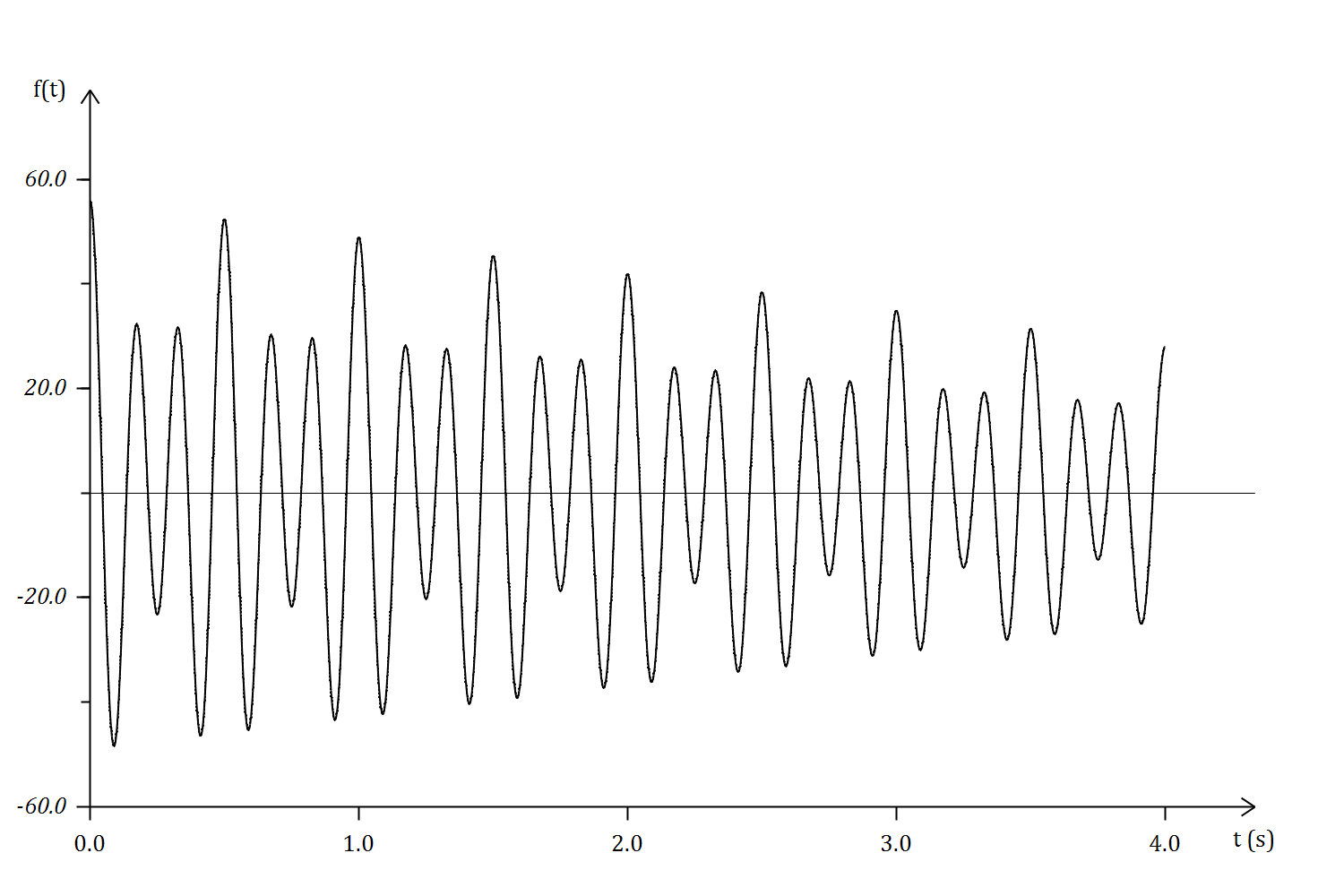}\label{f5b}}\qquad
\caption[]{}\end{figure}

Approximations to $f(t)$ specified in equation~\eqref{e1} when dividing $\Delta t$ into $100$ and into $1000$ intervals respectively, are displayed in~\subref{f5a} and \subref{f5b} in figure 5.

The SWM cannot be considered a branch of Fourier analysis; the trains of square waves $S_i$, where $i=1,2,3,\ldots,n$, do not make up a system of orthogonal functions.

The results obtained upon carrying out the type of analysis described of a function characterized in an interval $\Delta t$ divided into $n$ sub-intervals with an equal duration can be presented in a sequence of dyads (ordered pairs) such that the first element of the first dyad is the frequency $f_1$ corresponding to $S_1$ and the second element of the first dyad is the coefficient $C_1$; the first element of the second dyad is the frequency $f_2$ corresponding to $S_2$ and the second element of that dyad is the coefficient $C_2$; and so on successively, such that the first element of the $n^{\text{\tiny th}}$ dyad is the frequency $f_n$ corresponding to $S_n$ and the second element of that $n^{\text{\tiny th}}$ dyad is the coefficient $C_n$.

Consider, for example, the sequence of 18 dyads obtained when carrying out the type of analysis described of the $f(t)$ specified in~\eqref{e1}, if the interval $\Delta t$ is divided into 18 sub-intervals:
\begin{equation*}
\begin{aligned}[c]
(f_{1}; C_{1}) &= (0.1250000; 117.12980) \\
(f_{3}; C_{3}) &= (0.1406250; -210.98830) \\
(f_{5}; C_{5}) &= (0.1607143; 9.35088) \\
(f_{7}; C_{7}) &=(0.1875000; 61.27212)\\
(f_{9}; C_{9}) &=(0.2250000; 12.81335)\\
(f_{11}; C_{11}) &=(0.2812500; -85.68506)\\
(f_{13}; C_{13}) &= (0.3750000; 8.51973)\\
(f_{15}; C_{15}) &=(0.5625000; -69.86421)\\
(f_{17}; C_{17}) &=(1.1250000; -9.26140)\\
\end{aligned}
\qquad
\begin{aligned}[c]
(f_{2}; C_{2})&=,(0.1323529; 50.27631)\\
(f_{4}; C_{4}) &= (0.1500000; -53.27896) \\
(f_{6}; C_{6}) &=(0.1730769; 12.58025) \\
(f_{8}; C_{8}) &=(0.2045455; 49.80105)\\
(f_{10}; C_{10}) &=(0.2500000; 4.03101)\\
(f_{12}; C_{12}) &=(0.3214286; 12.88482)\\
(f_{14}; C_{14}) &=(0.4500000; 60.38772)\\
(f_{16}; C_{16}) &=(0.7500000; 28.08997)\\
(f_{18}; C_{18}) &=(2.2500000; -32.60758)\\
\end{aligned}
\end{equation*}

This approximation to function $f(t)$ specified in~\eqref{e1} can be expressed in the frequency domain. To achieve this objective, for each of the frequencies considered $f_1,f_2, f_3,\ldots,f_{18}$, the corresponding coefficients  $C_1,C_2, C_3,\ldots, C_{18}$ must be indicated. 

\begin{figure}[H]
\centering
\includegraphics[width=4.8in]{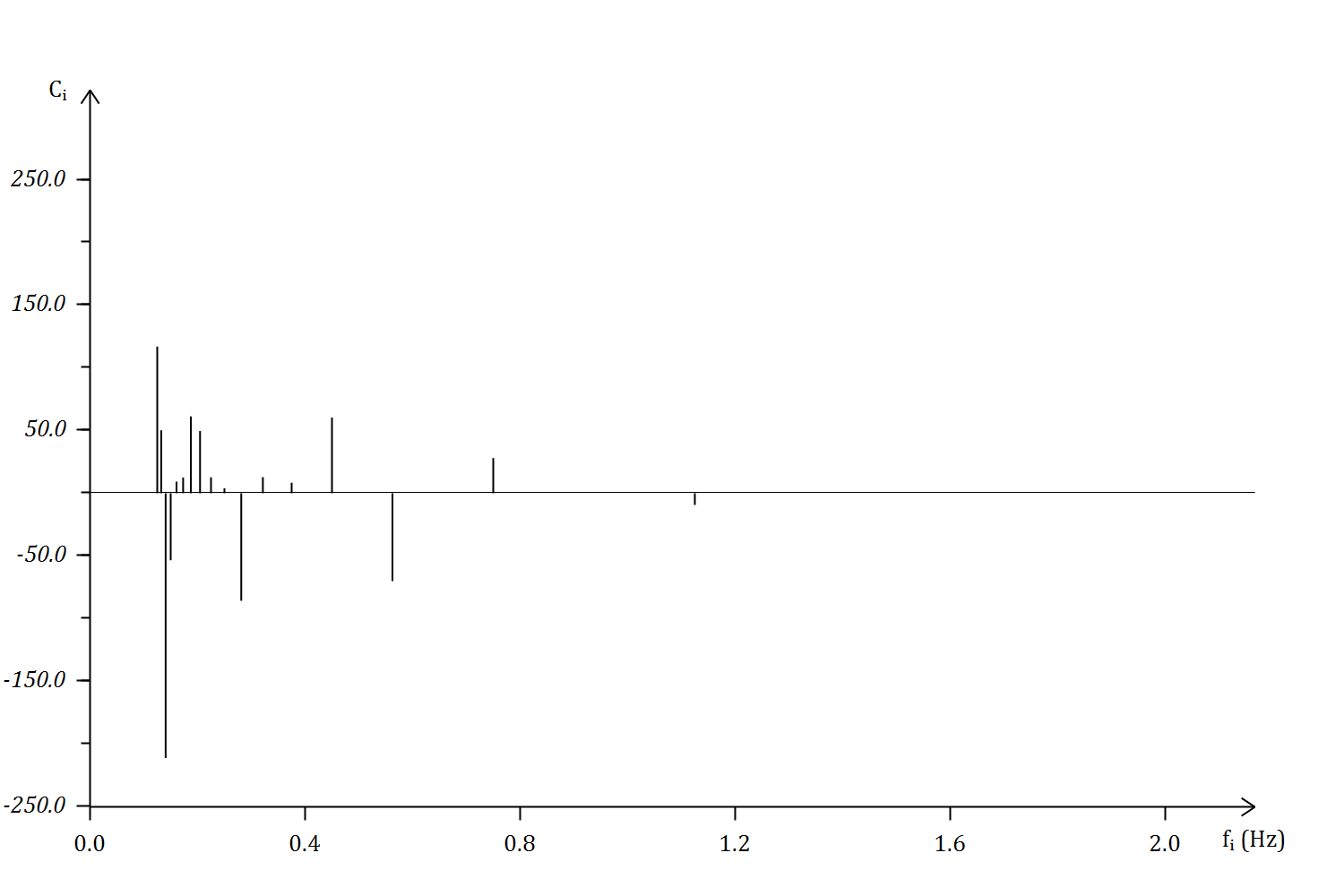}
\caption{SWT of the approximation to $f(t)$ obtained by dividing $\Delta t$ into 18 sub-intervals.}
\label{f6}
\end{figure}

The expression in the frequency domain of this approximation to $f(t)$ will be called the Square Wave Transform (SWT) of that approximation to $f(t)$. Note that previously (in figures~\ref{f4} and 5, for example), each approximation to $f(t)$ was represented in the time domain. This SWT is displayed in figure~\ref{f6}.

Of course, the SWTs corresponding to numbers as large as desired of equal sub-intervals into which $\Delta t$ is divided can be obtained for the $f(t)$ specified in~\eqref{e1}, or for any other function of the time which, in a particular interval $\Delta t$, satisfies the conditions of Dirichlet.

In (7.1), (7.2), and (7.3) of figure 7, the SWTs obtained for the approximations to the $f(t)$ specified in~\eqref{e1} are shown for $n=100$, $n=1000$, and $n=2000$, respectively. 

\newpage
\begin{figure}[H]
\subfloat[: SWT of the approximation to the $f(t)$, specified in equation~\eqref{e1}, for $n=100$.]{\includegraphics[width=5in]{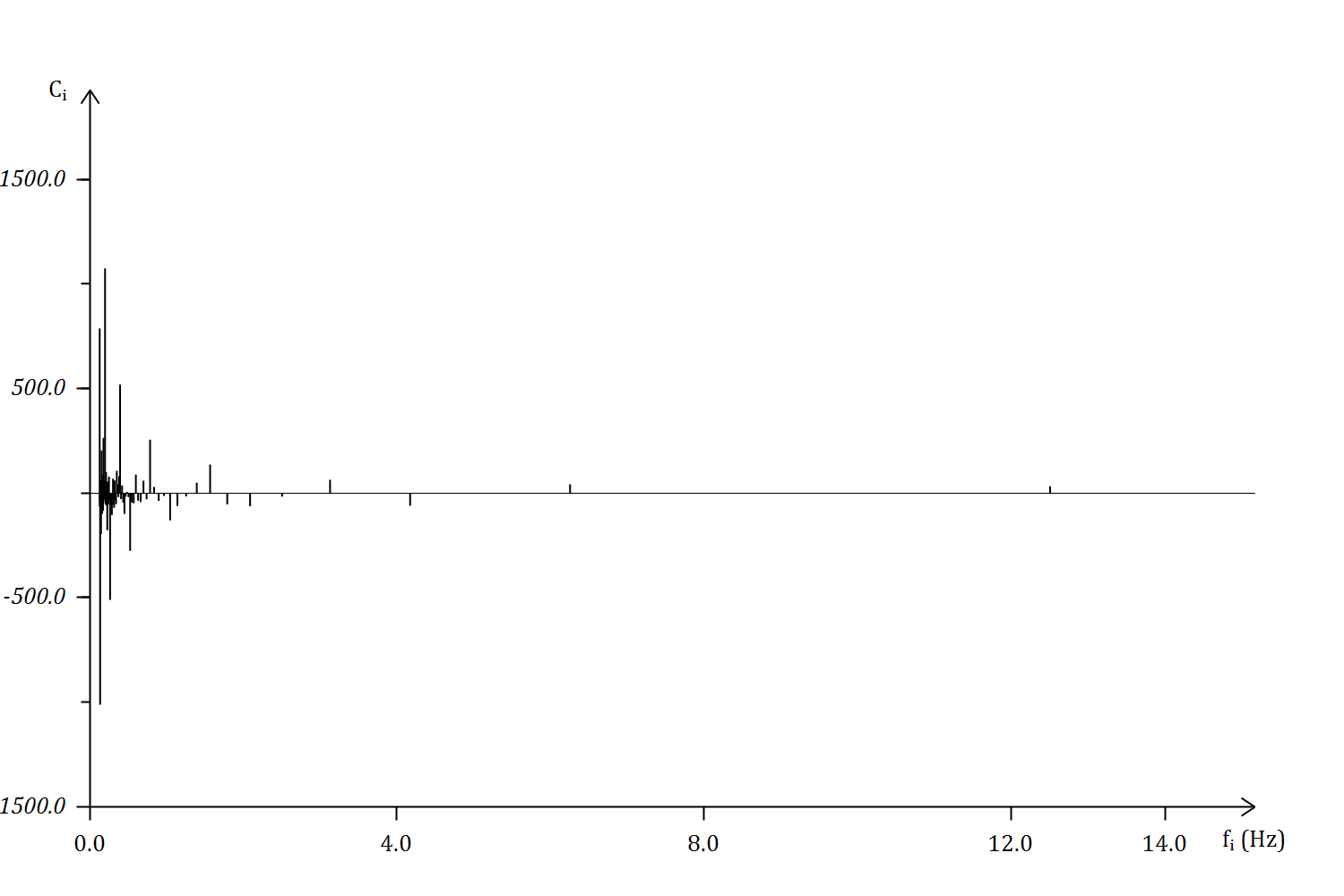}\label{f7a}}\qquad
\subfloat[: SWT of the approximation to the $f(t)$, specified in equation~\eqref{e1}, for $n=1000$.]{\includegraphics[width=5in]{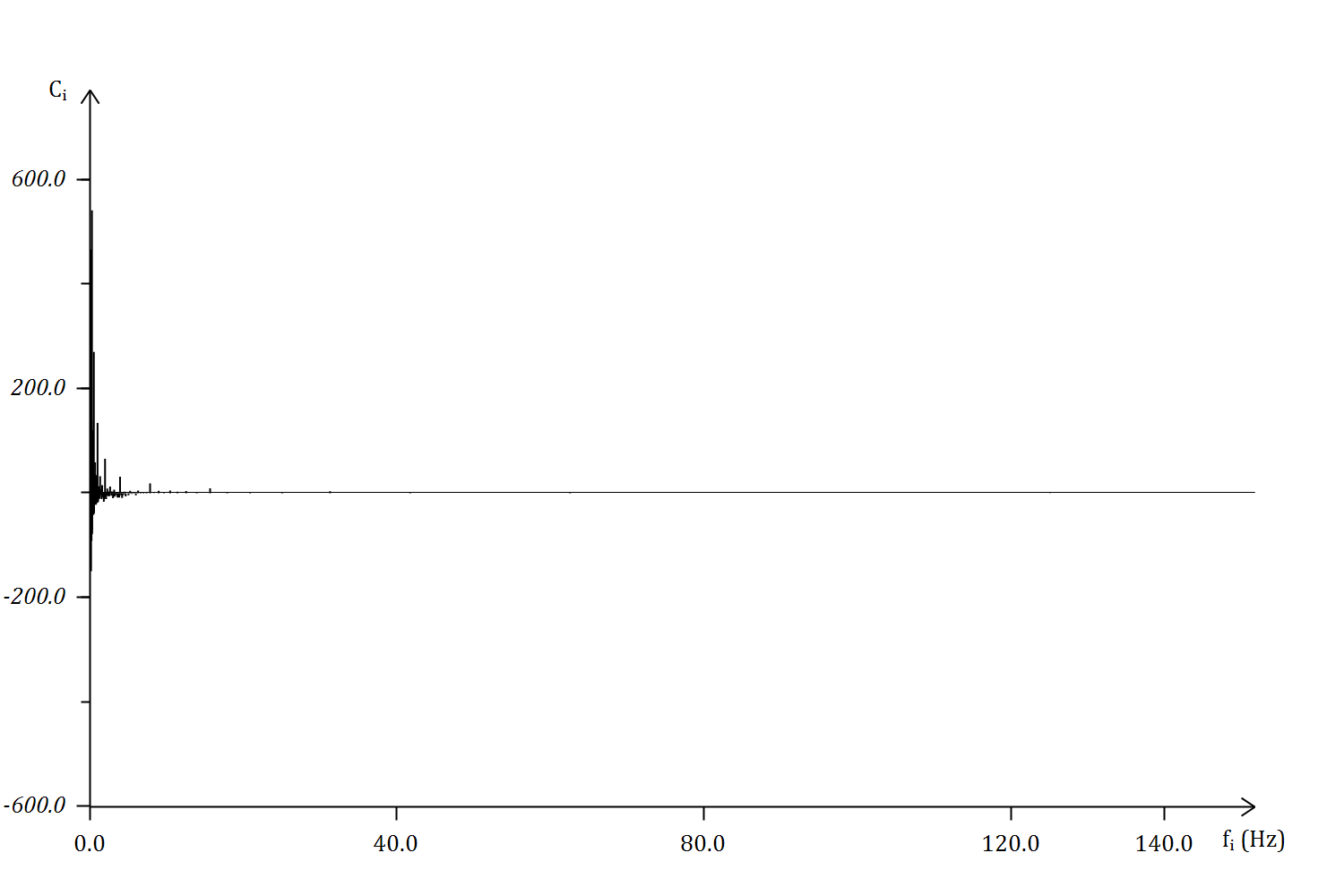}\label{f7b}}\qquad
\caption[]{}\end{figure}

\newpage
\begin{figure}
\ContinuedFloat
\subfloat[: SWT of the approximation to the $f(t)$, specified in equation~\eqref{e1}, for $n=2000$.]{\includegraphics[width=5in]{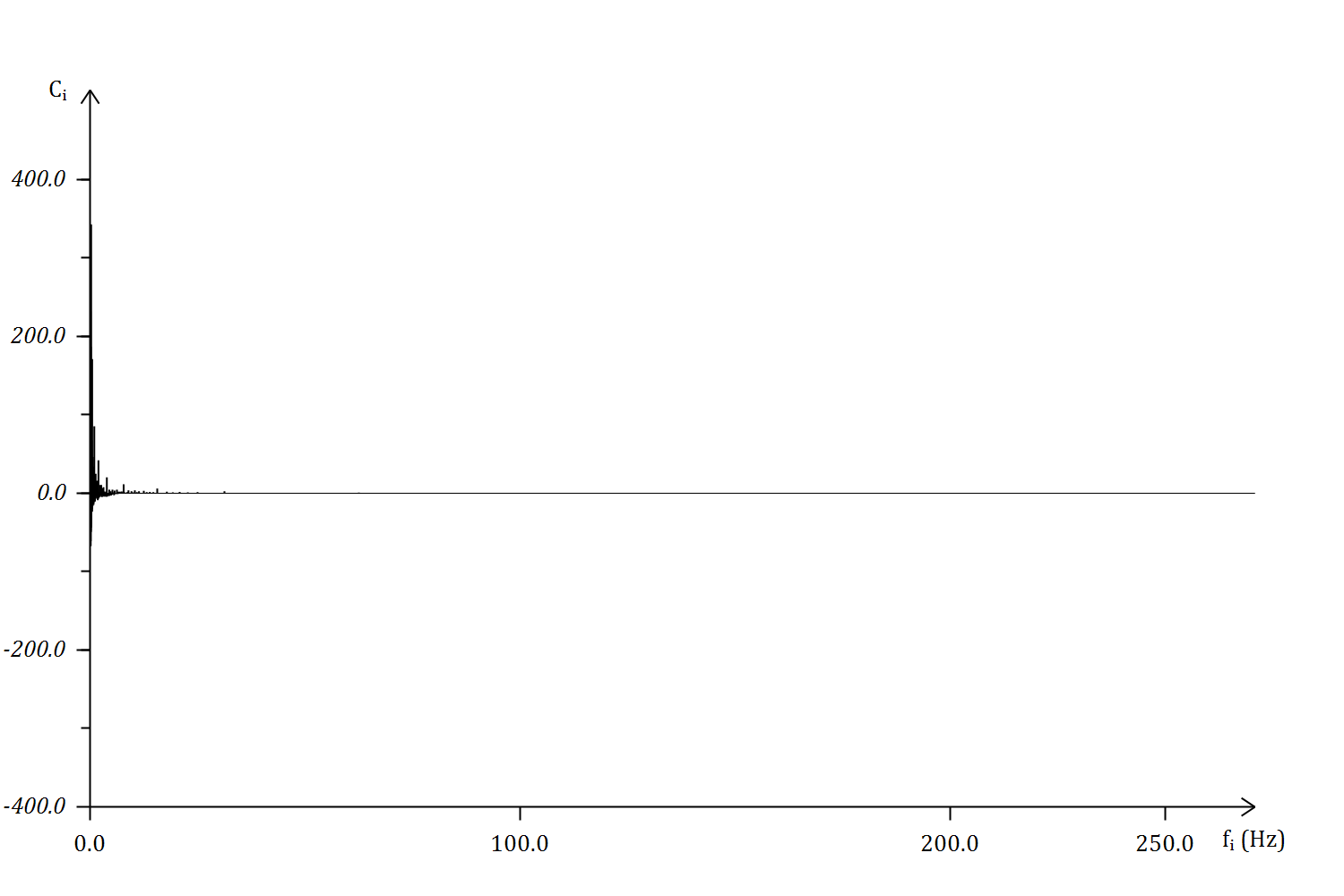}\label{f7c}}\qquad
\caption[]{}\end{figure}
In (8.1), (8.2), (8.3), and (8.4) of figure 8, partial representations can be seen (up to frequency $f_i=2$) of the SWTs of the different approximations to the $f(t)$ specified in~\eqref{e1}.

\newpage
\begin{figure}[H]
\subfloat[: Partial representation of the SWT of the approximation to the $f(t)$ specified in~\eqref{e1} for $n=1000$]{\includegraphics[width=5in]{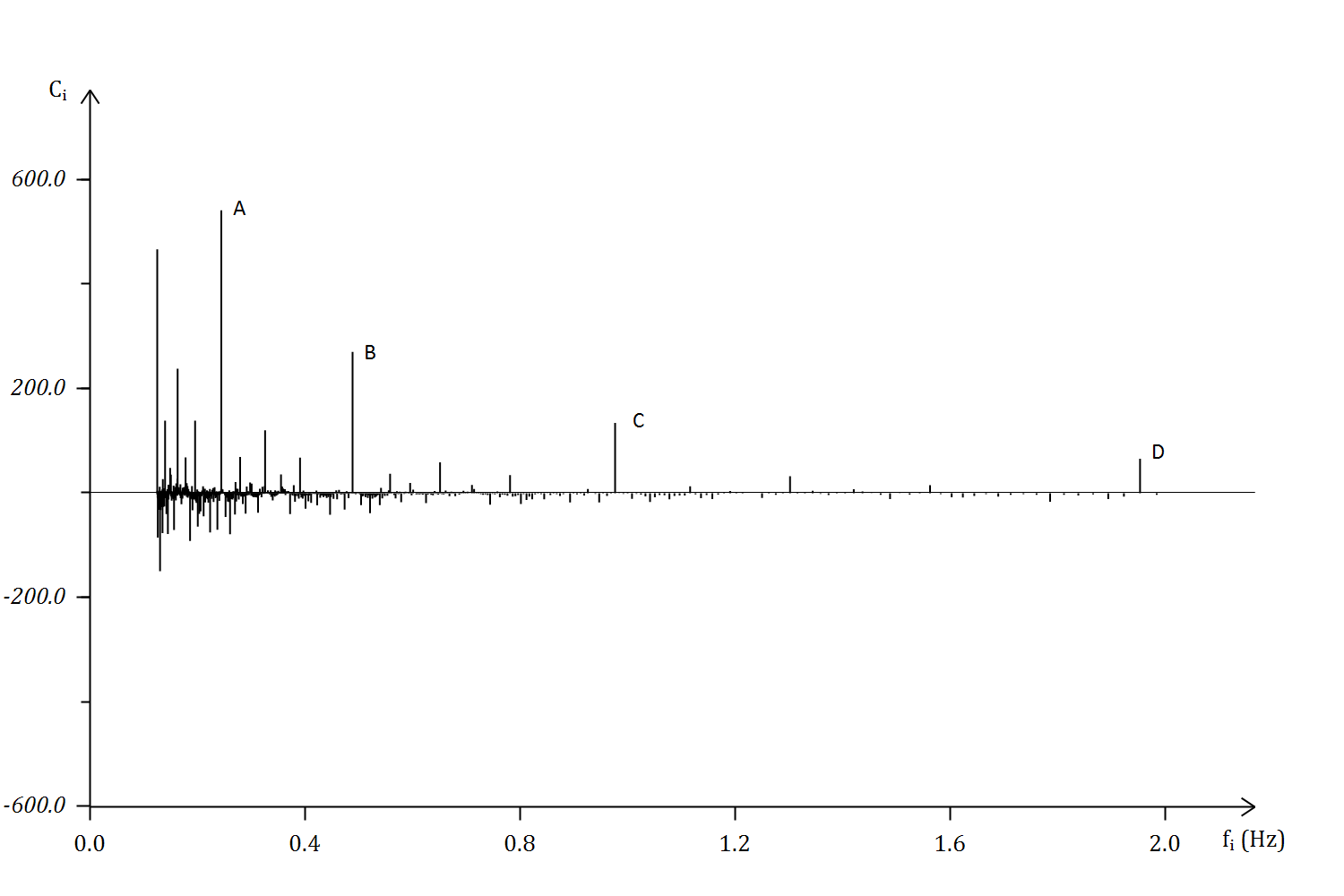}\label{f8a}}\qquad
\subfloat[: Partial representation of the SWT of the approximation to the $f(t)$ specified in~\eqref{e1} for $n=2000$]{\includegraphics[width=5in]{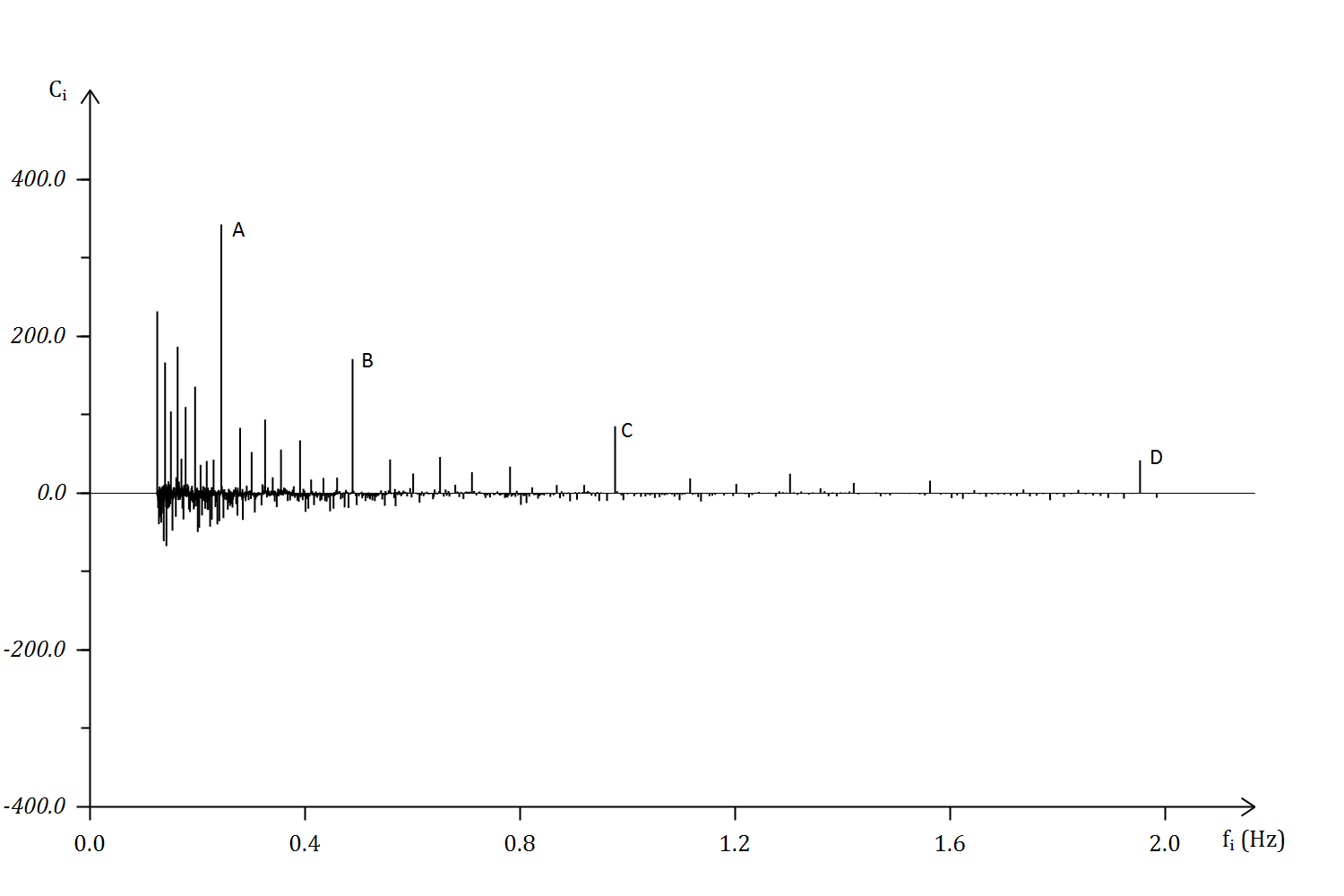}\label{f8b}}\qquad
\caption[]{}\end{figure}

\begin{figure}[H]
\ContinuedFloat
\subfloat[: Partial representation of the SWT of the approximation to the $f(t)$ specified in~\eqref{e1} for $n=4000$]{\includegraphics[width=5in]{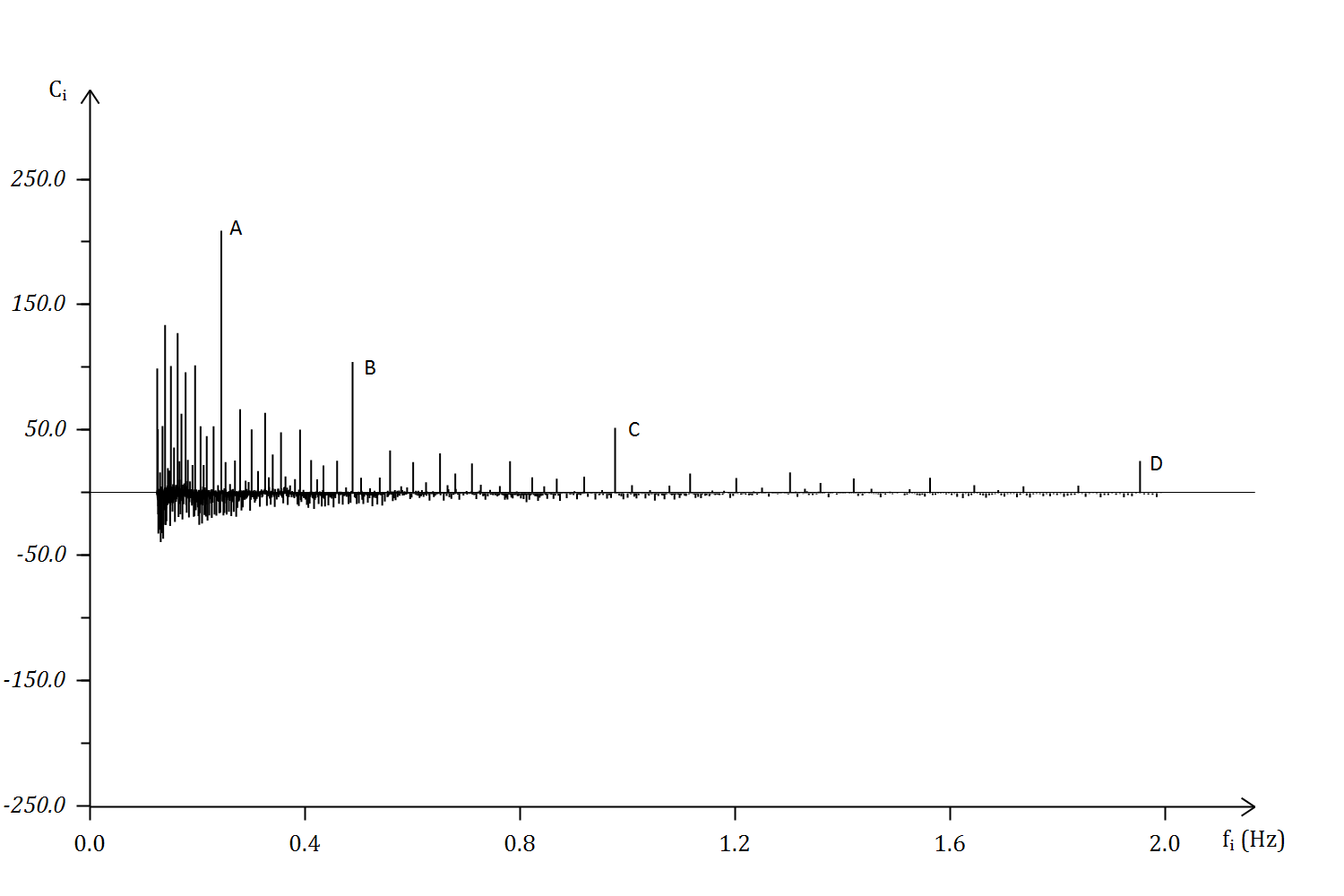}\label{f8c}}\qquad
\subfloat[: Partial representation of the SWT of the approximation to the $f(t)$ specified in~\eqref{e1} for $n=8000$]{\includegraphics[width=5in]{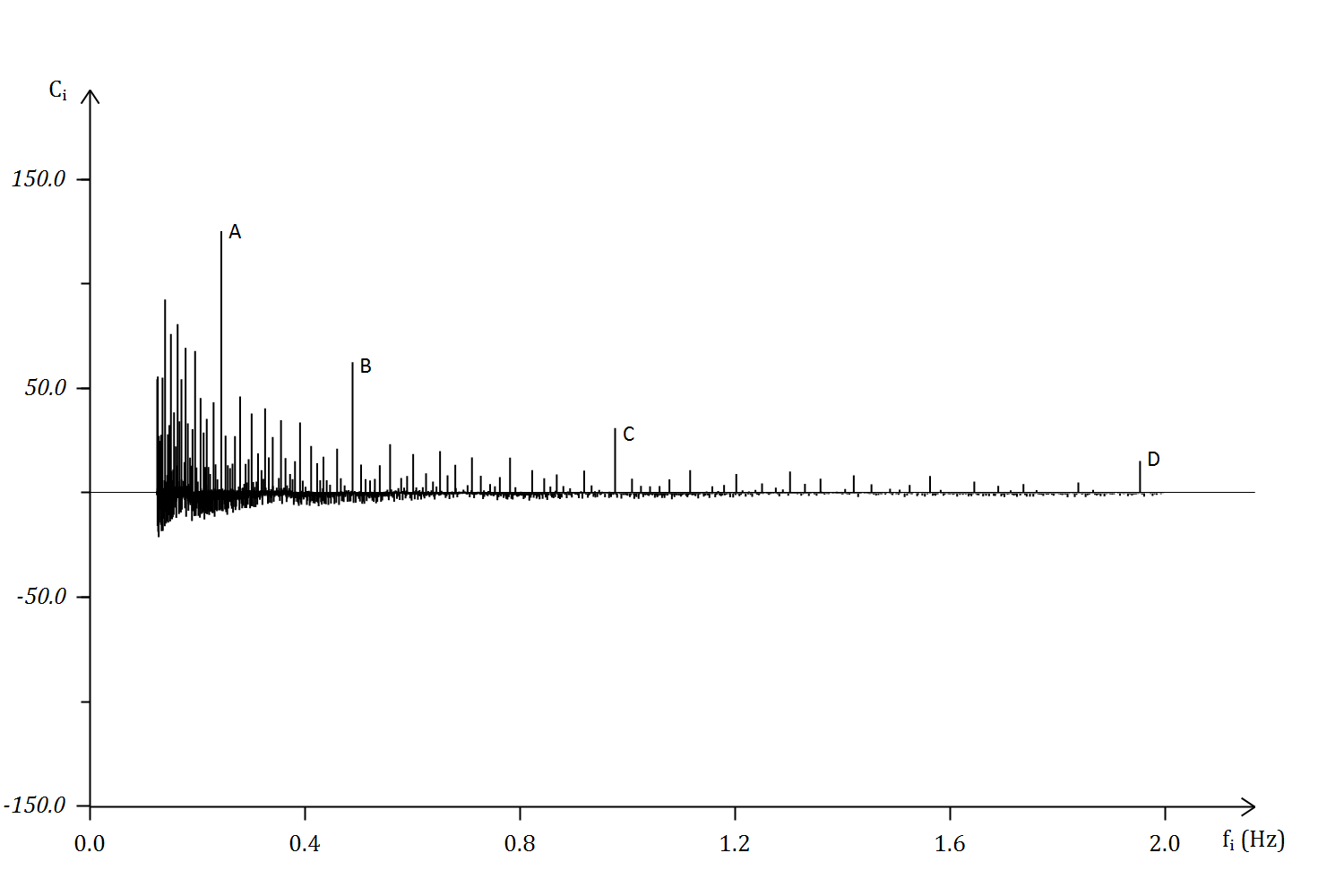}\label{f8d}}\qquad
\caption[]{}\end{figure}

In figure 8, it can be observed that in all 4 cases considered the ``prominent coefficients'' correspond to certain frequencies. (They are \textit{prominent} in the sense that their moduli are quite larger than the moduli of the coefficients corresponding to frequencies near to those considered.) They have been indicated by the letters A, B, C, and D. The dyads corresponding to the coefficients for these cases are the following:
\begin{align*}
& n = 1000 \\
A\colon&(0.2441410; 541.50054)\\ 
B\colon&(0.4882812; 270.06817)\\
C\colon&(0.9765633; 134.80741)\\
D\colon&(1.9531250;  66.39768)\\
& n=2000 \\
A\colon&(0.2441410; 342.97162)\\ 
B\colon&(0.4882812; 171.42003)\\
C\colon&(0.9765633;  85.45887)\\
D\colon&(1.9531250;  42.22607)\\
& n=4000 \\
A\colon&(0.2441410; 209.36325)\\ 
B\colon&(0.4882812; 104.64873)\\
C\colon&(0.9765633;  52.19877)\\
D\colon&(1.9531250;  25.84766)\\
& n=8000 \\
A\colon&(0.2441410; 541.40777)\\ 
B\colon&(0.4882812;  62.68744)\\
C\colon&(0.9765633;  31.28092)\\
D\colon&(1.9531250;  15.51459)
\end{align*}

Note that in these 4 cases, the frequency corresponding to coefficient A is equal to 0.2441410; the frequency corresponding to the coefficient B is equal to 0.4882812; the frequency corresponding to the coefficient C is equal to 0.9765633; and the frequency corresponding to the coefficient D is equal to 1.9531250. In other words, the frequencies of those ``prominent coefficients'' are invariable in the changes specified for the value of $n$ from 1000 to 2000, from 2000 to 4000, and from 4000 to 8000. 

\section{Analysis of Sequences of Samples (Measured Values) from an Electromyographic Recording}

The SWM and the corresponding expression of the results obtained with the SWT can be used for the analysis of sequences of ``samples'', or measured values, from different types of recordings. In particular, they can be used for the analysis of recordings which are important in medicine, as are those of the electrocardiogram (ECG), the electroencephalogram (EEG), and the electromyogram (EMG). Results obtained by using the SWM for the analysis of a sequence of samples from an electroencephalographic recording were described in \cite{b2}.

Suppose that, to take samples from a recording, a sampling frequency $f_s$ of 250 Hz is used. In this case let it be admitted that every one-second lapse (1 s) has been divided into 250 sub-intervals of equal duration and that the values of the samples ``correspond'' to the values (mentioned in section 2) of the function analyzed, at the midpoints of those different sub-intervals. In general, for any numeric value of the $f_s$ used, the number of sub-intervals into which the unit of time 1 s has been divided is equal to that numeric value.

Thus if with a specific $f_s$, such as one of 250 Hz, a sequence of samples is taken for a certain $\Delta t$, equal to 5 s, for example, to compute the total number $n$ of samples taken during that $\Delta t$, the $f_s$ must be multiplied by that $\Delta t$:  
\begin{equation}
n=f_s\cdot\Delta t
\label{e4}
\end{equation}

For the values of $f_s$ and $\Delta t$ specified in the above paragraph, the following equality is found:
\begin{equation*}
1250=250\,\mathrm{s}^{-1}\cdot 5\,\mathrm{s}
\end{equation*}

The SWM may be used to analyze a function $f(t)$ characterized analytically in an interval $\Delta=5$ s. If $n$ (the number of sub-intervals into which $\Delta t$ is divided) is equal to 1250, for example, the sequence of values $V_1, V_2, V_3, \ldots, V_{1250}$ will be obtained for the midpoints of the specific sub-intervals as follows: $V_1$ will be the value of $f(t)$ at the midpoint of the first sub-interval, $V_2$ will be the value of $f(t)$ at the midpoint of the second sub-interval, and so on successively, with $V_{1250}$ the value of $f(t)$ for the last of the sub-intervals considered.

When operating with $f_s = 250$ Hz, and taking samples for $\Delta t = 5$ s, 1250 samples will be obtained. To apply the SWM to this sequence of samples, it is considered that the first ``corresponds'' to $V_1$, the second to $V_2$, and so on, successively such that the last sample ``corresponds'' to $V_{1250}$. In other words, the sequence of samples is treated the same as the sequence $V_1, V_2, V_3, \ldots, V_{1250}$, when applying the SWM to an analytically characterized function.

Equation \eqref{e2}, which makes it possible to compute each $f_i$ corresponding to each $S_i$ (where $i=1,2,3,\ldots,n$), remains valid for sequences of samples of a recording:

\begin{equation}
f_{i}=\frac{1}{2\Delta t}\left(\frac{n}{n-(i-1)}\right)\:\mathrm{s}^{-1};\quad i=1,2,\ldots,n 
\tag{\ref{e2}}
\end{equation}

If \eqref{e4} is substituted in \eqref{e2}, the following equation~\ref{e5} is obtained:
\begin{equation}
\label{e5}
f_i=\frac{1}{2}\left(\frac{f_s}{f_s\Delta t - (i-1)}\right);\quad i=1,2,3,\ldots,n
\end{equation}
 
Note the equations obtained if in~\eqref{e5}, a) $i=1$; and b) $i=n$. 

For a), $(i=1) \longrightarrow f_1=\frac{1}{2\Delta t}$. In this first case, $f_1$ does not depend on $f_s$, but rather only on $\Delta t$.

For b), $(i=n) \longrightarrow f_n=\frac{f_s}{2}$. In this second case, $f_n$ does not depend on $\Delta t$, but rather only on $f_s$.

\begin{figure}[H]
\centering
\includegraphics[width=5in]{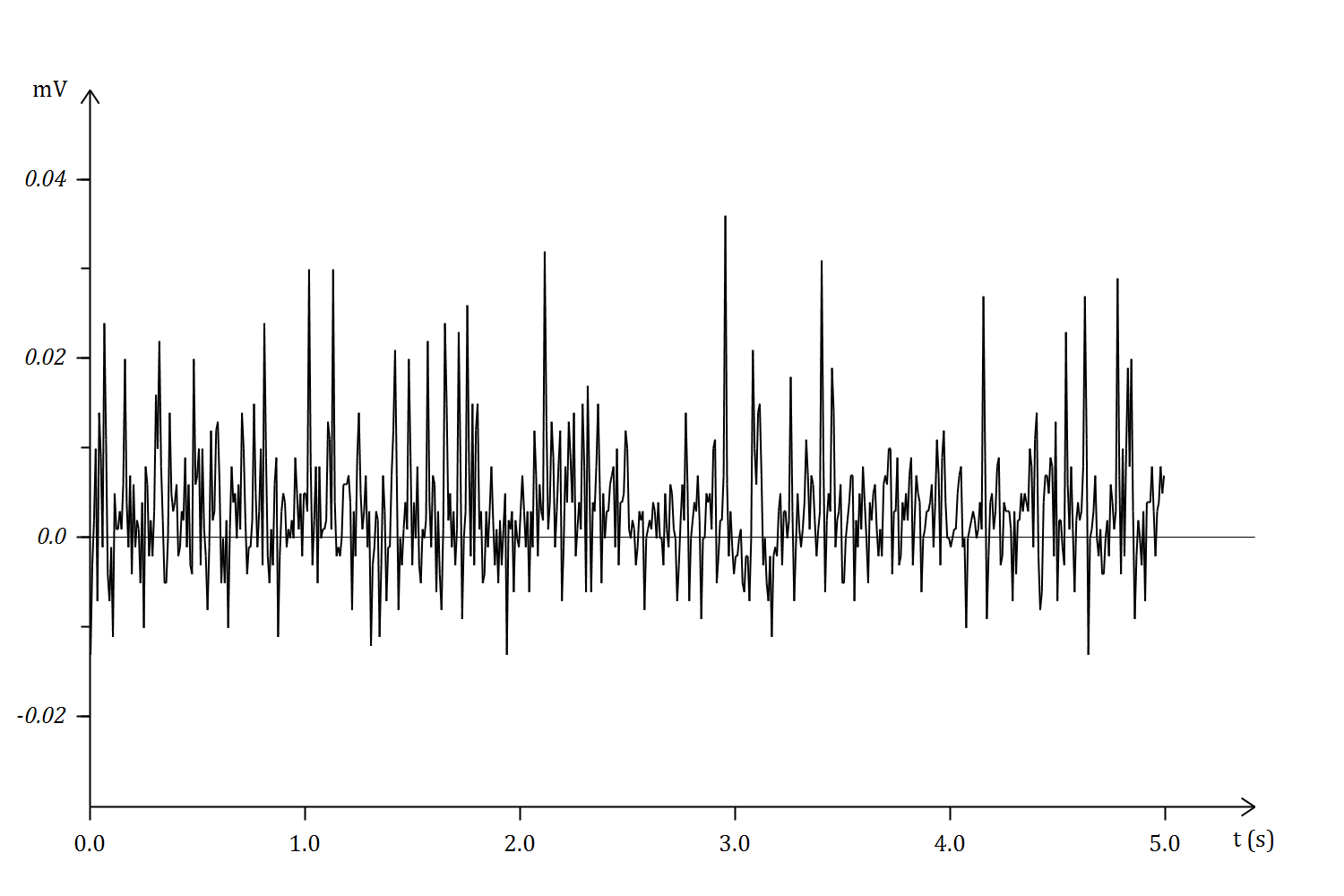}
\caption{Electromyographic recording made during $\Delta t=5\,\mathrm{s}$ with $f_s=250\,\mathrm{Hz}$.}
\label{f9}
\end{figure}

The graph corresponding to an electromyographic recording made during $\Delta t=5\,\mathrm{s}$ with $f_s=250\,\mathrm{Hz}$ is presented in figure \ref{f9}~\cite{b8}.

\begin{figure}[H]
\centering
\includegraphics[width=5in]{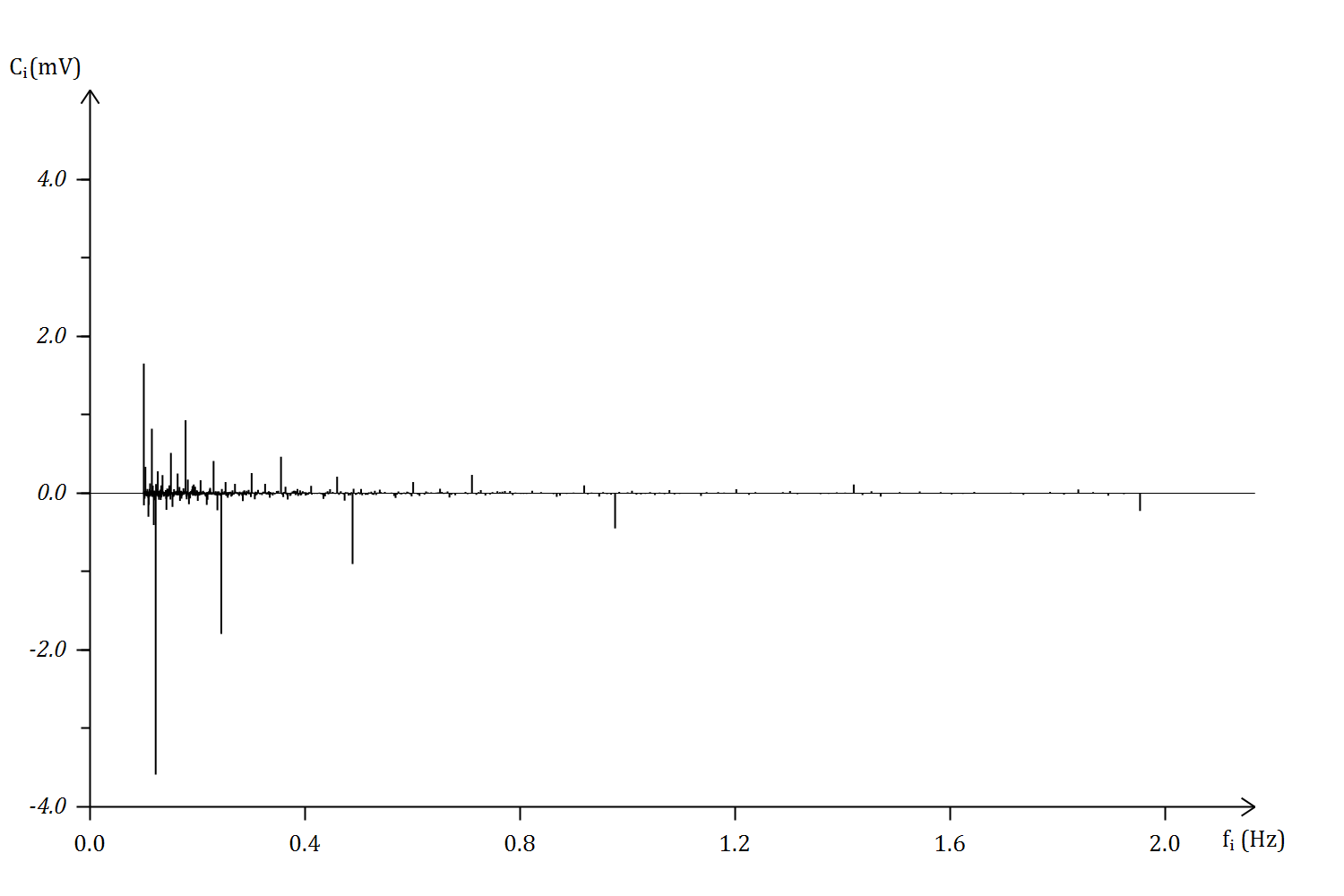}
\caption{Partial plot of the SWT corresponding to the electromyographic recording in figure~\ref{f9}. The moduli of the coefficients not represented (the $C_1$ for which $f_i > 2$) are much smaller than the moduli of $C_1$ which were represented.}
\label{f10}
\end{figure}

Extending only to $f_i=2 \:\mathrm{Hz}$, the graph in figure~\ref{f10} is a section of the SWT corresponding to the recording represented in figure~\ref{f9}.

The sequence of dyads for values of $f_i < 0.12$ is as follows:
\begin{align*}
(f_{1},C_{1})&=(0.1000000, 1.655)  &\qquad (f_{2},C_{2})&=(0.1000801, 0.002) \\
(f_{3},C_{3})&=(0.1001603, -0.152)  &\qquad (f_{4},C_{4})&=(0.1002406, -0.013) \\
(f_{5},C_{5})&=(0.1003210, 0.025)  &\qquad (f_{6},C_{6})&=(0.1004016, 0.004) \\
(f_{7},C_{7})&=(0.1004823, 0.017)  &\qquad (f_{8},C_{8})&=(0.1005632, -0.026) \\
(f_{9},C_{9})&=(0.1006441, 0.019)  &\qquad (f_{10},C_{10})&=(0.1007252, -0.027) \\
(f_{11},C_{11})&=(0.1008065, -0.033)  &\qquad (f_{12},C_{12})&=(0.1008878, -0.021) \\
(f_{13},C_{13})&=(0.1009693, 0.007)  &\qquad (f_{14},C_{14})&=(0.1010509, 0.011) \\
(f_{15},C_{15})&=(0.1011327, 0.023)  &\qquad (f_{16},C_{16})&=(0.1012146, 0.012) \\
(f_{17},C_{17})&=(0.1012966, 0.004)  &\qquad (f_{18},C_{18})&=(0.1013788, 0.006) \\
(f_{19},C_{19})&=(0.1014610, -0.067)  &\qquad (f_{20},C_{20})&=(0.1015435, 0.011) \\
(f_{21},C_{21})&=(0.1016260, -0.010)  &\qquad (f_{22},C_{22})&=(0.1017087, 0.011) \\
(f_{23},C_{23})&=(0.1017915, 0.016)  &\qquad (f_{24},C_{24})&=(0.1018745, -0.007) \\
(f_{25},C_{25})&=(0.1019576, 0.019)  &\qquad (f_{26},C_{26})&=(0.1020408, -0.005) \\
(f_{27},C_{27})&=(0.1021242, -0.009)  &\qquad (f_{28},C_{28})&=(0.1022077, 0.005) \\
(f_{29},C_{29})&=(0.1022913, -0.021)  &\qquad (f_{30},C_{30})&=(0.1023751, 0.028) \\
(f_{31},C_{31})&=(0.1024590, 0.010)  &\qquad (f_{32},C_{32})&=(0.1025431, -0.005) \\
(f_{33},C_{33})&=(0.1026273, -0.001)  &\qquad (f_{34},C_{34})&=(0.1027116, 0.010) \\
(f_{35},C_{35})&=(0.1027961, 0.340)  &\qquad (f_{36},C_{36})&=(0.1028807, 0.024) \\
(f_{37},C_{37})&=(0.1029654, 0.006)  &\qquad (f_{38},C_{38})&=(0.1030503, 0.016) \\
(f_{39},C_{39})&=(0.1031353, -0.001)  &\qquad (f_{40},C_{40})&=(0.1032205, 0.001) \\
(f_{41},C_{41})&=(0.1033058, 0.017)  &\qquad (f_{42},C_{42})&=(0.1033912, 0.006) \\
(f_{43},C_{43})&=(0.1034768, 0.006)  &\qquad (f_{44},C_{44})&=(0.1035626, -0.019) \\
(f_{45},C_{45})&=(0.1036484, 0.010)  &\qquad (f_{46},C_{46})&=(0.1037344, -0.015) \\
(f_{47},C_{47})&=(0.1038206, -0.031)  &\qquad (f_{48},C_{48})&=(0.1039069, -0.011) \\
(f_{49},C_{49})&=(0.1039933, 0.006)  &\qquad (f_{50},C_{50})&=(0.1040799, 0.015) \\
(f_{51},C_{51})&=(0.1041667, -0.020)  &\qquad (f_{52},C_{52})&=(0.1042535, -0.031) \\
(f_{53},C_{53})&=(0.1043406, 0.005)  &\qquad (f_{54},C_{54})&=(0.1044277, 0.011) \\
(f_{55},C_{55})&=(0.1045151, -0.010)  &\qquad (f_{56},C_{56})&=(0.1046025, -0.004) \\
(f_{57},C_{57})&=(0.1046901, -0.012)  &\qquad (f_{58},C_{58})&=(0.1047779, 0.010) \\
(f_{59},C_{59})&=(0.1048658, 0.027)  &\qquad (f_{60},C_{60})&=(0.1049538, -0.005) \\
(f_{61},C_{61})&=(0.1050420, 0.014)  &\qquad (f_{62},C_{62})&=(0.1051304, -0.026) \\
(f_{63},C_{63})&=(0.1052189, 0.020)  &\qquad (f_{64},C_{64})&=(0.1053075, 0.013) \\
(f_{65},C_{65})&=(0.1053963, 0.006)  &\qquad (f_{66},C_{66})&=(0.1054852, 0.010) \\
(f_{67},C_{67})&=(0.1055743, 0.018)  &\qquad (f_{68},C_{68})&=(0.1056636, -0.002) \\
(f_{69},C_{69})&=(0.1057530, 0.003)  &\qquad (f_{70},C_{70})&=(0.1058425, 0.011) \\
(f_{71},C_{71})&=(0.1059322, -0.034)  &\qquad (f_{72},C_{72})&=(0.1060221, -0.005) \\
(f_{73},C_{73})&=(0.1061121, -0.025)  &\qquad (f_{74},C_{74})&=(0.1062022, -0.039) \\
(f_{75},C_{75})&=(0.1062925, 0.055)  &\qquad (f_{76},C_{76})&=(0.1063830, -0.005) \\
(f_{77},C_{77})&=(0.1064736, 0.004)  &\qquad (f_{78},C_{78})&=(0.1065644, 0.034) \\
(f_{79},C_{79})&=(0.1066553, 0.018)  &\qquad (f_{80},C_{80})&=(0.1067464, 0.010) \\
(f_{81},C_{81})&=(0.1068376, -0.004)  &\qquad (f_{82},C_{82})&=(0.1069290, -0.011) \\
(f_{83},C_{83})&=(0.1070205, 0.021)  &\qquad (f_{84},C_{84})&=(0.1071123, -0.010) \\
(f_{85},C_{85})&=(0.1072041, -0.038)  &\qquad (f_{86},C_{86})&=(0.1072961, -0.009) \\
(f_{87},C_{87})&=(0.1073883, 0.002)  &\qquad (f_{88},C_{88})&=(0.1074807, 0.013) \\
(f_{89},C_{89})&=(0.1075731, 0.005)  &\qquad (f_{90},C_{90})&=(0.1076658, -0.009) \\
(f_{91},C_{91})&=(0.1077586, -0.012)  &\qquad (f_{92},C_{92})&=(0.1078516, 0.007) \\
(f_{93},C_{93})&=(0.1079447, 0.002)  &\qquad (f_{94},C_{94})&=(0.1080380, 0.003) \\
(f_{95},C_{95})&=(0.1081315, -0.058)  &\qquad (f_{96},C_{96})&=(0.1082251, -0.010) \\
(f_{97},C_{97})&=(0.1083189, 0.000)  &\qquad (f_{98},C_{98})&=(0.1084128, 0.013) \\
(f_{99},C_{99})&=(0.1085069, -0.300)  &\qquad (f_{100},C_{100})&=(0.1086012, 0.013) \\
(f_{101},C_{101})&=(0.1086957, 0.014)  &\qquad (f_{102},C_{102})&=(0.1087903, -0.013) \\
(f_{103},C_{103})&=(0.1088850, -0.037)  &\qquad (f_{104},C_{104})&=(0.1089799, -0.008) \\
(f_{105},C_{105})&=(0.1090750, 0.001)  &\qquad (f_{106},C_{106})&=(0.1091703, -0.005) \\
(f_{107},C_{107})&=(0.1092657, -0.146)  &\qquad (f_{108},C_{108})&=(0.1093613, 0.000) \\
(f_{109},C_{109})&=(0.1094571, 0.019)  &\qquad (f_{110},C_{110})&=(0.1095530, -0.022) \\
(f_{111},C_{111})&=(0.1096491, -0.014)  &\qquad (f_{112},C_{112})&=(0.1097454, -0.034) \\
(f_{113},C_{113})&=(0.1098418, 0.005)  &\qquad (f_{114},C_{114})&=(0.1099384, 0.003) \\
(f_{115},C_{115})&=(0.1100352, 0.021)  &\qquad (f_{116},C_{116})&=(0.1101322, 0.019) \\
(f_{117},C_{117})&=(0.1102293, -0.037)  &\qquad (f_{118},C_{118})&=(0.1103266, -0.023) \\
(f_{119},C_{119})&=(0.1104240, 0.056)  &\qquad (f_{120},C_{120})&=(0.1105217, 0.019) \\
(f_{121},C_{121})&=(0.1106195, 0.003)  &\qquad (f_{122},C_{122})&=(0.1107174, 0.008) \\
(f_{123},C_{123})&=(0.1108156, 0.007)  &\qquad (f_{124},C_{124})&=(0.1109139, -0.009) \\
(f_{125},C_{125})&=(0.1110124, -0.009)  &\qquad (f_{126},C_{126})&=(0.1111111, 0.002) \\
(f_{127},C_{127})&=(0.1112100, 0.019)  &\qquad (f_{128},C_{128})&=(0.1113090, 0.016) \\
(f_{129},C_{129})&=(0.1114082, 0.043)  &\qquad (f_{130},C_{130})&=(0.1115076, -0.008) \\
(f_{131},C_{131})&=(0.1116071, 0.128)  &\qquad (f_{132},C_{132})&=(0.1117069, -0.015) \\
(f_{133},C_{133})&=(0.1118068, -0.042)  &\qquad (f_{134},C_{134})&=(0.1119069, 0.019) \\
(f_{135},C_{135})&=(0.1120072, 0.000)  &\qquad (f_{136},C_{136})&=(0.1121076, -0.008) \\
(f_{137},C_{137})&=(0.1122083, 0.015)  &\qquad (f_{138},C_{138})&=(0.1123091, 0.008) \\
(f_{139},C_{139})&=(0.1124101, -0.010)  &\qquad (f_{140},C_{140})&=(0.1125113, -0.023) \\
(f_{141},C_{141})&=(0.1126126, -0.005)  &\qquad (f_{142},C_{142})&=(0.1127142, 0.009) \\
(f_{143},C_{143})&=(0.1128159, 0.035)  &\qquad (f_{144},C_{144})&=(0.1129178, -0.002) \\
(f_{145},C_{145})&=(0.1130199, 0.007)  &\qquad (f_{146},C_{146})&=(0.1131222, 0.032) \\
(f_{147},C_{147})&=(0.1132246, -0.024)  &\qquad (f_{148},C_{148})&=(0.1133273, -0.002) \\
(f_{149},C_{149})&=(0.1134301, -0.020)  &\qquad (f_{150},C_{150})&=(0.1135332, 0.000) \\
(f_{151},C_{151})&=(0.1136364, -0.008)  &\qquad (f_{152},C_{152})&=(0.1137398, 0.006) \\
(f_{153},C_{153})&=(0.1138434, 0.015)  &\qquad (f_{154},C_{154})&=(0.1139471, 0.011) \\
(f_{155},C_{155})&=(0.1140511, -0.014)  &\qquad (f_{156},C_{156})&=(0.1141553, 0.014) \\
(f_{157},C_{157})&=(0.1142596, 0.001)  &\qquad (f_{158},C_{158})&=(0.1143641, 0.021) \\
(f_{159},C_{159})&=(0.1144689, -0.002)  &\qquad (f_{160},C_{160})&=(0.1145738, 0.013) \\
(f_{161},C_{161})&=(0.1146789, -0.014)  &\qquad (f_{162},C_{162})&=(0.1147842, 0.005) \\
(f_{163},C_{163})&=(0.1148897, 0.824)  &\qquad (f_{164},C_{164})&=(0.1149954, 0.012) \\
(f_{165},C_{165})&=(0.1151013, -0.003)  &\qquad (f_{166},C_{166})&=(0.1152074, 0.005) \\
(f_{167},C_{167})&=(0.1153137, 0.010)  &\qquad (f_{168},C_{168})&=(0.1154201, 0.005) \\
(f_{169},C_{169})&=(0.1155268, 0.005)  &\qquad (f_{170},C_{170})&=(0.1156337, -0.005) \\
(f_{171},C_{171})&=(0.1157407, -0.044)  &\qquad (f_{172},C_{172})&=(0.1158480, -0.021) \\
(f_{173},C_{173})&=(0.1159555, -0.012)  &\qquad (f_{174},C_{174})&=(0.1160631, -0.008) \\
(f_{175},C_{175})&=(0.1161710, 0.005)  &\qquad (f_{176},C_{176})&=(0.1162791, -0.001) \\
(f_{177},C_{177})&=(0.1163873, 0.004)  &\qquad (f_{178},C_{178})&=(0.1164958, -0.018) \\
(f_{179},C_{179})&=(0.1166045, 0.097)  &\qquad (f_{180},C_{180})&=(0.1167134, -0.004) \\
(f_{181},C_{181})&=(0.1168224, -0.004)  &\qquad (f_{182},C_{182})&=(0.1169317, 0.014) \\
(f_{183},C_{183})&=(0.1170412, 0.017)  &\qquad (f_{184},C_{184})&=(0.1171509, -0.032) \\
(f_{185},C_{185})&=(0.1172608, -0.017)  &\qquad (f_{186},C_{186})&=(0.1173709, 0.003) \\
(f_{187},C_{187})&=(0.1174812, -0.013)  &\qquad (f_{188},C_{188})&=(0.1175917, 0.010) \\
(f_{189},C_{189})&=(0.1177024, 0.009)  &\qquad (f_{190},C_{190})&=(0.1178134, 0.008) \\
(f_{191},C_{191})&=(0.1179245, -0.027)  &\qquad (f_{192},C_{192})&=(0.1180359, -0.009) \\
(f_{193},C_{193})&=(0.1181474, 0.034)  &\qquad (f_{194},C_{194})&=(0.1182592, 0.000) \\
(f_{195},C_{195})&=(0.1183712, -0.404)  &\qquad (f_{196},C_{196})&=(0.1184834, -0.009) \\
(f_{197},C_{197})&=(0.1185958, -0.044)  &\qquad (f_{198},C_{198})&=(0.1187085, -0.023) \\
(f_{199},C_{199})&=(0.1188213, 0.032)  &\qquad (f_{200},C_{200})&=(0.1189343, 0.007) \\
(f_{201},C_{201})&=(0.1190476, 0.018)  &\qquad (f_{202},C_{202})&=(0.1191611, 0.011) \\
(f_{203},C_{203})&=(0.1192748, 0.025)  &\qquad (f_{204},C_{204})&=(0.1193887, -0.008) \\
(f_{205},C_{205})&=(0.1195029, 0.006)  &\qquad (f_{206},C_{206})&=(0.1196172, 0.032) \\
(f_{207},C_{207})&=(0.1197318, 0.024)  &\qquad (f_{208},C_{208})&=(0.1198466, -0.012) \\
(f_{209},C_{209})&=(0.1199616, 0.005)  &\qquad 
\end{align*}

In \cite{b9} potential users can find a tool which makes it possible to obtain the SWT of sequences of samples of different recordings, for $n \leq 3000$.

\section{Analysis of an Image}

The SWM is also useful for the analysis of images and the corresponding results can be expressed with the SWM.

\begin{figure}[H]
\centering
\begin{tikzpicture}[
axis/.style={thick, ->, >=stealth'}]
\draw[axis]( -0.1,0)  -- (4.4,0) node(xline)[below] {$x$};
\draw[axis] (0,-0.1) -- (0,4.4) node(yline)[left] {$y$};
\draw[step=1cm,color=gray] (0,0) grid (4,4);
\node at (+0.5,+0.5) {\begin{tabular}{c}$P_{1,1}$ \\ 100 \end{tabular} };
\node at (+1.5,+0.5) {\begin{tabular}{c}$P_{2,1}$ \\ 38 \end{tabular} };
\node at (+2.5,+0.5) {\begin{tabular}{c}$P_{3,1}$ \\ 25 \end{tabular} };
\node at (+3.5,+0.5) {\begin{tabular}{c}$P_{4,1}$ \\ 214 \end{tabular} };
\node at (+0.5,+1.5) {\begin{tabular}{c}$P_{1,2}$ \\ 98 \end{tabular} };
\node at (+1.5,+1.5) {\begin{tabular}{c}$P_{2,2}$ \\ 3 \end{tabular} };
\node at (+2.5,+1.5) {\begin{tabular}{c}$P_{3,2}$ \\ 77 \end{tabular} };
\node at (+3.5,+1.5) {\begin{tabular}{c}$P_{4,2}$ \\ 12 \end{tabular} };
\node at (+0.5,+2.5) {\begin{tabular}{c}$P_{1,3}$ \\ 195 \end{tabular} };
\node at (+1.5,+2.5) {\begin{tabular}{c}$P_{2,3}$ \\ 6 \end{tabular} };
\node at (+2.5,+2.5) {\begin{tabular}{c}$P_{3,3}$ \\ 249 \end{tabular} };
\node at (+3.5,+2.5) {\begin{tabular}{c}$P_{4,3}$ \\ 255 \end{tabular} };
\node at (+0.5,+3.5) {\begin{tabular}{c}$P_{1,4}$ \\ 55 \end{tabular} };
\node at (+1.5,+3.5) {\begin{tabular}{c}$P_{2,4}$ \\ 4 \end{tabular} };
\node at (+2.5,+3.5) {\begin{tabular}{c}$P_{3,4}$ \\ 69 \end{tabular} };
\node at (+3.5,+3.5) {\begin{tabular}{c}$P_{4,4}$ \\ 81 \end{tabular} };
\node at (+0.5,-0.2) {1};
\node at (+1.5,-0.2) {2};
\node at (+2.5,-0.2) {3};
\node at (+3.5,-0.2) {4};
\node at (-0.2,+0.5) {1};
\node at (-0.2,+1.5) {2};
\node at (-0.2,+2.5) {3};
\node at (-0.2,+3.5) {4};
\end{tikzpicture}
\caption{Representation of a pseudo-image. The name of each pixel is given at the top, and the level of gray corresponding to that pixel is given at the bottom.}
\label{f11}
\end{figure}
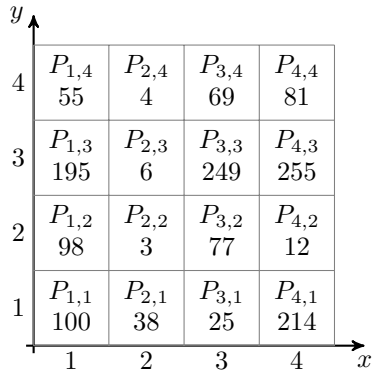

First, a ``pseudo-image'' (shown in figure~\ref{f11}) has been chosen to illustrate how the SWM may be used for analysis purposes. The method applied for the analysis of a pseudo-image is the same as that used for the analysis of genuine images. The term ``pseudo-image'' indicates that no reference is made to a genuine image of something real, such as an object or a living being, or abstract art with aesthetic value. Given that the method used, the SWM, for the analysis of a genuine image and of a pseudo-image is the same, the latter will serve to provide a clear and simple example of how that method can be applied.

The name of each pixel in the pseudo-image displayed in figure~\ref{f11} is given at the top. The numerical value (ranging from 0 for black to 255 for white) corresponding to the level of gray of the pixel is given at the bottom. (The levels of gray were chosen arbitrarily.)

Although figure~\ref{f11} can be considered as a representation of a matrix of pixels, the diverse elements in that matrix ($P_{i,j}$, where $i=1,2,3,4$; and $j=1,2,3,4$) are referred to differently from what is usually found. Here it has been taken into account that in the system of orthogonal Cartesian coordinates used (in which the $x$-axis is that of the abscissas, and the $y$-axis is that of the ordinates), the point $(i,j)$ (where $i=1,2,3,4$; and $j=1,2,3,4$) whose abscissa is $i$ and whose ordinate is $j$, is the center point of the pixel $P_{i,j}$.

In figure~\ref{f12}, it has been shown how each of the axes ($x$ and $y$) in figure 11 has been processed with the same approach as that of figure~\ref{f2} with the $x$-axis (the only axis of coordinates used in this case).

\begin{figure}[H]
\centering
\begin{tikzpicture}[
axis/.style={thick, ->, >=stealth'},dashed line/.style={dashed, thin}]
\draw[axis]( 0,0)  -- (4.4,0) node(xline)[below] {$x$};
\draw[axis] (0,0) -- (0,4.4) node(yline)[left] {$y$};
\draw[step=1cm,color=gray] (0,0) grid (4,4);
\draw[dashed] (-4,0) -- (0,0);
\draw[dashed] (-4,1) -- (0,1);
\draw[dashed] (-4,2) -- (0,2);
\draw[dashed] (-4,3) -- (0,3);
\draw[dashed] (-4,4) -- (0,4);
\draw[dashed] (0,0) -- (0,-4);
\draw[dashed] (1,0) -- (1,-4);
\draw[dashed] (2,0) -- (2,-4);
\draw[dashed] (3,0) -- (3,-4);
\draw[dashed] (4,0) -- (4,-4);
\node at (+0.5,+0.5) {\begin{tabular}{c}$P_{1,1}$ \\ 100 \end{tabular} };
\node at (+1.5,+0.5) {\begin{tabular}{c}$P_{2,1}$ \\ 38 \end{tabular} };
\node at (+2.5,+0.5) {\begin{tabular}{c}$P_{3,1}$ \\ 25 \end{tabular} };
\node at (+3.5,+0.5) {\begin{tabular}{c}$P_{4,1}$ \\ 214 \end{tabular} };
\node at (+0.5,+1.5) {\begin{tabular}{c}$P_{1,2}$ \\ 98 \end{tabular} };
\node at (+1.5,+1.5) {\begin{tabular}{c}$P_{2,2}$ \\ 3 \end{tabular} };
\node at (+2.5,+1.5) {\begin{tabular}{c}$P_{3,2}$ \\ 77 \end{tabular} };
\node at (+3.5,+1.5) {\begin{tabular}{c}$P_{4,2}$ \\ 12 \end{tabular} };
\node at (+0.5,+2.5) {\begin{tabular}{c}$P_{1,3}$ \\ 195 \end{tabular} };
\node at (+1.5,+2.5) {\begin{tabular}{c}$P_{2,3}$ \\ 6 \end{tabular} };
\node at (+2.5,+2.5) {\begin{tabular}{c}$P_{3,3}$ \\ 249 \end{tabular} };
\node at (+3.5,+2.5) {\begin{tabular}{c}$P_{4,3}$ \\ 255 \end{tabular} };
\node at (+0.5,+3.5) {\begin{tabular}{c}$P_{1,4}$ \\ 55 \end{tabular} };
\node at (+1.5,+3.5) {\begin{tabular}{c}$P_{2,4}$ \\ 4 \end{tabular} };
\node at (+2.5,+3.5) {\begin{tabular}{c}$P_{3,4}$ \\ 69 \end{tabular} };
\node at (+3.5,+3.5) {\begin{tabular}{c}$P_{4,4}$ \\ 81 \end{tabular} };
\node at (-0.5,+0.5) {$\tensor*[^2]{C}{_1}$};
\node at (-0.5,+1.5) {$\tensor*[^2]{C}{_1}$};
\node at (-0.5,+2.5) {$\tensor*[^2]{C}{_1}$};
\node at (-0.5,+3.5) {$\tensor*[^2]{C}{_1}$};
\node at (-1.5,+0.5) {$\tensor*[^2]{C}{_2}$};
\node at (-1.5,+1.5) {$\tensor*[^2]{C}{_2}$};
\node at (-1.5,+2.5) {$\tensor*[^2]{C}{_2}$};
\node at (-1.5,+3.5) {$-\tensor*[^2]{C}{_2}$};
\node at (-2.5,+0.5) {$\tensor*[^2]{C}{_3}$};
\node at (-2.5,+1.5) {$\tensor*[^2]{C}{_3}$};
\node at (-2.5,+2.5) {$-\tensor*[^2]{C}{_3}$};
\node at (-2.5,+3.5) {$-\tensor*[^2]{C}{_3}$};
\node at (-3.5,+0.5) {$\tensor*[^2]{C}{_4}$};
\node at (-3.5,+1.5) {$-\tensor*[^2]{C}{_4}$};
\node at (-3.5,+2.5) {$\tensor*[^2]{C}{_4}$};
\node at (-3.5,+3.5) {$-\tensor*[^2]{C}{_4}$};
\node at (+0.5,-0.5) {$\tensor*[^1]{C}{_1}$};
\node at (+1.5,-0.5) {$\tensor*[^1]{C}{_1}$};
\node at (+2.5,-0.5) {$\tensor*[^1]{C}{_1}$};
\node at (+3.5,-0.5) {$\tensor*[^1]{C}{_1}$};
\node at (+0.5,-1.5) {$\tensor*[^1]{C}{_2}$};
\node at (+1.5,-1.5) {$\tensor*[^1]{C}{_2}$};
\node at (+2.5,-1.5) {$\tensor*[^1]{C}{_2}$};
\node at (+3.5,-1.5) {$-\tensor*[^1]{C}{_2}$};
\node at (+0.5,-2.5) {$\tensor*[^1]{C}{_3}$};
\node at (+1.5,-2.5) {$\tensor*[^1]{C}{_3}$};
\node at (+2.5,-2.5) {$-\tensor*[^1]{C}{_3}$};
\node at (+3.5,-2.5) {$-\tensor*[^1]{C}{_3}$};
\node at (+0.5,-3.5) {$\tensor*[^1]{C}{_4}$};
\node at (+1.5,-3.5) {$-\tensor*[^1]{C}{_4}$};
\node at (+2.5,-3.5) {$\tensor*[^1]{C}{_4}$};
\node at (+3.5,-3.5) {$-\tensor*[^1]{C}{_4}$};
\Dimline[(0,-4)][(4,-4)][$\Delta x$][above]
\Dimline[(-4,4)][(-4,0)][$\Delta y$][above]
\end{tikzpicture}
\caption{How to apply the SWM to the analysis of the pseudo-image represented in figure 11. (See indications in text.) The superscript 1 located on the left of the coefficient indicates that the coefficient corresponds to the $x$-axis and the superscript 2 on the left of the coefficient indicates that it corresponds to the $y$-axis.}
\label{f12}
\end{figure}

Note, for example, how the equation corresponding to pixel $P_{4,3}$ is formed: that pixel belongs to both column 4 and row 3 in figure~\ref{f12}:

\begin{enumerate}
\item Take the sequence of these two sets: a) the set of elements underneath column 4 in figure 12; and b) the set of elements at the left of row 3 in the same figure:
\begin{equation*}
\left(\tensor*[^1]{C}{_1},-\tensor*[^1]{C}{_2},-\tensor*[^1]{C}{_3},-\tensor*[^1]{C}{_4}\right)\quad
\left(\tensor*[^2]{C}{_1},\tensor*[^2]{C}{_2},-\tensor*[^2]{C}{_3},\tensor*[^2]{C}{_4}\right)
\end{equation*}

Note that the above expression is an ordered pair of sets.

\item The Cartesian product is found for the two sets corresponding to that ordered pair of sets:
\begin{multline*}
\left(\tensor*[^1]{C}{_1},-\tensor*[^1]{C}{_2},-\tensor*[^1]{C}{_3},-\tensor*[^1]{C}{_4}\right)\times
\left(\tensor*[^2]{C}{_1},\tensor*[^2]{C}{_2},-\tensor*[^2]{C}{_3},\tensor*[^2]{C}{_4}\right)=P^{*}_{4,3}=\\
\{
(\enspace\tensor*[^1]{C}{_1},\;\tensor*[^2]{C}{_1}),\: (\;\tensor*[^1]{C}{_1},\;\tensor*[^2]{C}{_2}),\: (\;\tensor*[^1]{C}{_1},-\tensor*[^2]{C}{_3}),\: (\;\tensor*[^1]{C}{_1},\tensor*[^2]{C}{_4}),\\ (-\tensor*[^1]{C}{_2},\tensor*[^2]{C}{_1}),\: (-\tensor*[^1]{C}{_2},\tensor*[^2]{C}{_2}),\: (-\tensor*[^1]{C}{_2},-\tensor*[^2]{C}{_3}),\:(-\tensor*[^1]{C}{_2},\tensor*[^2]{C}{_4}),\\ (-\tensor*[^1]{C}{_3},\tensor*[^2]{C}{_1}),\: (-\tensor*[^1]{C}{_3},\tensor*[^2]{C}{_2}),\: (-\tensor*[^1]{C}{_3},-\tensor*[^2]{C}{_3}),\: (-\tensor*[^1]{C}{_3},\tensor*[^2]{C}{_4}),\\ (-\tensor*[^1]{C}{_4},\tensor*[^2]{C}{_1}),\: (-\tensor*[^1]{C}{_4},\tensor*[^2]{C}{_2}),\: (-\tensor*[^1]{C}{_4},-\tensor*[^2]{C}{_3}),\: (-\tensor*[^1]{C}{_4},\tensor*[^2]{C}{_4})\}
\end{multline*}

The preceding Cartesian product corresponding to pixel $P_{4,3}$ has been designated as $P^{*}_{4,3}$.

\item The set $P^{*}_{4,3}$ is acted on by an operator O such that:

\begin{enumerate}

\item it converts each element $P^{*}_{4,3}$ (a certain ordered pair) into a single element which will be a coefficient with two subscripts: the first is the same as that of the first element of the ordered pair, and the second subscript is the same as that of the second element of the ordered pair; and
\item if the signs preceding the two elements of the ordered pair have the same sign ($+$ and $+$; or $-$ and $-$), the coefficient is positive; if they are different ($+$ and $-$; or $-$ and $+$), the coefficient obtained is negative.

\end{enumerate}
Therefore, for the case considered, the following result is obtained:
\begin{align*}
\mathrm{O}(P^{*}_{4,3}) = & (C_{1,1}, C_{1,2},-C_{1,3},C_{1,4},\\
 & -C_{2,1},-C_{2,2},C_{2,3},-C_{2,4}, \\
 & -C_{3,1},-C_{3,2},C_{3,3},-C_{3,4}, \\
 & -C_{4,1},-C_{4,2},C_{4,3},-C_{4,4})
\end{align*}

The right-hand member of the above equation is a set of coefficients, each of which has two subscripts and is preceded by the $+$ sign or by the $-$ sign. That set will be called $C^{*}_{4,3}$. Of course, $C^{*}_{4,3}$ must not be confused with the coefficient $C_{4,3}$.

\item The algebraic sum is figured for all the elements in $C^{*}_{4,3}$, and the result is equated to the numeric value of the level of gray corresponding to pixel $P_{4,3}$ (i.e., $255$).

Therefore, the linear algebraic equation corresponding to pixel $P_{4,3}$ is obtained:
\begin{multline*}
C_{1,1}+C_{1,2}-C_{1,3}+C_{1,4}-C_{2,1}-C_{2,2}+C_{2,3}-C_{2,4} \\
-C_{3,1}-C_{3,2}+C_{3,3}-C_{3,4} -C_{4,1},-C_{4,2}+C_{4,3}-C_{4,4}=255
\end{multline*}
\end{enumerate}
If the same procedure is applied to the set of elements under column 3 in figure 12 and to the set of elements appearing at the left of row 1 in that same figure, the following linear algebraic equation may be specified for pixel $P_{3,1}$:
\begin{multline*}
C_{1,1}+C_{1,2}+C_{1,3}+C_{1,4}+C_{2,1}+C_{2,2}+C_{2,3}+C_{2,4} \\ 
-C_{3,1}-C_{3,2}-C_{3,3}-C_{3,4}+C_{4,1}+C_{4,2}+C_{4,3}+C_{4,4}= 25 
\end{multline*}

If the same procedure is carried out for each of the remaining pixels of the pseudo-image in figure 11, the following system of linear algebraic equations is obtained:

\begin{equation}
  \left.\begin{aligned}
C_{1,1}&+C_{1,2}+C_{1,3}+C_{1,4}+C_{2,1}+C_{2,2}+C_{2,3}+C_{2,4}& \\ & +C_{3,1}+C_{3,2}+C_{3,3}+C_{3,4}+C_{4,1}+C_{4,2}+C_{4,3}+C_{4,4}& = 100 \\
 C_{1,1}&+C_{1,2}+C_{1,3}-C_{1,4}+C_{2,1}+C_{2,2}+C_{2,3}-C_{2,4}& \\ & +C_{3,1}+C_{3,2}+C_{3,3}-C_{3,4}+C_{4,1}+C_{4,2}+C_{4,3}-C_{4,4}& = 98 \\
 C_{1,1}&+C_{1,2}-C_{1,3}+C_{1,4}+C_{2,1}+C_{2,2}-C_{2,3}+C_{2,4}& \\ & +C_{3,1}+C_{3,2}-C_{3,3}+C_{3,4}+C_{4,1}+C_{4,2}-C_{4,3}+C_{4,4}& = 195 \\
 C_{1,1}&-C_{1,2}-C_{1,3}-C_{1,4}+C_{2,1}-C_{2,2}-C_{2,3}-C_{2,4}& \\ & +C_{3,1}-C_{3,2}-C_{3,3}-C_{3,4}+C_{4,1}-C_{4,2}-C_{4,3}-C_{4,4}& = 55 \\
 C_{1,1}&+C_{1,2}+C_{1,3}+C_{1,4}+C_{2,1}+C_{2,2}+C_{2,3}+C_{2,4}& \\ & +C_{3,1}+C_{3,2}+C_{3,3}+C_{3,4}-C_{4,1}-C_{4,2}-C_{4,3}-C_{4,4}& = 38 \\
 C_{1,1}&+C_{1,2}+C_{1,3}-C_{1,4}+C_{2,1}+C_{2,2}+C_{2,3}-C_{2,4}& \\ & +C_{3,1}+C_{3,2}+C_{3,3}-C_{3,4}-C_{4,1}-C_{4,2}-C_{4,3}+C_{4,4}& = 3 \\
 C_{1,1}&+C_{1,2}-C_{1,3}+C_{1,4}+C_{2,1}+C_{2,2}-C_{2,3}+C_{2,4}& \\ & +C_{3,1}+C_{3,2}-C_{3,3}+C_{3,4}-C_{4,1}-C_{4,2}+C_{4,3}-C_{4,4}& = 6 \\
 C_{1,1}&-C_{1,2}-C_{1,3}-C_{1,4}+C_{2,1}-C_{2,2}-C_{2,3}-C_{2,4}& \\ & +C_{3,1}-C_{3,2}-C_{3,3}-C_{3,4}-C_{4,1}+C_{4,2}+C_{4,3}+C_{4,4}& = 4 \\
 C_{1,1}&+C_{1,2}+C_{1,3}+C_{1,4}+C_{2,1}+C_{2,2}+C_{2,3}+C_{2,4}& \\ & -C_{3,1}-C_{3,2}-C_{3,3}-C_{3,4}+C_{4,1}+C_{4,2}+C_{4,3}+C_{4,4}& = 25 \\
 C_{1,1}&+C_{1,2}+C_{1,3}-C_{1,4}+C_{2,1}+C_{2,2}+C_{2,3}-C_{2,4}& \\ & -C_{3,1}-C_{3,2}-C_{3,3}+C_{3,4}+C_{4,1}+C_{4,2}+C_{4,3}-C_{4,4}& = 77 \\
 C_{1,1}&+C_{1,2}-C_{1,3}+C_{1,4}+C_{2,1}+C_{2,2}-C_{2,3}+C_{2,4}& \\ & -C_{3,1}-C_{3,2}+C_{3,3}-C_{3,4}+C_{4,1}+C_{4,2}-C_{4,3}+C_{4,4}& = 249 \\
 C_{1,1}&-C_{1,2}-C_{1,3}-C_{1,4}+C_{2,1}-C_{2,2}-C_{2,3}-C_{2,4}& \\ & -C_{3,1}+C_{3,2}+C_{3,3}+C_{3,4}+C_{4,1}-C_{4,2}-C_{4,3}-C_{4,4}& = 69 \\
 C_{1,1}&+C_{1,2}+C_{1,3}+C_{1,4}-C_{2,1}-C_{2,2}-C_{2,3}-C_{2,4}& \\ & -C_{3,1}-C_{3,2}-C_{3,3}-C_{3,4}-C_{4,1}-C_{4,2}-C_{4,3}-C_{4,4}& = 55 \\
 C_{1,1}&+C_{1,2}+C_{1,3}-C_{1,4}-C_{2,1}-C_{2,2}-C_{2,3}+C_{2,4}& \\ & -C_{3,1}-C_{3,2}-C_{3,3}+C_{3,4}-C_{4,1}-C_{4,2}-C_{4,3}+C_{4,4}& = 12 \\
 C_{1,1}&+C_{1,2}-C_{1,3}+C_{1,4}-C_{2,1}-C_{2,2}+C_{2,3}-C_{2,4}& \\ & -C_{3,1}-C_{3,2}+C_{3,3}-C_{3,4}-C_{4,1}-C_{4,2}+C_{4,3}-C_{4,4}& = 255 \\
 C_{1,1}&-C_{1,2}-C_{1,3}-C_{1,4}-C_{2,1}+C_{2,2}+C_{2,3}+C_{2,4}& \\ & -C_{3,1}+C_{3,2}+C_{3,3}+C_{3,4}-C_{4,1}+C_{4,2}+C_{4,3}+C_{4,4}& = 81
   \end{aligned}
  \right\}  
\label{e6}
\end{equation}

The 16 unknowns of the preceding system of linear algebraic equations are the 16 coefficients $C_{i,j}$, where $i=1,2,3,4$; and $j=1,2,3,4$. 

If the system of equations~\eqref{e6} is solved, the following values are found for the 16 unknowns:
\begin{align*}
C_{1,1}&=112.50000 & C_{3,1}&=-34.00000 \\
C_{1,2}&=-78.50000 & C_{3,2}&=-14.00000 \\
C_{1,3}&=15.25000  & C_{3,3}&=32.25000 \\
C_{1,4}&=28.25000  & C_{3,4}&=-31.75000 \\
C_{2,1}&=27.50000  & C_{4,1}&=51.00000  \\
C_{2,2}&=22.25000  & C_{4,2}&=-55.25000 \\
C_{2,3}&=-23.50000 & C_{4,3}&=13.50000 \\
C_{2,4}&=42.75000  & C_{4,4}&=-8.25000 
\end{align*}

Information will be provided later about the analysis of the image of Lenna mentioned above, regarding the meaning of this type of coefficients. For now, it will suffice to say that the system of equations (6) obtained for the pseudo-image considered ``corresponds'' to the system of equations (3) obtained by analyzing the function of the single variable specified in (1). 

In general, to determine the equation corresponding to any pixel $P_{i,j}$ (where $i=1,2,3,\ldots,n$; and $j=1,2,3,\ldots,n$) of an image composed of $n$ rows and $n$ columns of pixels (i.e., with $n^2$ pixels), the following procedure may be used:

\begin{enumerate}
\item Using the same criterion as that of the case of the pseudo-image discussed above, it is possible to determine the set of elements underneath column $i$, and the set of elements at the left of row $j$ of the image analyzed. Recall that the abscissa $i$ of the center of pixel $P_{i,j}$ indicates the column where that pixel is located; and $j$, the ordinate of the center of $P_{i,j}$ indicates the row in which that pixel is found.
\item The Cartesian product is found for these two sets, and the result obtained is referred to as $P_{i,j}^{*}$.
\item An action is carried out on the set $P^{*}_{i,j}$ by an operator O such that:

\begin{enumerate}
\item it converts each element of $P^{*}_{i,j}$ (a particular ordered pair) into a unique element which will be a coefficient  with two subscripts: the first is the same as the subscript of the first element of that ordered pair and the second subscript is the same as the subscript of the second element of that ordered pair; and
\item if the signs preceding the two elements of the ordered pair are the same (both positive or both negative), then the sign of the coefficient obtained is positive. If, on the other hand, the signs are different (positive and negative, or negative and positive), then the sign of the coefficients is negative.
\end{enumerate}
\item $C^{*}_{i,j}$ is used to refer to the set obtained by having the operator O act on the set $P_{i,j}^{*}$:
\begin{equation*}
C^{*}_{i,j}=\mathrm{O}(P^{*}_{i,j})
\end{equation*}
\item The algebraic sum of all the elements is of $C^{*}_{i,j}$ is calculated and the result is equated to the numeric value of the level of gray corresponding to pixel $P_{i,j}$.
\end{enumerate}

When doing the same with each pixel $P_{i,j}$ (where $i=1,2,3,\ldots,n$; and $j=1,2,3,\ldots,n$), the result is a system of $n^2$ linear algebraic equations with $n^2$ unknowns. Those unknowns, which \textit{can} be found, are the $n^2$ coefficients $C_{i,j}$ (where $i=1,2,3,\ldots,n$; and $j=1,2,3,\ldots,n$).

The classic image of Lenna to be analyzed using the SWM is displayed in figure 13.

\begin{figure}[H]
\centering
\begin{tikzpicture}[
axis/.style={thick, ->, >=stealth'}]
\draw[axis]( -0.1,0)  -- (8.4,0) node(xline)[below] {$x$};
\draw[axis] (0,-0.1) -- (0,8.4) node(yline)[left] {$y$};
\node (img) 	at (4,4) {\includegraphics[width=8cm]{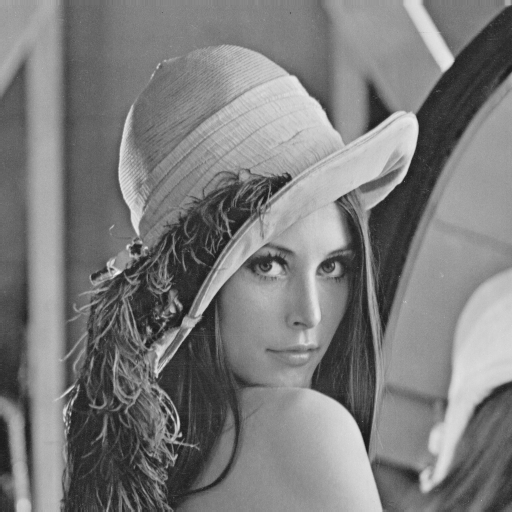}};
\end{tikzpicture}
\caption{Image of Lenna}
\label{f13}
\end{figure}

The image in figure~\ref{f13} is composed of $262,144$ pixels. It can be considered a matrix of pixels with 512 rows, each made up of 512 pixels, and 512 columns, each of which is also made up of 512 pixels. Each of the $262,144$ pixels belongs to only one of the 512 rows and only one of the 512 columns. In other words, the image of Lenna analyzed can be viewed as a matrix arranged in a $512\times512$ square grid of pixels.

If one applies the SWM directly to that image of Lenna, a system of $262,144$ linear algebraic equations is obtained. Because the authors of this paper cannot solve this system of linear equations due to the limitations of the computational tools available, the image has been further divided into 256 sub-images, each of which is a $32\times32$ matrix, consisting of:   
\begin{enumerate}
\item 16 rows of sub-images, each composed of 16 sub-images; or
\item 16 columns of sub-images, each composed of 16 sub-images.
\end{enumerate}

In other words, the image of Lenna analyzed can also be considered as a matrix arranged in a $16\times16$ square grid of sub-images.

The computational tools available to the authors \textit{do} have the capacity to apply the SWM to the analysis of each of these 256 sub-images, and this process was carried out successfully.  

The columns of sub-images have been numbered from left to right, from the first column (1) to the last column (16). The rows of sub-images have also been numbered from bottom up, row 1 to row 16. The sub-image belonging to column $k$ ($k=1,2,3,\ldots,16$) and row $l$ ($l=1,2,3,\ldots,16$) is called $I_{k,l}$.

In figure 14, it can be seen how the image of Lenna has been divided into 256 sub-images, so that each sub-image belongs only to one column and one row of the matrix.

\begin{figure}[H]
\centering
\begin{tikzpicture}[
axis/.style={thick, ->, >=stealth'}]
\draw[axis]( -0.1,0)  -- (8.4,0) node(xline)[below] {$x$};
\draw[axis] (0,-0.1) -- (0,8.4) node(yline)[left] {$y$};
\node (img) 	at (4,4) {\includegraphics[width=8cm]{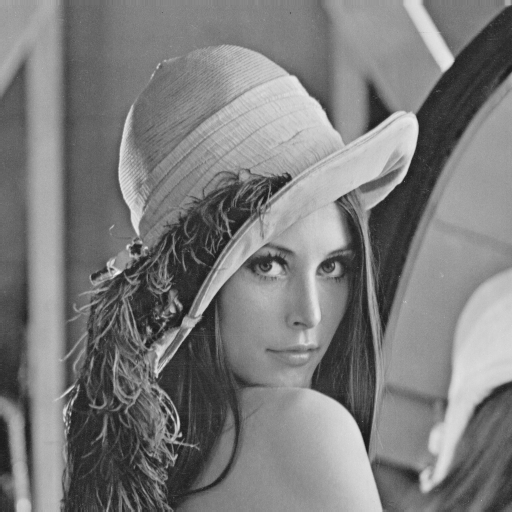}};
\draw[step=0.5cm,color=black] (0,0) grid (8,8);
\foreach \x [count=\xi] in {0.25,0.75,...,8} \node at(\x,-0.2) {$\scriptstyle\xi$};
\foreach \x [count=\xi] in {0.25,0.75,...,8} \node at(-0.2,\x) {$\scriptstyle\xi$};
\end{tikzpicture}
\caption{The image of Lenna is divided here into 256 sub-images.}
\label{f14}
\end{figure}

The notations introduced above for the pixels and for the sub-images can be used to refer to any pixel of the image of Lenna, by using the letter ``P'' for pixel and a tetrad of subscripts of that letter. Therefore, $P_{4,17,8,14}$ refers to the pixel located in the fourth column and the seventeenth row of the sub-image located, in turn, in the eighth column and the fourteenth row of the matrix of sub-images of Lenna.
 
The analysis of sub-image $I_{9,8}$ of the image of Lenna (i.e., the sub-image located in column 9, row 8 of the matrix of the sub-images of the image analyzed) is discussed below. The same approach was used for the analysis of all 255 remaining sub-images of the image of Lenna. Sub-image $I_{9,8}$ of the image of Lenna is displayed in figure 15.

\begin{figure}[H]
\centering
\begin{tikzpicture}[
axis/.style={thick, ->, >=stealth'}]
\draw[axis]( -0.1,0)  -- (4.4,0) node(xline)[below] {$x$};
\draw[axis] (0,-0.1) -- (0,4.4) node(yline)[left] {$y$};
\node (img) 	at (2,2) {\includegraphics[width=4cm]{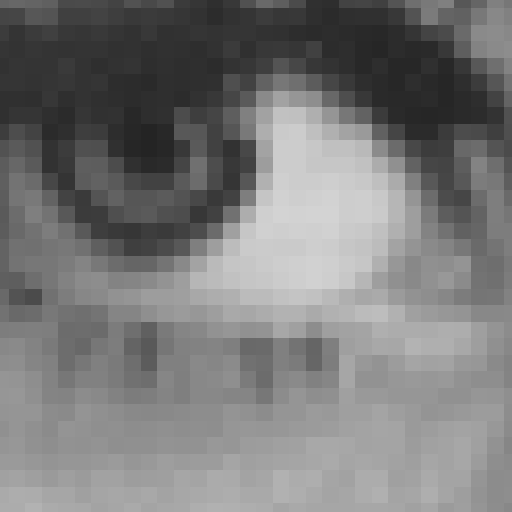}};
\end{tikzpicture}
\caption{Sub-image $I_{9,8}$ of the image of Lenna}
\label{f15}
\end{figure}

When using the SWM to analyze $I_{9,8}$, a system of 1024 linear algebraic equations was solved. Recall that each sub-image was composed of a matrix of $32\times32$ pixels. For each of those 1024 pixels, a linear algebraic equation was obtained using the procedure described.

The left-hand member of the linear algebraic equation obtained for each of the 1024 pixels of $I_{9,8}$ is an algebraic sum of the 1024 coefficients $C_{i,j}$ (where $i=1,2,3,\ldots,32$; and $j=1,2,3,\ldots,32$); and the right-hand member of the equation is the numeric value of the level of gray corresponding to the pixel considered. When solving this system of 1024 equations, each of the 1024 coefficients (the unknowns of the system of equations) can be positive, negative or null. Hence, according to the algebraic equation corresponding to each pixel, each coefficient can make a positive, negative or null contribution to that algebraic sum. Of course, to determine the contribution of each coefficient to that algebraic sum, both the computed value of each coefficient and the preceding sign in the algebraic sum must be taken into account. Thus, for example, when solving the system of equations, if it is determined that the value of a coefficient is $-250$, and the sign preceding it is negative, the contribution of that coefficient to that algebraic sum is equal to a positive numeric value: $250$. However, if the sign preceding that coefficient is positive, then the contribution of that coefficient to the algebraic sum is negative: $-250$.

Figure 16 represents the type of contribution (either positive or negative) which the coefficient $C_{17,25} = 14.500$ gives to the value of gray of each pixel of $I_{9,8}$, according to its numeric value (which then remains unchanged for all of the pixels of $I_{9,8}$) and to the preceding sign, which \textit{can} change according to the algebraic equation corresponding to each pixel of $I_{9,8}$. Blue is used for the positive contributions of $C_{17,25}$, and red for the negative contributions for each pixel of $I_{9,8}$. In figure 16, consideration was not given to the absolute value of these contributions to the numeric value of the level of gray of each pixel; it is only shown whether they are positive or negative.
\begin{figure}[H]
\centering
\begin{tikzpicture}
\draw[thick, ->]( -0.1,0)  -- (8.4,0) node(xline)[below] {$x$};
\draw[thick, ->] (0,-0.1) -- (0,8.4) node(yline)[left] {$y$};
\node (img) 	at (4,4) {\includegraphics[width=8cm]{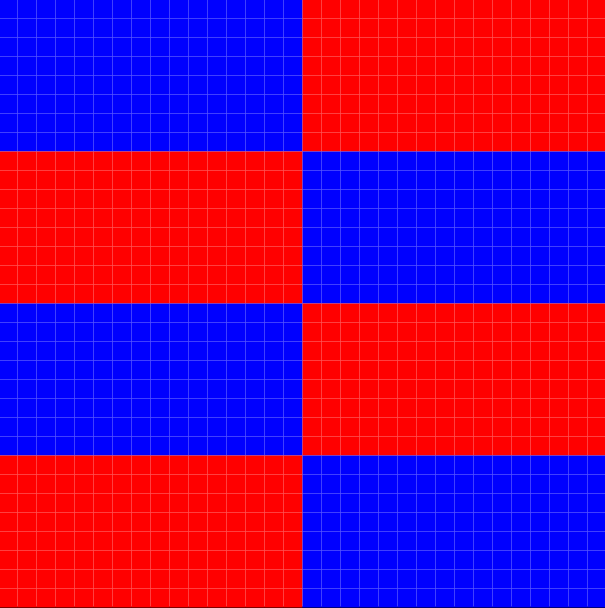}};
\end{tikzpicture}
\caption{Type of contribution (positive or negative) of $C_{17,25}$ to each pixel $I_{9,8}$. The positive contributions are represented in blue; the negative, in red.}
\label{f16}
\end{figure}

\begin{figure}[H]
\centering
\begin{tikzpicture}
\draw[thick, ->]( -0.1,0)  -- (8.4,0) node(xline)[below] {$x$};
\draw[thick, ->] (0,-0.1) -- (0,8.4) node(yline)[left] {$y$};
\node (img) 	at (4,4) {\includegraphics[width=8cm]{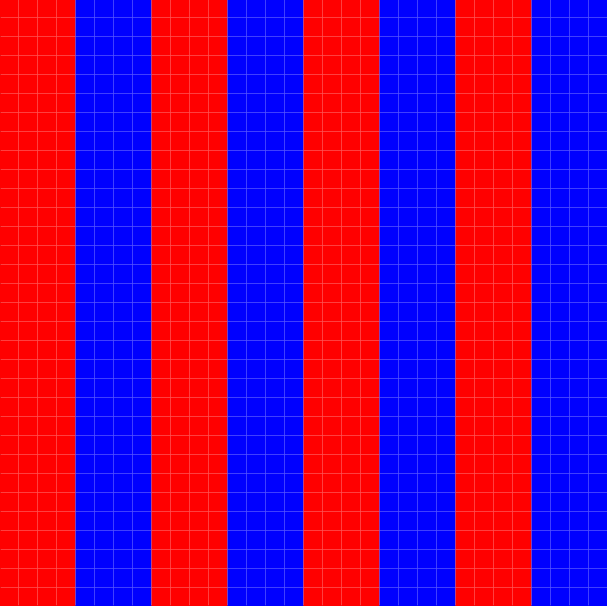}};
\end{tikzpicture}
\caption{Type of contribution (positive or negative) of $C_{29,1}$ to each pixel of $I_{9,8}$. The positive contributions are represented in blue; and the negative, in red.}
\label{f17}
\end{figure}

Like figure~\ref{f16}, figure~\ref{f17} indicates whether the coefficient $C_{29,1}$ provides a positive or negative contribution to each pixel $I_{9,8}$.

The contributions of each coefficient $C_{i,j}$ (where $i=1,2,3,\ldots,32$; and $j=1,2,3,\ldots,32$) to each pixel of $I_{9,8}$ have certain bidimensional patterns. If the modulus of each contribution is taken into account in addition to its positive or negative sign, these bidimensional patterns are the elements ``corresponding'' to the trains of square waves $S_1,S_2,\ldots, S_n$. 

An image of the type analyzed can be considered as a function of two variables: $x$ and $y$. In this case, to obtain the different approximations to the image analyzed, one must add, for each pixel in every sub-image, the patterns mentioned, taking into account the modulus of the coefficient corresponding to each pattern (just as in the cases of the functions of one variable, the trains of square waves are added to obtain the approximations to the functions analyzed).

In the case of a one-variable function, each coefficient $C_i$ (where $i=1,2,3,\ldots,n$) is associated to a specific frequency $f_i$. Given the way in which the coefficients $C_{i,j}$ (where $i=1,2,3,\ldots,n$; and $j=1,2,3,\ldots,n$) have been characterized, it is clear that each of the latter coefficients is similarly associated to the two spatial frequencies $f_{i,x}$ and $f_{j,y}$, corresponding to the $x$-axis and $y$-axis respectively. The first of these frequencies $f_{i,x}$ characterizes the periodicity made evident by the pattern corresponding to each coefficient $C_{i,j}$, according to the $x$-axis, and the second of these frequencies $f_{j,y}$ characterizes the periodicity shown by the same pattern, according to the $y$-axis. 

The computations of $f_{i,x}$ and $f_{j,y}$ are carried out by using the following equations (7) and (8), obtained from equation (2), substituting $\Delta t$ by $\Delta x$ and $\Delta y$, respectively (given that $f_{i,x}$ and $f_{j,y}$ are spatial frequencies) and equating $n=32$:

\begin{equation}
f_{i,x}=\frac{1}{2\Delta x} \left(\frac{32}{32-(i-1)}\right); \quad i=1,2,3,\ldots, 32
\label{e7}
\end{equation}

\begin{equation}
f_{j,y}=\frac{1}{2\Delta y} \left(\frac{32}{32-(j-1)}\right); \quad j=1,2,3,\ldots, 32
\label{e8}
\end{equation}

In this case, rather than taking the length of the side of a pixel as the unit of length, $\Delta x$ (which has the same length as $\Delta y$) can be used for that purpose. Both the length of $\Delta x$ and that of $\Delta y$ are equal to the product of $32$ times the length of one side of a pixel. If that is done, the following equations are obtained for $f_{i,x}$ and $f_{j,y}$:
\begin{equation}
f_{i,x}=\frac{1}{2} \left(\frac{32}{32-(i-1)}\right); \quad i=1,2,3,\ldots, 32
\label{e9}
\end{equation}

\begin{equation}
f_{j,y}=\frac{1}{2} \left(\frac{32}{32-(j-1)}\right); \quad j=1,2,3,\ldots, 32
\label{e10}
\end{equation}

In the first place, equations~\eqref{e9} and \eqref{e10} are used to consider the case presented in figure~\ref{f16}, which represents the contribution of $C_{17,25}$ (either positive or negative) to each pixel of $I_{9,8}$. Given that in this case $i=17$ and $j=25$, the following equations are obtained:

\begin{equation*}
f_{i,x}=\frac{1}{2} \left(\frac{32}{32-(17-1)}\right) = \frac{1}{2}  \left(\frac{32}{32-16)}\right) =\frac{1}{2}\left(\frac{32}{16}\right)=1
\end{equation*}

\begin{equation*}
f_{j,y}=\frac{1}{2} \left(\frac{32}{32-(25-1)}\right) = \frac{1}{2}  \left(\frac{32}{32-24)}\right) =\frac{1}{2}\left(\frac{32}{8}\right)=\frac{4}{2}=2
\end{equation*}

The meaning of the first of the preceding equations is as follows: In interval $\Delta x$ there is one wave on the $x$-axis, of the pattern represented in figure~\ref{f16}. (See figure to confirm.)

The meaning of the second of the preceding equations is as follows: In interval $\Delta y$ there are two waves on the $y$-axis, of the pattern represented in figure~\ref{f16}. (See figure to confirm.)

Secondly, equations~\eqref{e9} and \eqref{e10} will be used to consider the case in figure~\ref{f17}, which represents the contribution of $C_{29,1}$ (indicating either positive or negative) to each pixel $I_{9,8}$. Given that in this case $i=29$ and $j=1$, the following equations are obtained:

\begin{equation*}
f_{i,x}=\frac{1}{2} \left(\frac{32}{32-(29-1)}\right) = \frac{1}{2}\left(\frac{32}{32-28)}\right) =\frac{1}{2}\left(\frac{32}{4}\right)=\frac{8}{2}=4
\end{equation*}

\begin{equation*}
f_{j,y}=\frac{1}{2} \left(\frac{32}{32-(1-1)}\right) = \frac{1}{2}  \left(\frac{32}{32-0)}\right) =\frac{1}{2}\left(\frac{32}{32}\right)=\frac{1}{2}
\end{equation*}

The meaning of the first of the preceding equations is as follows: In interval $\Delta x$, there are 4 waves on the $x$-axis, of the pattern represented in figure 17. (See figure to confirm.)

The meaning of the second of the preceding equations is as follows: In the interval $\Delta y$, there is a half-wave on the $y$-axis, of the pattern represented in figure 17. (See figure to confirm.)

In the case of the images, the SWM generates two spatial frequencies, $f_{i,x}$ and $f_{j,y}$, associated with $C_{i,j}$, where $i=1,2,3,\ldots,n$, and $j=1,2,3,\ldots,n$. Hence, when applying the method it is also possible in the frequency domain to represent the results obtained. This expression in the frequency domain of the results obtained when analyzing images with SWM will be known as SWT, just as when functions of one variable were analyzed. The SWT will be used here to express some of the results produced.

If one adds, for each pixel of every sub-image of Lenna, the positive (or negative) contributions corresponding to the different coefficients, one obtains once again the image analyzed of Lenna. The analysis process carried out with the SWM, as indicated above, has made it possible to determine, unambiguously, for each pixel of the sub-image of Lenna, which patterns must be added in order to obtain the image of Lenna again. (Those patterns, as shown above, ``correspond'' to the trains of square waves obtained when analyzing functions with one variable with the SWM.)

Suppose that not all but only part of the patterns considered are added. For instance, let us admit that for each pixel of every sub-image of Lenna, the contributions corresponding to the coefficients  $C_{i,j}$ such that $i\leq8$ and $j\leq8$ are the only ones added. The image thus obtained will be called ``approximation $8\_32$" to the image of Lenna. On the other hand, if for each pixel of every sub-image of Lenna, the contributions corresponding to the coefficients $C_{i,j}$ such that $i\leq25$ and $j\leq25$ are the only ones added, the image thus obtained will be referred to as ``approximation $25\_32$" to the image of Lenna. Of course, if for each pixel of every sub-image of Lenna, the contributions corresponding to all the coefficients $C_{i,j}$ are added, the ``approximation $32\_32$", which is the same as the image itself, will be produced.

Several approximations of the image of Lenna have been displayed in \subref{f18a} to \subref{f18e}.

\newpage
\begin{figure}[H]
\centering
\subfloat[: Approximation 8\_32 to image of Lenna]{
\begin{tikzpicture}
\draw[thick, ->]( -0.1,0)  -- (8.4,0) node(xline)[below] {$x$};
\draw[thick, ->] (0,-0.1) -- (0,8.4) node(yline)[left] {$y$};
\node (img) 	at (4,4) {\includegraphics[width=8cm]{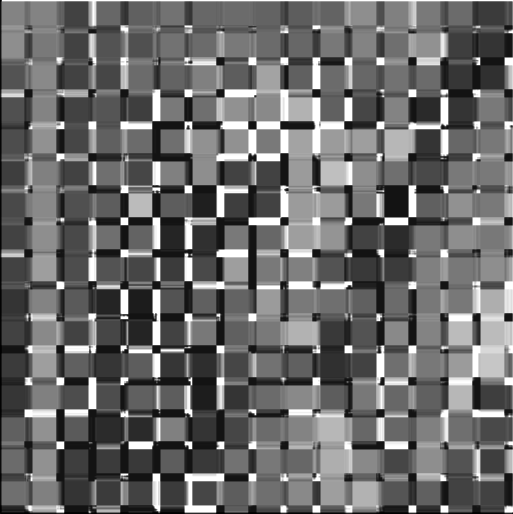}};
\end{tikzpicture}\label{f18a}}\\
\subfloat[: Approximation 16\_32 to image of Lenna]{
\begin{tikzpicture}
\draw[thick, ->]( -0.1,0)  -- (8.4,0) node(xline)[below] {$x$};
\draw[thick, ->] (0,-0.1) -- (0,8.4) node(yline)[left] {$y$};
\node (img) 	at (4,4) {\includegraphics[width=8cm]{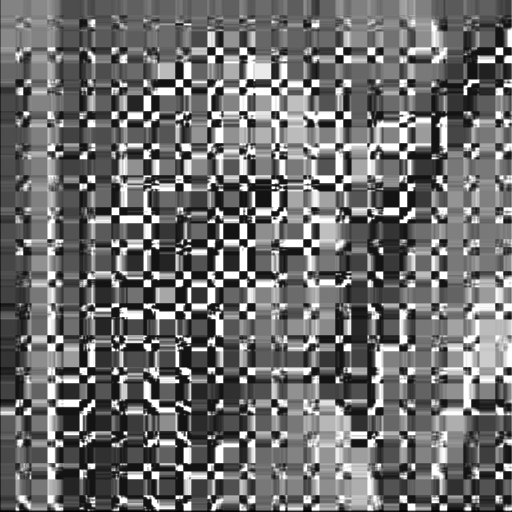}};
\end{tikzpicture}\label{f18b}}
\caption[]{}
\end{figure}

\begin{figure}[H]
\ContinuedFloat
\centering
\subfloat[: Approximation 24\_32 to image of Lenna]{\begin{tikzpicture}
\draw[thick, ->]( -0.1,0)  -- (8.4,0) node(xline)[below] {$x$};
\draw[thick, ->] (0,-0.1) -- (0,8.4) node(yline)[left] {$y$};
\node (img) 	at (4,4) {\includegraphics[width=8cm]{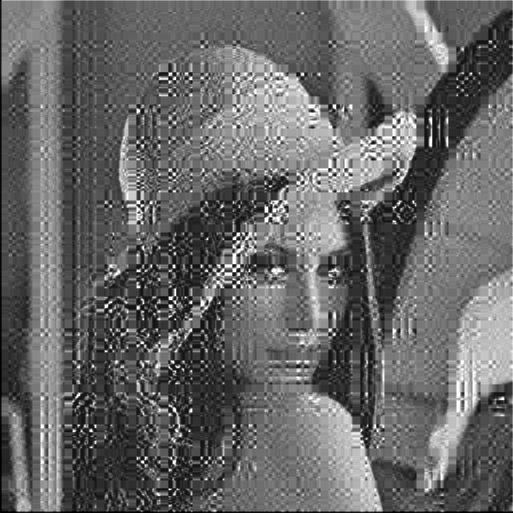}};
\end{tikzpicture}\label{f18c}}\\
\subfloat[: Approximation 28\_32 to image of Lenna]{\begin{tikzpicture}
\draw[thick, ->]( -0.1,0)  -- (8.4,0) node(xline)[below] {$x$};
\draw[thick, ->] (0,-0.1) -- (0,8.4) node(yline)[left] {$y$};
\node (img) 	at (4,4) {\includegraphics[width=8cm]{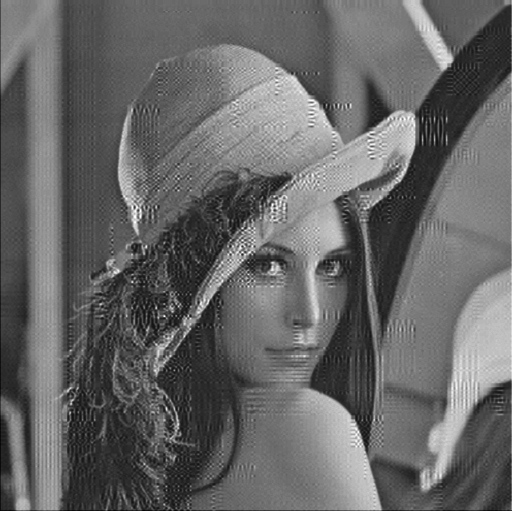}};
\end{tikzpicture}\label{f18d}}
\caption[]{}\end{figure}
\begin{figure}[H]
\ContinuedFloat
\centering
\subfloat[: Approximation 32\_32 to image of Lenna]{\begin{tikzpicture}
\draw[thick, ->]( -0.1,0)  -- (8.4,0) node(xline)[below] {$x$};
\draw[thick, ->] (0,-0.1) -- (0,8.4) node(yline)[left] {$y$};
\node (img) 	at (4,4) {\includegraphics[width=8cm]{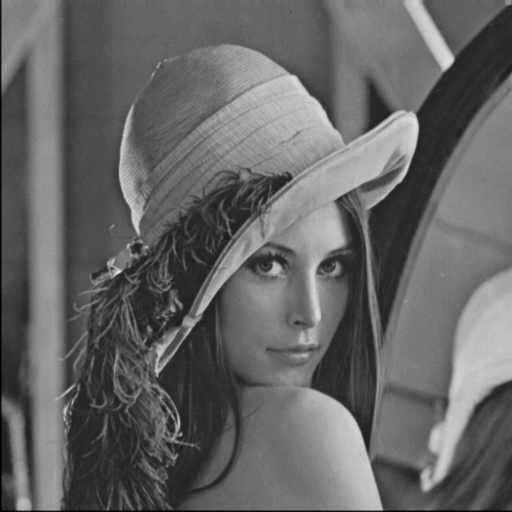}};
\end{tikzpicture}\label{f18e}}
\caption[]{}\end{figure}

In the examples presented below the results obtained by applying the SWM are given in the frequency domain. 

Take, for instance, the approximation $8\textunderscore 32$ to sub-image $I_{9,8}$ of the image of Lenna. The corresponding SWT can be expressed with a sequence of 64 triads. The third element of each triad is the coefficient considered. The first element of that triad is the frequency $f_{i,x}$ (where $i=1,2,3,\ldots,8$) on the $x$-axis, specified by the first subscript of that coefficient. The second element of that triad is the frequency $f_{j,y}$ (where $j=1,2,3,\ldots,8$) on the $y$-axis, specified by the second subscript of that coefficient. The complete list of these 64 triads is as follows:

\begin{align*}
(f_{1,x},f_{1,y};C_{1,1})&= (0.5000000, 0.5000000; 161.250)\\
(f_{2,x},f_{1,y};C_{2,1})&= (0.5161290, 0.5000000; 14.500) \\
(f_{3,x},f_{1,y};C_{3,1})&= (0.5333333, 0.5000000; -64.000) \\
(f_{4,x},f_{1,y};C_{4,1})&= (0.5517241, 0.5000000; 46.000) \\
(f_{5,x},f_{1,y};C_{5,1})&= (0.5714286, 0.5000000; 52.750) \\
(f_{6,x},f_{1,y};C_{6,1})&= (0.5925926, 0.5000000; -70.250) \\
(f_{7,x},f_{1,y};C_{7,1})&= (0.6153846, 0.5000000; 9.500) \\
(f_{8,x},f_{1,y};C_{8,1})&= (0.6400000, 0.5000000; -11.500) \\
(f_{1,x},f_{2,y};C_{1,2})&= (0.5000000, 0.5161290; -10.000) \\
(f_{2,x},f_{2,y};C_{2,2})&= (0.5161290, 0.5161290; -2.000) \\
(f_{3,x},f_{2,y};C_{3,2})&= (0.5333333, 0.5161290; 3.250) \\
(f_{4,x},f_{2,y};C_{4,2})&= (0.5517241, 0.5161290; -1.750) \\
(f_{5,x},f_{2,y};C_{5,2})&= (0.5714286, 0.5161290; -2.000) \\
(f_{6,x},f_{2,y};C_{6,2})&= (0.5925926, 0.5161290; 12.000) \\
(f_{7,x},f_{2,y};C_{7,2})&= (0.6153846, 0.5161290; -4.750) \\
(f_{8,x},f_{2,y};C_{8,2})&= (0.6400000, 0.5161290; -3.250) \\
(f_{1,x},f_{2,y};C_{1,3})&= (0.5000000, 0.5333333; -9.500) \\
(f_{2,x},f_{3,y};C_{2,3})&= (0.5161290, 0.5333333; 7.000) \\
(f_{3,x},f_{3,y};C_{3,3})&= (0.5333333, 0.5333333; -18.500) \\
(f_{4,x},f_{3,y};C_{4,3})&= (0.5517241, 0.5333333; 2.000) \\
(f_{5,x},f_{3,y};C_{5,3})&= (0.5714286, 0.5333333; 3.250) \\
(f_{6,x},f_{3,y};C_{6,3})&= (0.5925926, 0.5333333; -4.250) \\
(f_{7,x},f_{3,y};C_{7,3})&= (0.6153846, 0.5333333; 14.750) \\
(f_{8,x},f_{3,y};C_{8,3})&= (0.6400000, 0.5333333; 6.750) \\
(f_{1,x},f_{4,y};C_{1,4})&= (0.5000000, 0.5517241; -8.750) \\
(f_{2,x},f_{4,y};C_{2,4})&= (0.5161290, 0.5517241; -0.500) \\
(f_{3,x},f_{4,y};C_{3,4})&= (0.5333333, 0.5517241; 4.750) \\
(f_{4,x},f_{4,y};C_{4,4})&= (0.5517241, 0.5517241; -4.250) \\
(f_{5,x},f_{4,y};C_{5,4})&= (0.5714286, 0.5517241; -11.250) \\
(f_{6,x},f_{4,y};C_{6,4})&= (0.5925926, 0.5517241; 10.000) \\
(f_{7,x},f_{4,y};C_{7,4})&= (0.6153846, 0.5517241; -1.500) \\
(f_{8,x},f_{4,y};C_{8,4})&= (0.6400000, 0.5517241; 2.000) \\
(f_{1,x},f_{5,y};C_{1,5})&= (0.5000000, 0.5714286; -4.750) \\
(f_{2,x},f_{5,y};C_{2,5})&= (0.5161290, 0.5714286; -6.000) \\
(f_{3,x},f_{5,y};C_{3,5})&= (0.5333333, 0.5714286; -4.500) \\
(f_{4,x},f_{5,y};C_{4,5})&= (0.5517241, 0.5714286; -13.000) \\
(f_{5,x},f_{5,y};C_{5,5})&= (0.5714286, 0.5714286; -19.000) \\
(f_{6,x},f_{5,y};C_{6,5})&= (0.5925926, 0.5714286; 20.000) \\
(f_{7,x},f_{5,y};C_{7,5})&= (0.6153846, 0.5714286; -0.500) \\
(f_{8,x},f_{5,y};C_{8,5})&= (0.6400000, 0.5714286; -3.000) \\
(f_{1,x},f_{6,y};C_{1,6})&= (0.5000000, 0.5925926; -14.250) \\
(f_{2,x},f_{6,y};C_{2,6})&= (0.5161290, 0.5925926; 5.250) \\
(f_{3,x},f_{6,y};C_{3,6})&= (0.5333333, 0.5925926; 4.500) \\
(f_{4,x},f_{6,y};C_{4,6})&= (0.5517241, 0.5925926; -4.000) \\
(f_{5,x},f_{6,y};C_{5,6})&= (0.5714286, 0.5925926; 7.000) \\
(f_{6,x},f_{6,y};C_{6,6})&= (0.5925926, 0.5925926; 3.500) \\
(f_{7,x},f_{6,y};C_{7,6})&= (0.6153846, 0.5925926; 5.250) \\
(f_{8,x},f_{6,y};C_{8,6})&= (0.6400000, 0.5925926; 3.750) \\
(f_{1,x},f_{7,y};C_{1,7})&= (0.5000000, 0.6153846; -2.250) \\
(f_{2,x},f_{7,y};C_{2,7})&= (0.5161290, 0.6153846; 0.000) \\
(f_{3,x},f_{7,y};C_{3,7})&= (0.5333333, 0.6153846; 28.750) \\
(f_{4,x},f_{7,y};C_{4,7})&= (0.5517241, 0.6153846; -13.750) \\
(f_{5,x},f_{7,y};C_{5,7})&= (0.5714286, 0.6153846; -8.750) \\
(f_{6,x},f_{7,y};C_{6,7})&= (0.5925926, 0.6153846; 3.500) \\
(f_{7,x},f_{7,y};C_{7,7})&= (0.6153846, 0.6153846; 5.000) \\
(f_{8,x},f_{7,y};C_{8,7})&= (0.6400000, 0.6153846; 12.500) \\
(f_{1,x},f_{8,y};C_{1,8})&= (0.5000000, 0.6400000; 2.750) \\
(f_{2,x},f_{8,y};C_{2,8})&= (0.5161290, 0.6400000; -5.250) \\
(f_{3,x},f_{8,y};C_{3,8})&= (0.5333333, 0.6400000; 9.750) \\
(f_{4,x},f_{8,y};C_{4,8})&= (0.5517241, 0.6400000; 3.500) \\
(f_{5,x},f_{8,y};C_{5,8})&= (0.5714286, 0.6400000; 2.250) \\
(f_{6,x},f_{8,y};C_{6,8})&= (0.5925926, 0.6400000; 2.250) \\
(f_{7,x},f_{8,y};C_{7,8})&= (0.6153846, 0.6400000; -1.750) \\
(f_{8,x},f_{8,y};C_{8,8})&= (0.6400000, 0.6400000; -2.250) 
\end{align*}

The SWTs of different approximations to sub-image I9, 8 of Lenna have been displayed graphically in figure 19. Sub-figure 19.1 is the graphic presentation corresponding to the above sequence of 64 triads. (Some coefficients cannot be seen clearly here due to the small size of the values of their moduli.)

\newpage
\begin{figure}[H]
\resetsubfigs
\centering
\subfloat[: SWT of the approximation 8\_32 to sub-image $I_{9,8}$ of the image of Lenna. The positive coefficients are shown in blue, and the negative in red.]{\includegraphics[width=4in]{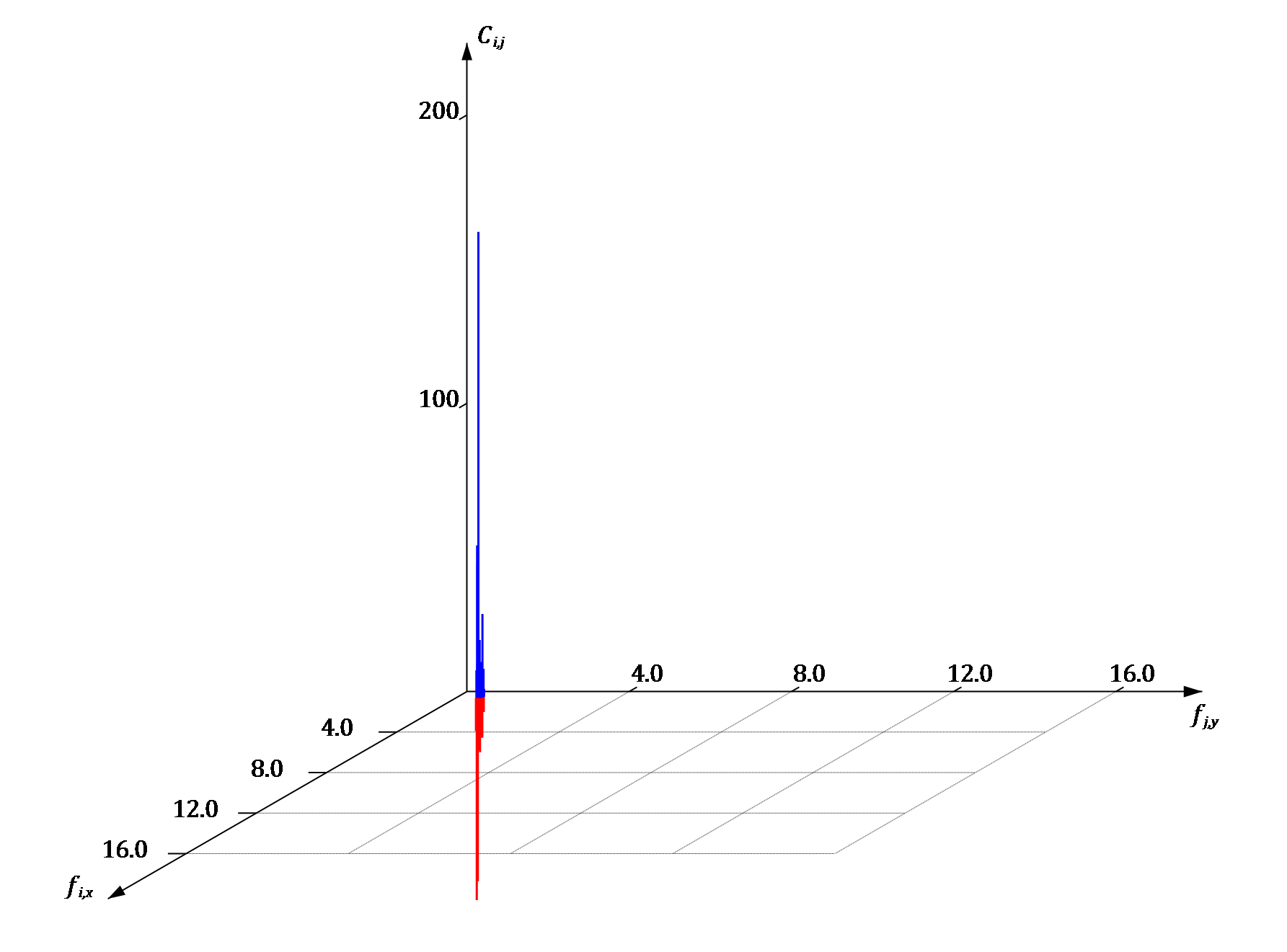}\label{f19a}}\qquad
\subfloat[: SWT of the approximation 16\_32 to sub-image $I_{9,8}$ of the image of Lenna. The positive coefficients are shown in blue, and the negative in red.]{\includegraphics[width=4in]{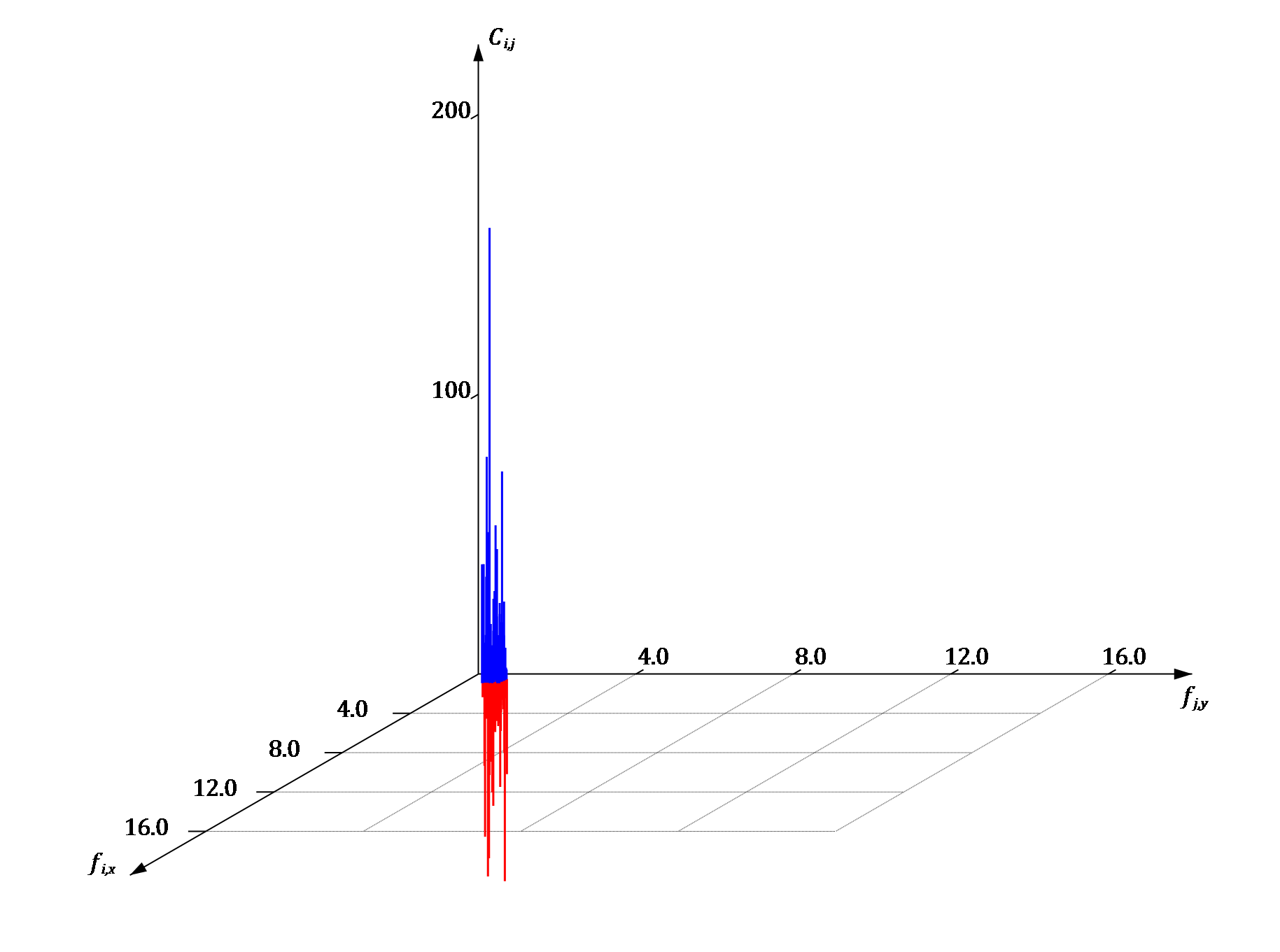}\label{f19b}}\qquad
\caption[]{Graphic display of SWTs of different approximations to sub-image I9, 8 of Lenna}\end{figure}
\begin{figure}[H]
\ContinuedFloat
\centering
\subfloat[: SWT of the approximation 24\_32 to sub-image $I_{9,8}$ of the image of Lenna. The positive coefficients are shown in blue, and the negative in red.]{\includegraphics[width=4in]{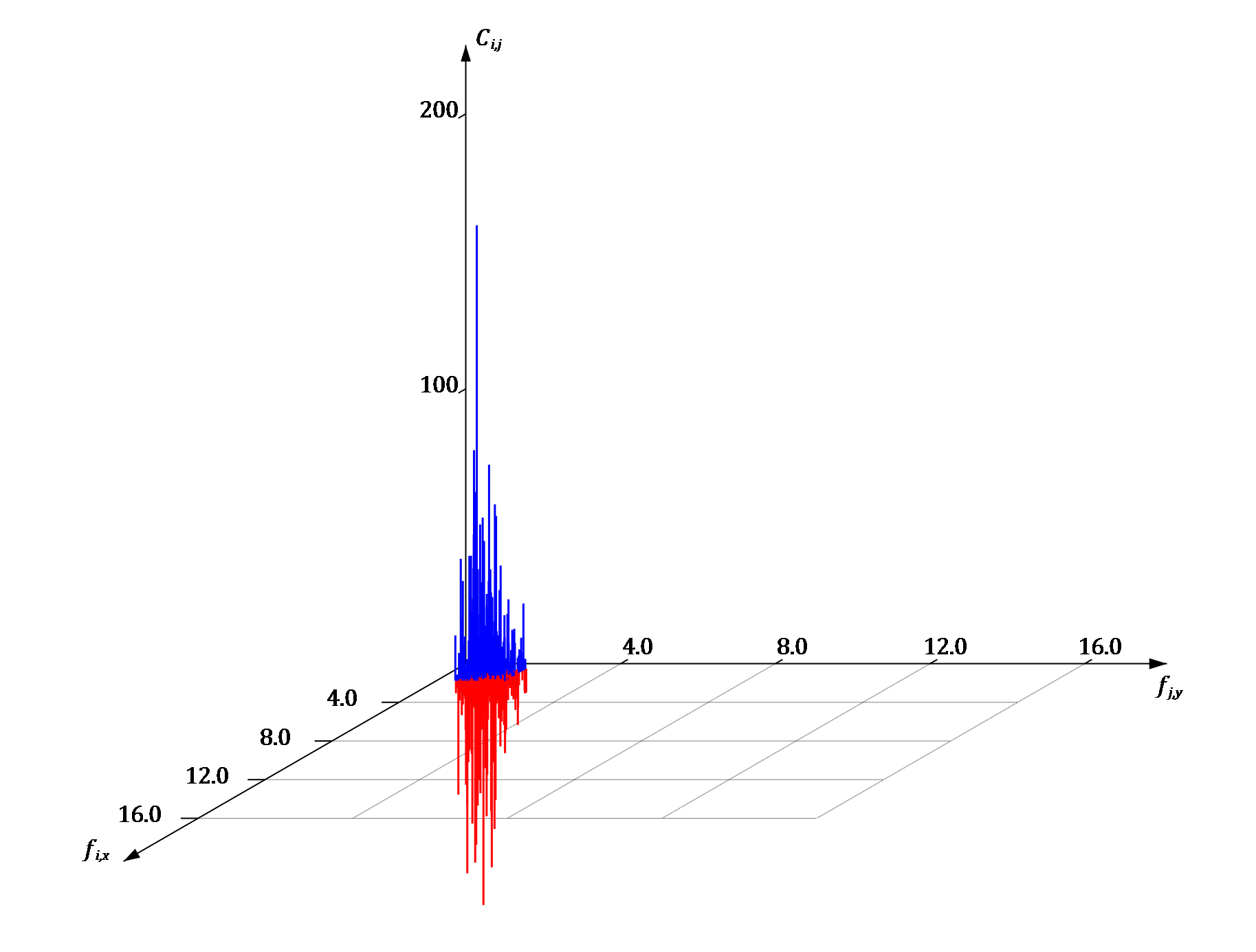}\label{f19c}}\qquad
\subfloat[: SWT of the approximation 32\_32 to sub-image $I_{9,8}$ of the image of Lenna. The positive coefficients are shown in blue, and the negative in red.]{\includegraphics[width=4in]{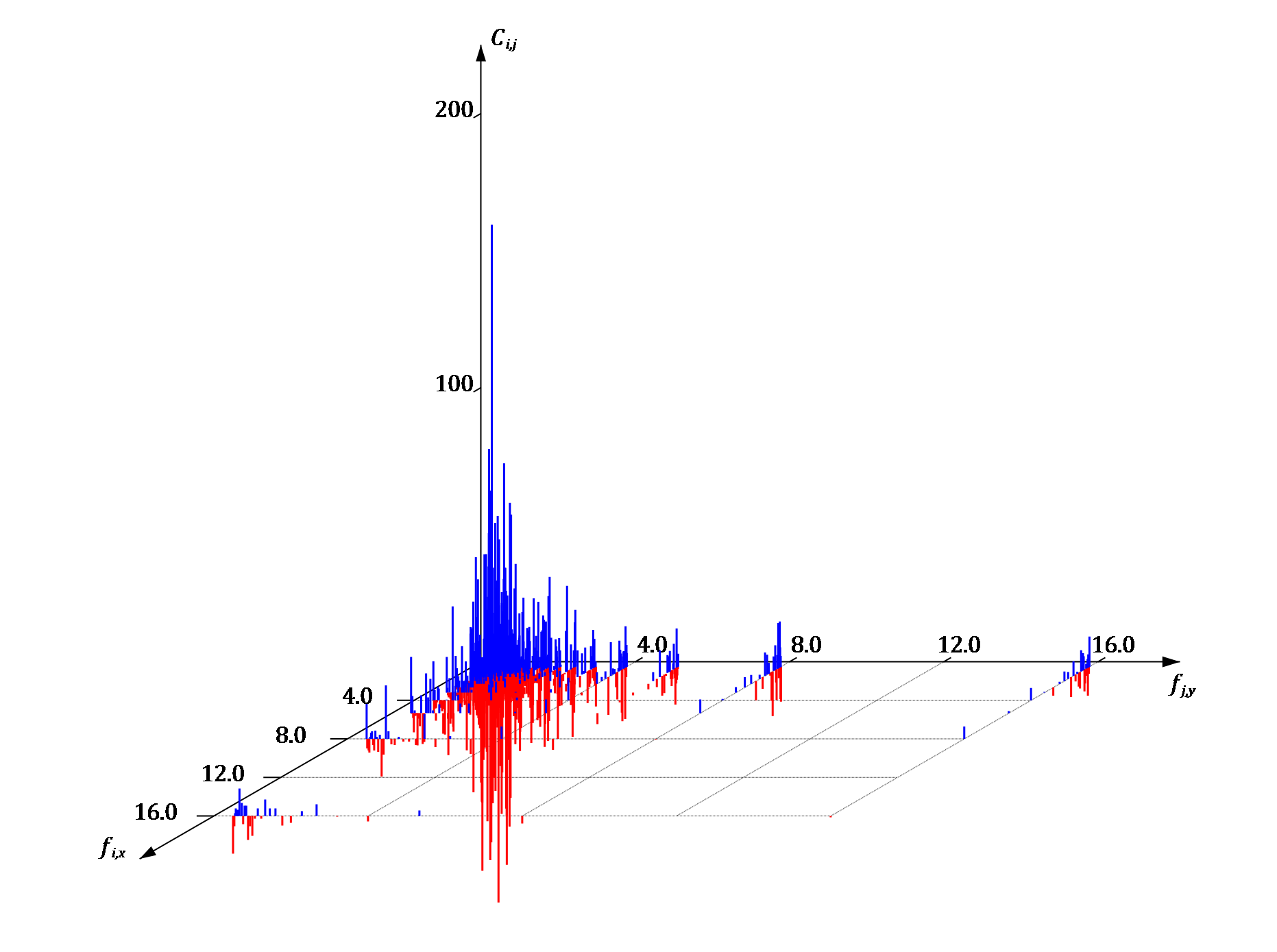}\label{f19d}}\qquad
\caption[]{Graphic display of SWTs of different approximations to sub-image I9, 8 of Lenna}\end{figure}

\section{Discussion and Prospects}

Not only has it been shown how the SWM can be applied to the analysis of 1) functions of one variable which in the intervals where they are analyzed satisfy the conditions of Dirichlet; and 2) sequences of samples from different types of biomedical recordings (ECGs, EEGs, and EMGs); but it has also been shown how that method can be applied to the analysis of images.

Fourier's outstanding and influential contributions, along with other approaches based on the development of his ideas, have played a very important role in a number of scientific and technological fields. In particular, mathematical tools such as the Fourier series, the Fourier transform, the discrete Fourier transform (DFT), the fast Fourier transform (FFT) -- an algorithm to compute the DFT -- and wavelets are used in the field of signal and image analysis. Perhaps in this field, the SWM can become a useful tool to complement the techniques provided by the Fourier approach. 

Among the advantages which the SWM offers are its simplicity and efficiency with which it is possible to make the necessary computations for its application, using software and hardware readily available; and 2) the systematic way in which it is applied: In all cases its use implies the solution of a system of linear algebraic equations which can be specified unambiguously.

An algorithm making it possible to apply the SWM to functions of $n$ variables, where $n=1,2,3,\dots$, will be addressed in another article.

The authors of this paper are particularly interested in the application of the SWM to the analysis of signals and images essential to the field of medicine and expect to devote several articles to go into greater depth on the topic.


\end{document}